\pdfoutput=1
\RequirePackage{ifpdf}
\ifpdf % We are running pdfTeX in pdf mode
\documentclass[pdftex]{sigma}
\else
\documentclass{sigma}
\fi

\usepackage{nccmath}
\usepackage[mathcal]{eucal}

\usepackage{tikz}\usetikzlibrary{backgrounds,matrix,arrows,shapes,decorations.markings}
\newcommand{\btp}{\begin{tikzpicture}[baseline=-5pt,scale=0.25,line width=0.7pt]}
\newcommand{\etp}{\end{tikzpicture}}

\numberwithin{equation}{section}

\newtheorem{Theorem}{Theorem}[section]
\newtheorem{Lemma}[Theorem]{Lemma}
\newtheorem{Proposition}[Theorem]{Proposition}
\theoremstyle{definition}
\newtheorem{Definition}[Theorem]{Definition}
\newtheorem{Example}[Theorem]{Example}
\newtheorem{Remark}[Theorem]{Remark}

\def\cA{{\mathcal A}} \def\cB{{\mathcal B}} 
 
 \def\cH{{\mathcal H}} 
\def\cJ{{\mathcal J}} \def\cK{{\mathcal K}} \def\cL{{\mathcal L}}
\def\cM{{\mathcal M}}  \def\cO{{\mathcal O}}
  
  \def\cU{{\mathcal U}}
  
\def\cY{{\mathcal Y}}  

\newcommand{\nc}{\newcommand}
\newcommand{\rnc}{\renewcommand}
\nc{\ul}{\underline}
\nc{\nn}{\nonumber}
\nc{\goth}{\mathfrak}
\rnc{\bold}{\mathbf}
\renewcommand{\frak}{\mathfrak}

\nc{\J}{{\mathcal J}}

\nc{\N}{{\Bbb N}}
\nc\lan{\langle}
\nc\ran{\rangle}
\nc\sub{\subset}
\nc\ti{\tilde}
\nc\wt{\widetilde}
\nc\nl{\newline}
\nc\ot{\otimes}
\nc\op{\oplus}
\nc\ad{\text{\rm ad}}
\nc\al{\alpha}
\nc\bet{\beta}
\nc\ga{\gamma}
\nc\de{\delta}
\nc\ep{\epsilon}
\nc\io{\iota}
\nc\om{\omega}
\nc\si{\sigma}
\rnc\th{\theta}
\nc\ka{\kappa}
\nc\la{\lambda}
\nc\ze{\zeta}
\nc\Ga{\Gamma}
\nc\De{\Delta}
\nc\Om{\Omega}
\nc\Si{\Sigma}
\nc\Th{\Theta}
\nc\La{\Lambda}
\nc\fa{\frak a}
\nc\fb{\frak b}
\nc\fc{\frak c}
\nc\fd{\frak d}
\nc\fe{\frak e}
\nc\ff{\frak f}
\nc\fg{\frak g}
\nc\fh{\frak h}
\nc\fj{\frak j}
\nc\fk{\frak k}
\nc\fl{\frak l}
\nc\fm{\frak m}
\nc\fn{\frak n}
\nc\fo{\frak o}
\nc\fp{\frak p}
\nc\fq{\frak q}
\nc\fr{\frak r}
\nc\fs{\frak s}
\nc\ft{\frak t}
\nc\fu{\frak u}
\nc\fv{\frak v}
\nc\fz{\frak z}
\nc\fx{\frak x}
\nc\fy{\frak y}
\nc\qu{\quad}
\nc\qq{\qquad}
\nc{\hs}[1]{\hspace{#1 mm}}
\nc{\mb}[1]{\hs{4}\mbox{#1}\hs{4}}
\nc\RR{\mathbb R}
\nc\CC{\mathbb C}
\nc\NN{\mathbb N}
\nc\ZZ{\mathbb Z}
\nc\II{\mathbb I}
\nc\cfg{{\bf c}_{\mathfrak{g}}}
\nc\cfh{{\bf c}_{\mathfrak{h}}}
\nc\cfm{{\bf c}_{\mathfrak{m}}}
\nc\mysum{\textstyle\sum\limits}
\nc{\fsl}{\mathfrak{sl}}
\nc{\fso}{\mathfrak{so}}
\nc{\fsp}{\mathfrak{sp}}
\nc{\fgl}{\mathfrak{gl}}
\nc{\sk}{\;\;\,}
\nc{\To}{\qq\Longrightarrow\qq}
\nc\B{\mathcal{B}}
\nc\G{\mathcal{G}}
\nc{\el}{\nonumber\\}
\nc{\Cx}{\mathsf x}
\nc{\Cy}{\mathsf y}
\nc{\Ce}{\mathsf e}
\nc{\Cf}{\mathsf f}
\nc{\Ch}{\mathsf h}
\nc{\Ck}{\mathsf k}
\nc{\CE}{\mathsf E}
\nc{\CF}{\mathsf F}
\nc{\CH}{\mathsf H}
\nc{\CY}{\mathsf Y}
\nc{\CX}{\mathsf X}
\nc{\Dh}{\Delta_\hbar}
\nc{\LH}{\!\!\mbox{
\btp
\draw[thick] (0,-.65) -- (.8,-.65) -- (.8,.3);
\draw[thin] (0,-.65) -- (.8,.3);
\etp
}\hspace{.25mm}}

\nc{\LHb}{\!\!\mbox{
\btp
\draw[thick] (0,-.65) -- (.8,-.65) -- (.8,.3);
\draw[thin] (0,-.65) -- (.8,.3);
\etp
}\hspace{.25mm}_\hbar}

\nc{\RHb}{\!\!\mbox{
\btp
\draw[thick] (0,.3) -- (0,-.65) -- (.8,-.65);
\draw[thin] (.8,-.65) -- (0,.3);
\etp
}\hspace{.25mm}_\hbar}

\nc{\LD}{\!\!\mbox{
\btp
\draw[thick] (0,-.65) -- (.8,-.65) -- (.8,.3);
\draw[thin] (0,-.65) -- (.8,.3);
\etp
}}

\nc{\RD}{\!\!\mbox{
\btp
\draw[thick] (0,.3) -- (0,-.65) -- (.8,-.65);
\draw[thin] (.8,-.65) -- (0,.3);
\etp
}}

\begin{document}

\allowdisplaybreaks

\newcommand{\arXivNumber}{1401.2143}

\renewcommand{\thefootnote}{}

\renewcommand{\PaperNumber}{011}

\FirstPageHeading

\ShortArticleName{Drinfeld J Presentation of Twisted Yangians}

\ArticleName{Drinfeld J Presentation of Twisted Yangians\footnote{This paper is a~contribution to the Special Issue on Recent Advances in Quantum Integrable Systems. The full collection is available at \href{http://www.emis.de/journals/SIGMA/RAQIS2016.html}{http://www.emis.de/journals/SIGMA/RAQIS2016.html}}}

\Author{Samuel BELLIARD~$^\dag$ and Vidas REGELSKIS~$^{\ddag\S}$}

\AuthorNameForHeading{S.~Belliard and V.~Regelskis}

\Address{$^\dag$~Institut de Physique Th\'eorique, Orme des Merisiers batiment 774, CEA/DSM/IPhT, \\
\hphantom{$^\dag$}~CEA/Saclay, F-91191 Gif-sur-Yvette Cedex, France}
\EmailD{\href{mailto:samuel.belliard@cea.fr}{samuel.belliard@cea.fr}}

\Address{$^\ddag$~Department of Mathematics, University of York, Heslington, York, YO10 5DD, UK}
\EmailD{\href{mailto:vidas.regelskis@york.ac.uk}{vidas.regelskis@york.ac.uk}}

\Address{$^\S$~Department of Mathematics, University of Surrey, Guildford, GU2 7XH, UK}

\ArticleDates{Received May 24, 2016, in final form February 21, 2017; Published online March 01, 2017}

\Abstract{We present a quantization of a Lie coideal structure for twisted half-loop algebras of finite-dimensional simple complex Lie algebras. We obtain algebra closure relations of twisted Yangians in Drinfeld~J presentation for all symmetric pairs of simple Lie algebras and for simple twisted even half-loop Lie algebras. We provide the explicit form of the closure relations for twisted Yangians in Drinfeld~J presentation for the ${\mathfrak{sl}}_3$ Lie algebra.}


\Classification{81R10; 81R50; 17B37}

\renewcommand{\thefootnote}{\arabic{footnote}}
\setcounter{footnote}{0}

\section{Introduction} \label{Sec:1}

The Yangian $\cY(\fg)$ is a flat quantization of the half-loop Lie algebra $\cL^+\cong\fg[u]$ of a finite-dimensional simple complex Lie algebra $\fg$ \cite{Dri85}. The name Yangian is due to V.G.~Drinfeld to honour C.N.~Yang who found the simplest solution of the Yang--Baxter equation, the rational $R$ matrix \cite{Yan67} (see also \cite{Bax72,Bax82}). This $R$ matrix and the Yang--Baxter equation were discovered in the studies of exactly solvable two-dimensional statistical models and quantum integrable systems. One of the most important results was the quantization of the inverse scattering method by Leningrad's school~\cite{FST79} that lead to the formulation of quantum groups in the so-called RTT formalism \cite{FRT90}. These quantum groups are deformations of semisimple Lie algebras and are closely associated to quantum integrable systems. In particular, the representation theory of the Yangian $\cY({\mathfrak{sl}}_2)$, which is one of the simplest examples of the infinite-dimensional quantum groups, is used to solve the rational $6$-vertex statistical model \cite{Bax82}, the Heisenberg (XXX) spin chain~\cite{FadTak81}, the principal chiral field model with the ${\rm SU}(2)$ symmetry group \cite{FadRes86, Mac04}.

The mathematical description of quantum groups and of quantization of Lie bi-algebras was presented by Drinfeld in his seminal work~\cite{Dri85} (see also~\cite{Dri87}). Drinfeld gave a quantization procedure for the universal enveloping algebra $\cU(\tilde\fg)$ for any semisimple Lie algebra~$\tilde\fg$.\footnote{Throughout this manuscript we use $\tilde\fg$ to denote any Lie algebra, while the undecorated~$\fg$ is reserved for finite-dimensional simple complex Lie algebras.} The quantization is based on the Lie bi-algebra structure on $\tilde\fg$ given by a skew symmetric map $\delta\colon \tilde\fg \to \tilde\fg \wedge \tilde\fg $, the cocommutator. A quantization of $(\tilde\fg,\delta)$ is a (topological) Hopf algebra $(\cU_\hbar(\tilde\fg),\Delta_\hbar)$, such that $\cU_\hbar(\tilde\fg)/\hbar \cU_\hbar(\tilde\fg)\cong\cU(\tilde\fg)$ as a Hopf algebra and
\begin{gather*}
\delta(x) =  \big(\Dh(X)-\sigma\circ\Dh(X)\big) / \hbar \quad ({\rm mod} \ \hbar) ,
\end{gather*}
where $\sigma (a\ot b)=b\ot a$ and $X$ is any lifting of $x\in\tilde\fg$ to $\cU_\hbar(\tilde\fg)$. The Lie bi-algebra structure on $\tilde\fg$ can be constructed from the Manin triple $(\tilde\fg,\tilde\fg_+,\tilde\fg_-)$, where~$\tilde\fg_\pm$ are the isotropic subalgebras of~$\tilde\fg$ such that $\tilde\fg_+\op\tilde\fg_-=\tilde\fg$ as a vector space and $\tilde\fg_-\cong\tilde\fg_+^*$, the dual of $\tilde\fg_+$. Then the commutation relations of the quantum group can be obtained by requiring $\Dh$ to be a homomorphism of algebras, $\Delta_\hbar \colon \cU_\hbar(\tilde\fg)\to \cU_\hbar(\tilde\fg)\ot \cU_\hbar(\tilde\fg)$. The question of the existence of such a quantization for any Lie bi-algebra was raised by Drinfeld in~\cite{Dri92} and was answered by P.~Etingof and D.~Kazhdan in~\cite{EtiKaz96}. They proved that any finite- or infinite-dimensional Lie bi-algebra admits a quantization. Here we will consider only the Yangian case, $\cU_\hbar(\tilde\fg)=\cY(\fg)$ with $\tilde\fg=\cL^+$. We will use the so-called {\it Drinfeld~J} presentation which is very convenient to approach the quantization problem.

In physics, quantum groups are related to quantum integrable models without boundaries and their extensions to models with boundaries. The underlying symmetry of models with boundaries is given by coideal subalgebras of quantum groups that were introduced in the context of (1+1)-dimensional quantum field theories on a half-line by Cherednik~\cite{Che84} and in the context of one-dimensional spin chains with boundaries by Sklyanin \cite{Skl88} in the so-called {\it reflection equation algebra} formalism. Mathematical aspects of reflection algebras in the RTT presentation, called twisted Yangians, were first considered by G.~Olshanskii in~\cite{Ols90} and were further explored in \cite{GR,GRW1,GRW2,Mo92,Mo97,MNO96, MolRag02} (also see references therein). A similar approach for the $q$-deformed universal enveloping and quantum loop Lie algebras, $\cU_q(\fg)$ and $\cU_q(\widehat\fg)$, was considered in~\cite{BasBel11, GCM13, MRS03}.

A slightly different approach to coideal subalgebras, using {\it Serre--Chevalley} presentation of~$\cU_q(\fg)$ was surveyed by G.~Letzter in~\cite{Let02}. Here coideal subalgebras of $\cU_q(\fg)$ were constructed using quantum symmetric pairs, see, e.g.,~\cite{Koo93}, and have been classified for all finite-dimensional semisimple Lie algebras. Examples of coideal subalgebras of quantum affine Lie algebras of type~$A$ in the Serre--Chevalley presentation were the $q$-Onsager algebra \cite{Bas05a,Bas05b} and the generalized $q$-Onsager algebras \cite{BasBel10}. A generalization of quantum symmetric pairs for infinite-dimensional Lie algebras of type $A$ was proposed in \cite{Reg12} and the general theory of quantum symmetric pairs for Kac--Moody Lie algebras was developed in~\cite{Kol12}. A~remarkably simple approach to coideal subalgebras based on the Drinfeld's original construction of Yangians, usually referred to as Drinfeld first or simply as Drinfeld J presentation, was introduced in~\cite{DMS01} (see also~\cite{Mac02}) and is now conveniently called the MacKay twisted Yangians~\cite{GCM13}.

Twisted Yangians, which we denote by $\cY(\fg,\fh)^{\rm tw}$, are in exact correspondence with the symmetric pair decomposition of $\fg$ given by a proper involution $\theta$ \cite{Ara62,Hel78}. The decomposition is given by $\fg=\fh\op\fm$ with $\theta(x)=x$ for all $x\in\fh$, $\theta(y)=-y$ for all $y\in\fm$, the positive and negative eigenspaces of $\theta$. For the half-loop Lie algebra this decomposition is given by $\cL^+=\cH^+\op \cM^+$, where the positive eigenspace of $\theta$ is the twisted current Lie algebra $\cH^+\cong(\fh\oplus u\fm)[u^2]$. The reflection equation algebra formalism provides evidence of the existence of a quantization of such a symmetric pair. More generally, P.~Etingof and D.~Kazhdan have shown in \cite{EtiKaz95} that any homogeneous space~$G/H$ admits a local quantization. In~\cite{BelCra12}, the notion of left and right Lie coideal structures\footnote{The term ``Lie bi-ideal structure'' is used in~\cite{BelCra12} instead.} $\theta$-invariant Lie subalgebras, denoted $\tau$ and $\tau'$, respectively, was introduced by one of the authors. The Lie coideal structure is in one-to-one correspondence with a twisted version of the Manin triple $(\tilde\fg,\tilde\fg_+,\tilde\fg_-)$ and gives a quantization scheme for the $\theta$-invariant Lie subalgebra leading to a quantum coideal subalgebra. The Lie coideal structure corresponds to the splitting $\delta=\tau+\tau'$ of the so-called coisotropic cocommutator~$\delta$ of a Lie bi-algebra satisfying $\delta(\cH^+)\in \cM^+ \wedge \cH^+$ \cite{Oha10,Zam11}. Here we consider the left Lie coideal structure, which allows us to construct a quantization of twisted half-loop Lie algebras leading to quantum coideal subalgebras called twisted Yangians in Drinfeld~J presentation.

A quantization of the twisted half-loop Lie algebras for $\operatorname{rank}(\fg)=1$ case was shown in~\cite{BelCra12}. In this paper we consider twisted Yangians $\cY(\fg,\fh)^{\rm tw}$ in the Drinfeld J presentation for any $\operatorname{rank}(\fg)\geq2$. In particular, we present a quantization procedure which holds for all symmetric pairs of the half-loop Lie algebras of finite-dimensional simple complex Lie algebras of rank $\fg \geq 2$. These symmetric pairs follow from the ones of the simple complex Lie algebras that have been classified by Araki~\cite{Ara62} (see also~\cite{Hel78}). Such symmetric pairs can be grouped into four classes: $\fh$ is simple, $\fh=\fa\op\fb$, $\fh=\fa\op\fk$ and $\fh=\fa\op\fb\op\fk$, where $\fa$ and $\fb$ are simple Lie subalgebras of $\fh$, and $\fk$ is a one-dimensional Lie algebra. We will further refer to all of these cases as the $\theta\neq id$ case. We also consider the trivial involution $\theta$ of~$\fg$. Being trivial at the Lie algebra level, this involution can be extended non-trivially to the half-loop Lie algebra. In this case the $\theta$-fixed subalgebra is the even half-loop Lie algebra~$\fg[u^2]$ isomorphic to~$\fg[u]$ as a~Lie algebra. We will refer to this setup as the $\theta=\operatorname{id}$ case. The main results of this paper are Theorems~\ref{T:2} and~\ref{T:3}, which are analogues of the Theorem~2 in~\cite{Dri85} for the twisted Yangians for $\theta\neq id$ and $\theta= \operatorname{id}$ cases, respectively. We call the defining relations of the twisted Yangians the {\it horrific} relations due to their complex form and similarity to the Drinfeld {\it terrific} relations. A~proof of the horrific relations of Theorem~\ref{T:2} is given in Section~\ref{sec:7.5}. A proof of the horrific relations of Theorem~\ref{T:3} is outlined in Section~\ref{sec:7.6}.

The results of this paper provide a uniform way of constructing twisted Yangians for all symmetric pairs $(\fg,\fh)$, when $\fg$ is a simple complex Lie algebra. Twisted Yangians in Drinfeld J presentation have important applications in quantum integrable models, where it is important to know the minimal set of the defining relations of the underlying symmetry algebra. For example, twisted Yangians emerge as the non-abelian symmetries of the principal chiral models defined on a half-line and can be used to find the scalar boundary $S$-matrices and thus solve the spectral problem \cite{DMS01, Mac02, Mac04}. The closure relations are necessary in constructing representation theory of these algebras and allow to classify algebraically all scalar and dynamical boundary conditions and to construct the corresponding dynamical boundary $S$-matrices; i.e., for the affine Toda models with open boundaries this was done in \cite{BasBel10,BelFom12}. Twisted Yangians in Drinfeld J presentation may also be used to solve the spectral problem of a semi-infinite XXX spin chain for an arbitrary simple Lie algebra using the ``Onsager method'' \cite{BasBel13}. We also remark, that twisted Yangians of this type were shown to play an important role in quantum integrable systems for which the RTT presentation of the underlying symmetries is not well understood, for example in the AdS/CFT correspondence \cite{MacReg10,MacReg11}.

The paper is organized as follows. In Section \ref{Sec:2} we recall the basic facts about simple complex Lie algebras and define the symmetric pair decomposition with respect to involution $\theta$. Then, in Section \ref{Sec:3}, we recall the description of the half-loop Lie algebra $\cL^+$ of $\fg$ and twisted half-loop Lie algebra $\cH^+$ with respect to the symmetric pair decomposition of $\cL^+$. In Section \ref{Sec:4} we construct the Lie bi-algebra structure on $\cL^+$ and Lie coideal structure on $\cH^+$ that provide the necessary setup for the quantization procedure presented in Section~\ref{Sec:5}. A low rank example is presented Section \ref{Sec:6}. We give the explicit form of twisted Yangians in Drinfeld J presentation when $\fg=\fsl_3$. Section \ref{Sec:7} contains the proofs which were omitted in the main part of the paper due to their length and for convenience of the reader.

\section{Definitions and preliminaries}\label{Sec:2}

\subsection{Lie algebra}

Consider a finite-dimensional simple complex Lie algebra $\fg$ of dimension $\dim (\fg)=n$, with a~basis~$\{x_a\}$ such that
\begin{gather}\label{Lie:g}
[x_a,x_b]= {\alpha}_{ab}^{\sk c} x_c, \qquad \alpha_{ab}^{\sk c}+\alpha_{ba}^{\sk c}=0, \qquad \alpha_{ab}^{\sk c} \alpha_{dc}^{\sk e}+ \alpha_{da}^{\sk c} \alpha_{bc}^{\sk e}+ \alpha_{bd}^{\sk c} \alpha_{ac}^{\sk e}=0 .
\end{gather}
Here $\alpha_{ab}^{\sk c}$ are the structure constants of $\fg$ and the Einstein summation rule of the repeated indices is assumed. Let $\eta_{ab}$ denote the non-degenerate invariant bilinear (Killing) form on~$\fg$ in the~$\{x_a\}$ basis
\begin{gather*}%\label{K-form-x}
(x_a,x_b)_\fg= \eta_{ab},
\end{gather*}
that can be used to lower indices \{$a,b,c,\ldots$\} of the structure constants
\begin{gather*}
{\alpha}_{ab}^{\sk d} \eta_{dc}={\alpha}_{abc} \qquad \text{with} \quad {\alpha}_{abc}+{\alpha}_{acb}=0.
\end{gather*}
The inverse of $\eta_{ab}$ is given by $\eta^{ab}$ and satisfies $\eta_{ab}\eta^{bc}=\delta_a^{\; c}$. Set $\{x,y\}=\tfrac{1}{2}(xy+yx)$. Let $C_\fg=\eta^{ab}\{x_a, x_b\}$ denote the second order Casimir operator associated to $\fg$ and let $\fc_{\fg}$ be its eigenvalue in the adjoint representation. For a simple Lie algebra it is non-zero and is given by
\begin{gather} \label{VPCx:g}
\fc_\fg \delta_{c}^{\;d}=\eta^{ab}{\alpha}_{ac}^{\sk e}{\alpha}_{be}^{\sk d}={\alpha}_{c}^{\;eb}{\alpha}_{be}^{\sk d} .
\end{gather}
Here $\alpha_a^{\;bc}$ are the structure constants in the dual basis, namely if $\{x^a\}$ denote the basis of $\fg$ defined by $(x_a,x^b)_\fg=\delta_a^{\;b}$, then $[x^b,x^c]=\alpha_a^{\;bc} x^a$. Constants $\alpha_a^{\;bc}$ satisfy dual co-Jacobi identity, which is obtained by raising one of the lower indices of the Jacobi identity in~\eqref{Lie:g}. Moreover, contracting $\al_a^{\;bc}$ with the Lie commutator in~\eqref{Lie:g} gives
\begin{gather} \label{x=[x,x]}
{\alpha}_{a}^{\;bc}[x_c,x_b] = \fc_\fg x_a .
\end{gather}

\subsection{Symmetric pair decomposition}

Let $\theta \colon \fg \to \fg$ be an involution, $\theta^2 = \operatorname{id}$. Then $\fg$ can be decomposed into positive and negative eigenspaces of $\theta$, i.e., $\fg=\fh\op \fm$ with $\theta(x)=x$ for all $x\in\fh$ and $\theta(y)=-y$ for all $y\in\fm$, here $\dim (\fh)=h$, $\dim (\fm)=m$ satisfying $h+m=n$. Numbers $h$ and $m$ correspond to the number of positive and negative eigenvalues of $\theta$. This decomposition leads to the symmetric pair relations
\begin{gather*}
[\fh,\fh] \subseteq \fh, \qquad [\fh,\fm] = \fm, \qquad [\fm,\fm] = \fh .
\end{gather*}

From the classification of symmetric pairs for simple complex Lie algebras it follows that the invariant subalgebra $\fh$ is a (semi)simple or reductive Lie algebra which can be decomposed into a direct sum of two simple complex Lie algebras $\fa$ and $\fb$, and a one-dimensional Lie algebra $\fk$, at most (see, e.g., \cite[Section~5]{Hel78}). We write $\fh=\fa \op \fb \op \fk$. Set $\dim (\fa)=a$ and $\dim (\fb)=b$. Let the elements
\begin{gather}
X_i \in\fa, \qquad X_{i'}\in \fb, \qquad X_z\in\fk, \qquad Y_p\in \fm,\nonumber\\
\text{with} \quad i=1,\dots, a,\quad i'=1,\dots, b \quad \text{and} \quad p=1,\dots, m ,\label{NewBasis}
\end{gather}
be a basis of $\fg$ such that $\theta(X_\alpha)=X_\alpha$ for any $\alpha\in\{i,i',z\}$, and $\theta(Y_p)=-Y_p$. Here we use indices $i(j,k,\dots)$ for elements $X_\alpha\in\fa$, primed indices $i'(j',k',\dots)$ for elements $X_\alpha\in\fb$, index $\alpha=z$ for the element $X_\alpha\in\fk$, and indices $p(q,r,\ldots)$ for $Y_p\in \fm$, when needed. Note that $X_z$ is central in~$\fh$. We will denote the commutators in this basis as follows:
\begin{gather}
 [X_\alpha,X_\beta] = f^{\sk \gamma}_{\alpha\beta} X_\gamma \qquad \text{with} \quad f^{\sk \gamma}_{\alpha\beta}=0 \ \text{if} \ \al\neq\bet \quad \text{or} \quad \al\neq\ga \quad \text{or} \quad \bet\neq\ga, \nonumber\\
[X_\alpha,Y_p]= g^{\sk q}_{\alpha p} Y_q, \qquad [Y_p, X_\alpha] = g^{\sk q}_{p\alpha} Y_q,
\qquad [Y_p,Y_q]=\mysum_{\alpha}w^{\sk \alpha}_{pq} X_\alpha.\label{SS:coms}
\end{gather}
The structure constants above are obtained from the ones of~$\fg$ by restricting to the appropriate elements. Here and further we will use the sum symbol $\sum_{\alpha}$ to denote the summation over all simple subalgebras of~$\fh$, e.g., $\sum_{\alpha}w^{\sk \alpha}_{pq} X_\alpha=w^{\sk i}_{pq} X_i+w_{pq}^{\sk i'}X_{i'}+w_{pq}^{\sk z}X_z$ for $\fh=\fa \op \fb \op \fk$. The Einstein summation rule for the Greek indices will be used in cases when the sum is over a single simple subalgabera of~$\fh$. The notation $\al\neq\ga$ means that indices~$\al$ and~$\ga$ correspond to different subalgebras of $\fh$. The structure constants given above satisfy the (anti-)symmetry relations
\begin{gather*}
f^{\sk \gamma}_{\alpha\beta}+f^{\sk \gamma}_{\beta\alpha} =0, \qquad g^{\sk q}_{\mu p} + g^{\sk q}_{p\mu} =0, \qquad w^{\sk\mu}_{pq} + w^{\sk \mu}_{qp} =0,
\end{gather*}
which follow from the properties of the Lie bracket, and the homogeneous and mixed Jacobi identities
\begin{gather}
f_{\alpha \beta}^{\sk \nu}f_{\gamma \nu}^{\sk \mu} + f_{\gamma \alpha }^{\sk \nu}f_{\beta \nu}^{\sk \mu} + f_{\beta \gamma}^{\sk \nu}f_{\alpha \nu}^{\sk \mu} = 0, \qquad f^{\sk \mu}_{\alpha\beta} g^{\sk s}_{p\mu} + g^{\sk q}_{\beta p} g^{\sk s}_{\alpha q} - g^{\sk q}_{\alpha p} g^{\sk s}_{\beta q} = 0, \nonumber\\
w^{\sk \beta}_{pq} f^{\sk \mu}_{\alpha\beta} + g^{\sk r}_{\alpha p} w^{\sk \mu}_{qr}- g^{\sk r}_{\alpha q} w^{\sk \mu}_{pr} = 0 ,\label{MixJac}
\end{gather}
with $\alpha(\beta, \gamma,\dots)=i(j,k,\dots)\in \fa$ or $\alpha(\beta, \gamma,\dots)=i'(j',k',\dots)\in \fb$ and
\begin{gather}
\mysum_{\alpha}(w^{\sk \alpha}_{pq} g^{\sk s}_{r\alpha} + w^{\sk \alpha}_{qr} g^{\sk s}_{p\alpha} + w^{\sk\alpha}_{rp} g^{\sk s}_{q\alpha}) = 0 , \qquad
g_{p\alpha}^{\sk r}w_{qr}^{\sk \beta} - g_{q\alpha}^{\sk r}w_{pr}^{\sk \beta} = 0 \quad\text{for} \ \al\neq\bet. \label{JacSym}
\end{gather}
We will further refer to the set $\{X_\al,Y_p\}=\{X_i, X_{i'},X_z, Y_p\}$ given by~\eqref{NewBasis} and satisfying relations~(\ref{SS:coms})--(\ref{JacSym}) as the {\it symmetric space basis} for a given Lie algebra~$\fg$ and involution~$\theta$.

Let $\kappa$ denote the Killing form on $\fg$ in the symmetric space basis. It has a block-diagonal form, namely
\begin{gather*}
(X_i,X_j)_\fg=(\kappa_\fa)_{ij}, \qquad ({X}_{i'},{X}_{j'})_\fg=(\kappa_\fb)_{i'j'}, \\
(X_z,X_z)_\fg=(\kappa_\fk)_{zz}, \qquad (Y_p,Y_q)_\fg=(\kappa_\fm)_{pq},
\end{gather*}
with the remaining entries being trivial. The Casimir element $C_\fg$ in this basis decomposes as
\begin{gather*}
C_\fg = C_X+C_{Y} = \mysum_{\al,\bet}(\ka_\fh)^{\alpha\beta}\{X_\alpha,X_\beta\}+(\kappa_\fm)^{pq}\{Y_p,Y_q\},\\
C_X = C+C'+C_{z}= (\kappa_\fa)^{ij}\{X_i,X_j\}+(\kappa_\fb)^{i'j'}\{X_{i'},X_{j'}\}+(\kappa_\fk)^{zz} \{ X_z,X_z \} .
\end{gather*}
Here $\kappa^{\alpha\beta}\in\{(\kappa_\fa)^{ij},(\kappa_\fb)^{i'j'},(\kappa_\fk)^{zz}\}$. The block diagonal decomposition of the inverse Killing form can be used to raise the indices of the structure constants. We set
\begin{gather*}
f_\alpha^{\;\beta\gamma} = \kappa^{\beta \mu} f_{\alpha\mu}^{\sk \gamma},\qquad
 g_p^{\;q\nu} = \kappa^{\nu\rho} g_{\rho p}^{\sk q} , \qquad w_\nu^{\;pq} = (\kappa_\fm)^{pr} g_{\nu r}^{\sk q}
\end{gather*}
with $\alpha(\beta, \gamma,\mu,\dots)=i(j,k,\dots)$ or $i'(j',k',\dots)$ and $\nu(\rho)= i(j)$ or $i'(j')$ or $z(z)$. When the two upper indices of $g$ are written in the opposite order, we set $g_p^{\;\nu q}=-g_p^{\;q\nu}$. Thus $f_\alpha^{\;\beta\gamma}$, $g_p^{\;q\nu}$, $g_p^{\;\nu q}$ and $w_\nu^{\;pq}$ are components of $\alpha_a^{\;bc}$, which is antisymmetric in its upper indices. Let us now consider the commutation relations. For generators $Y_p$ we have
\begin{gather*}
[Y_p,C_X] =2\mysum_{\al}g_{p}^{\;\alpha q}\{Y_q,X_\alpha\}, \qquad [Y_p,C_Y] =2\mysum_{\al}g_{p}^{\;q\alpha}\{Y_q,X_\alpha\}.
\end{gather*}
The remaining commutation relations are trivial. Let $\fc_\fa$, $\fc_\fb$, $\fc_z$ and $\fc_\fm$ be the eigenvalues of~$C$,~$C'$, $C_z$ and $C_{Y}$ in the adjoint representation, respectively. We have $\fc_\fg=\fc_\fa+\fc_\fb+\fc_\fm+\fc_z$. Using~\eqref{x=[x,x]} we find
\begin{gather}
f_\al^{\;\bet \nu} [X_\nu,X_\bet]=\fc_{(\al)} X_\al, \qquad w_\gamma^{\;qp} [Y_p,Y_q] = \bar\fc_{(\gamma)} X_\gamma, \label{inv:Y}
\end{gather}
with $\al=i(i')$, $\gamma=i(i',z)$, $\bar\fc_{(\alpha)}= \fc_\fg-\fc_{(\alpha)}$, $ \fc_{(i)}= \fc_\fa$, $ \fc_{(i')}= \fc_\fb$ and $ \fc_{(z)}=0$. Using~\eqref{VPCx:g},~\eqref{SS:coms} and the equality above we obtain
\begin{gather}
f_\al^{\;\mu\nu}f^{\sk \bet}_{\nu\mu} = \fc_{(\al)} \delta_\al^{\;\bet}, \qquad w_\al^{\;qp} w^{\sk \bet}_{pq} = \bar\fc_{(\al)} \delta_\al^{\;\bet}, \qquad w_{\alpha}^{\; qp} w^{\sk \gamma}_{pq} =0 \quad \text{for} \ \al\neq\ga, \label{wwid}
\end{gather}
and $\al=i(i',z)$, $\bet=j(j',z)$. Finally, for $Y_q$ we have
\begin{gather} \label{ggid}
{\mysum_{\alpha}} g^{\;p\alpha}_q [X_\alpha,Y_p] = \tfrac12 \fc_\fg Y_q,
\qquad {\mysum_{\alpha}} g^{\;r\alpha}_{p} g_{\alpha r}^{\sk q} = \tfrac12 \fc_\fg \delta_{p}^{\;q} .
\end{gather}

\section{Symmetric spaces and simple half-loop Lie algebras} \label{Sec:3}

\subsection{Half-loop Lie algebra}

Consider the half-loop Lie algebra $\cL^+$ generated by elements $x^{(k)}_a$ with $k\ge0$ and $a=1,\dots$, $\dim (\fg)$. It is an infinite-dimensional graded Lie algebra with the grading given by $\deg \big(x^{(k)}_a\big)=k$ and the defining relations
\begin{gather} \label{Lp[x,x]}
[x^{(k)}_a,x^{(\ell)}_b]={\alpha}^{\sk c}_{ab} x^{(k+\ell)}_c.
\end{gather}
This algebra can be identified with the set of polynomial maps $f \colon \CC \to \fg$ using the Lie algebra isomorphism $\cL^+\cong\fg[u] = \fg \ot \CC[u]$ with~$x^{(k)}_a \cong x_a\ot u^{k}$.

The algebra $\cL^+$ has another presentation, conveniently called Drinfeld~J presentation.

\begin{Proposition} \label{P:31}The half-loop Lie algebra $\cL^+$ is isomorphic to the algebra generated by elements~$x_a$, $J(x_b)$ satisfying
\begin{gather}
[x_a,x_b]={\alpha}^{\sk c}_{ab} x_c, \qquad J(\mu x_a+ \nu x_b) = \mu J(x_a)+\nu J(x_b), \qquad [x_a,J(x_b)]={\alpha}^{\sk c}_{ab} J(x_c) , \label{DL01}\\
[J(x_a),J([x_b,x_c])] + [J(x_b),J([x_c,x_a])] + [J(x_c),J([x_a,x_b])] = 0 , \label{DL2}\\
[[J(x_a),J(x_b)],J([x_c,x_d])] + [[J(x_c),J(x_d)],J([x_a,x_b])] = 0 , \label{DL3}
\end{gather}
for any $\mu,\nu\in\CC$. The grading is given by $\deg(x_a)=0$ and $\deg(J(x_a))=1$. The isomorphism is given by the map $x_a \mapsto x^{(0)}_a$ and $J(x_a) \mapsto x^{(1)}_a$.
\end{Proposition}

Relations \eqref{DL2} and \eqref{DL3} are homogeneous relations of degree $2$ and $3$. In the $\operatorname{rank}(\fg)=1$ case the relation \eqref{DL2} is trivial and for the $\operatorname{rank}(\fg)\geq 2$ cases the relation~\eqref{DL3} is implied by~\eqref{DL2}. Since we were unable to locate a complete proof of this proposition in the literature, we have given one in Section~\ref{sec:7.1}. (An outline of a proof is given in \cite[Section~12.1]{ChaPre94}.)

\subsection{Twisted half-loop Lie algebra}

Let $\theta(x_a) = \theta_a^{\; b} x_b$. We can extend the involution $\theta$ of $\fg$ to the whole of $\cL^+$ as follows:
\begin{gather*} %\label{thetaLp}
\theta\big(x^{(k)}_a\big) = (-1)^k \theta_a^{\; b} x_b^{(k)} \cong (-1)^k \theta(x_a) \ot u^k \qu\text{for}\qu k\ge 0,
\end{gather*}
or in the polynomial map point of view, $\theta(x(u))= \theta(x)(-u)$.

The twisted half-loop Lie algebra $\cH^+\cong\fg[u]^{\theta}$ is a fixed-point subalgebra of $\cL^+$ generated by the elements stable under the action of the (extended) involution $\theta$, namely $\cH^+=\operatorname{span}\{x\in\cL^+\,|\,\theta(x)=x\}$. In the physics literature it is often referred to as the {\it twisted current algebra}.

Let $\theta\neq \operatorname{id}$. Consider the symmetric space basis of~$\fg$. We write the half-loop Lie algebra $\cL^+$ in terms of the elements~$X^{(k)}_\alpha$ and $Y^{(k)}_q$ satisfying
\begin{gather}
 \big[X^{(k)}_\alpha,X^{(\ell)}_\beta\big] = {f}^{\sk \gamma}_{\alpha\beta} X^{(k+\ell)}_\gamma, \qquad
\big[X^{(k)}_\al,Y^{(\ell)}_p\big]=g^{\sk q}_{\al p} Y^{(k+\ell)}_q, \nonumber\\
\big[Y^{(k)}_p,Y^{(\ell)}_q\big] = \mysum_{\alpha}w^{\sk \alpha}_{pq} X^{(k+\ell)}_\alpha,\label{LpXY}
\end{gather}
for all $k,\ell\ge0$. Involution $\theta$ acts on these elements by $\theta\big(X_\al^{(k)}\big)=(-1)^kX_\al^{(k)}$ and $\theta\big(Y^{(k)}_p\big)=(-1)^{k+1}Y^{(k)}_p$. Thus the half-loop Lie algebra decomposes as $\cL^+=\cH^+\op\cM^+$, where $\cH^+=\operatorname{span}\big\{X^{(2k)}_\al, Y^{(2k+1)}_q\big\}$ is the subalgebra of $\cL^+$ generated by $\theta$-invariant elements, and
$\cM^+=\operatorname{span}\big\{X^{(2k+1)}_\al,Y^{(2k)}_q\big\}$ is the subset of $\cL^+$ of $\theta$-anti-invariant elements. In the case when $\theta = \operatorname{id}$ we have $\fg[u]^\theta \cong \fg[u^2]$, that is $\cH^+=\operatorname{span}\big\{x^{(2k)}_a\big\}$ and $\cM^+=\operatorname{span}\big\{x^{(2k+1)}_a\big\}$.

We now give another presentation of $\cH^+$, which an analogue of the Drinfeld J presentation of $\cL^+$.

\begin{Proposition} \label{P:32} Let $\operatorname{rank}(\fg)\geq2$ and $\theta\neq \operatorname{id}$. Then the twisted half-loop Lie algebra~$\cH^+$ is isomorphic to the algebra generated by elements~$X_\al$, $B(Y_p)$ satisfying
\begin{gather}
[X_\al,X_\bet]=f^{\sk\, \ga}_{\al\bet} X_\ga, \qquad [X_\al,B(Y_p)] = g^{\sk\, q}_{\al p} B(Y_q),\nonumber\\
 B(a Y_p + b Y_q) = a B(Y_p) + b B(Y_q), \label{L[X,BY]}\\
[B(Y_p),B(Y_q)] + {\mysum_\al} (\bar\fc_{(\al)})^{-1} w_{pq}^{\sk\al}w_{\al}^{\;rs} [B(Y_r),B(Y_s)] = 0 , \label{LH2} \\
[[B(Y_p),B(Y_q)],B(Y_r)] + 2 \fc_\fg^{-1} \mysum_{\al} (\ka_\fm)^{tu} w_{pq}^{\sk\al} g_{r\al}^{\sk s} [[B(Y_s),B(Y_t)],B(Y_u)] = 0 , \label{LH3}
\end{gather}
for all $a,b\in\CC$. The isomorphism is given by the map $X_\al\mapsto X^{(0)}_\al$ and $B(Y_p)\mapsto Y^{(1)}_p$.
\end{Proposition}

The proof is given in Section \ref{sec:7.2}. Note that in the contrast to $\cL^+$, the twisted algebra $\cH^+$ for $\operatorname{rank}(\fg)\geq2$ has degree-2 and degree-3 defining relations. The $\operatorname{rank}(\fg)=1$ case is exceptional. The Drinfeld~J presentation in this case has a degree-$4$ relation instead; see \cite[Section~4.2]{BelCra12}.

\begin{Proposition} \label{Prop35}
Let $\operatorname{rank}(\fg)\geq2$ and $\theta = \operatorname{id}$. Then the even half-loop Lie algebra is isomorphic to the algebra generated by elements~$x_i$, $G(x_j)$ satisfying
\begin{gather}
[x_i,x_j] = \al^{\sk k}_{ij} x_k, \!\qquad G(\la x_a + \mu x_b) = \la G(x_a) +\mu G(x_b) ,\! \qquad [x_i,G(x_j)]=\al^{\sk k}_{ij} G(x_k) , \!\!\\
[G(x_i),G([x_j,x_k])] + [G(x_j),G([x_k,x_i])] + [G(x_k),G([x_i,x_j])] = 0 , \label{LH4}
\end{gather}
for any $\la,\mu \in\CC$. The isomorphism is given by the map $x_i\mapsto x^{(0)}_i$ and $G(x_i)\mapsto x^{(2)}_i$.
\end{Proposition}

The proof is analogous to that of Proposition \ref{P:31}, since $\fg[u^2]\cong \fg[u]$ as a Lie algebra.

\section{Lie bi-algebras and coideals} \label{Sec:4}

\subsection{Lie bi-algebra structure of a half-loop Lie algebra}

A Lie bi-algebra structure on $\cL^+$ is a skew-symmetric linear map $ \delta \colon \cL^+ \to \cL^+ \ot \cL^+$, the cocommutator, such that $ \delta^*$ is a Lie bracket and $\delta $ is a 1-cocycle, $\delta([x,y])=x.\delta(y)-y.\delta(x)$, where dot denotes the adjoint action on $\cL^+ \ot \cL^+$. The cocommutator is given for the elements in the {Drinfeld J} presentation of $\cL^+$ by
\begin{gather} \label{deltaJ}
\delta(x_a)=0, \qquad \delta(J(x_a))=[x_a \ot 1, \Omega_\fg], \qquad \text{where}\quad \Omega_{\fg}=\eta^{ab} x_a \ot x_b.
\end{gather}
This cocommutator can be constructed from the Manin triple $(\cL, \cL^+, \cL^-)$, with $\cL=\fg((u^{-1}))$ the loop algebra generated by elements $x^{(n)}$ with $x\in \fg$, $n\in\ZZ$ and defining relations~(\ref{Lp[x,x]}) (but with $n,m\in\ZZ$), $\cL^+=\fg[u]$ the positive half-loop algebra (\ref{Lp[x,x]}), and $\cL^-=u^{-1}\fg[[u^{-1}]]$ the negative half-loop algebra (i.e., $n,m<0$; see, e.g., \cite[Example~1.3.9]{ChaPre94}). The triple $(\cL, \cL^+, \cL^-)$ satisfies axioms of the Manin triple.

\begin{Definition}[\cite{Dri87}]
A Manin triple is a triple of Lie algebras $(\tilde \fg, \tilde\fg^+,\tilde \fg^-)$ together with a~non-degenerate symmetric bilinear form $(\;,\;)_{\tilde \fg}$ on $\tilde \fg$ invariant under the adjoint action of $\tilde \fg$ such that:
\begin{itemize}\itemsep=0pt
\item ${\tilde \fg}^+$ and ${\tilde \fg}^-$ are Lie subalgebras of ${\tilde \fg}$;
\item ${\tilde \fg}={\tilde \fg}^+ \op{\tilde \fg}^-$ as a vector space;
\item $(\;,\;)_{\tilde \fg}$ is isotropic for $ {\tilde \fg}^\pm$ (i.e., $({\tilde \fg}^\pm,{\tilde \fg}^\pm)_{\tilde\fg}=0$);
\item $({\tilde \fg}^+)^* \cong {\tilde \fg}^-$.
\end{itemize}
\end{Definition}

The invariant bilinear form on $\tilde\fg = \cL$ satisfying the requirements above is given by
\begin{gather*}
\big(x^{(k)},y^{(l)}\big)_\cL = -(x,y)_\fg \delta_{k+l+1,0} , %\label{form_L}
\end{gather*}
or in the polynomial map point of view $(x(u),y(u))_\cL=-\operatorname{res}_0 (x(u),y(u))_\fg$ where `$\operatorname{res}_0$' means taking the coefficient of $u^{-1}$ in the Laurent series expansion. This data together with the definition of the Manin triple uniquely fixes the cocommutator~$\delta$ on~$\cL^+$.

\begin{Remark}
If $({\tilde \fg},{\tilde \fg}^+,{\tilde \fg}^-)$ is a Manin triple for $\dim ({\tilde \fg}^+)=\infty$, then in general one must take $(\tilde\fg^+)^* \cong \overline{\fg}{^-}$, where $\overline{\fg}{^-}$ is a suitable completion of ${\tilde \fg}^-$. However in our case $(\cL^+)^* \cong {\cL}^-$, as it is easy to see:
\begin{gather*}
\textstyle (\cL^+)^* \cong \big( \bigoplus_{k\geq0} \fg\ot u^k \big)^* = \prod_{k\geq0} (\fg\ot u^k)^* = \prod_{k\geq1} \fg\ot u^{-k} \cong \cL^- .
\end{gather*}
Here in the second equality we have used the identity $(\bigoplus_{i\geq0} V_i)^* = \prod_{i\geq0} V_i^*$, where $V_i$ denotes a finite-dimensional vector space; an equivalent identity is used in the last equality.
\end{Remark}

The cocomutator is obtained using the duality between $\cL^+$ and $\cL^-$. Recall that $\delta^* \colon \cL^-\ot \cL^- \to \cL^-$ is the Lie bracket of $\cL^-$. We can deduce the cocommutator $\delta$ of $\cL^+$ from the duality relation
\begin{gather} \label{DualR}
(\delta(x),y\ot z)_{\cL\ot\cL}=(x,[y,z])_\cL .
\end{gather}
The cocommutator of the degree zero generators $x^{(0)}_a=x_a$ is trivial since the degree of $[y,z]\in \cL^-$ is strictly less than $-1$, thus $(x_a,[y,z])_\cL=0$ for all $[y,z]\in \cL^-$. This implies that $(\delta(x_a),y\ot z)_{\cL\ot\cL}=0$ for all $y,z\in \cL^-$ and thus
\begin{gather*} \delta(x_a)=0.
\end{gather*}
The case of degree one generators $x^{(1)}_a=J(x_a)$ is considered in a similar way. For this case we have the non-trivial pairing $\big(J(x_a),x^{(-2)}_b\big)_\cL=-\eta_{ab}$. Using $x^{(-2)}_b=\fc_\fg^{-1}\alpha_b^{\;ji} \big[x^{(-1)}_i,x^{(-1)}_j\big]$ and the duality relation~\eqref{DualR} we obtain the constraint
\begin{gather*}
\big(\delta(J(x_a)),\alpha_b^{\;ji} x^{(-1)}_i\ot x^{(-1)}_j\big)_{\cL\ot\cL}=-\fc_\fg \eta_{ab} .
\end{gather*}
Writing $\delta(J(x_a))=v_a^{\;\ell k} x_k\ot x_\ell$ for some $v_a^{\;\ell k}\in\CC$ we deduce that $v_a^{\;\ell k}\alpha_{b\ell k}= -\fc_\fg \eta_{ab} $. By~\eqref{VPCx:g} we have that $\al_{a}^{\;k\ell}\al_{b\ell k} = \fc_\fg \eta_{ab}$ and thus $v_a^{\;\ell k}=-\al_a^{\;k\ell}$ giving
\begin{gather*}
\delta(J(x_a))=\alpha_a^{\;\ell k} x_k\ot x_\ell=[x_a \ot 1,\Omega_\fg].
\end{gather*}

\subsection{Lie coideal structure of a twisted half-loop Lie algebra}

The Lie coideal structure of a twisted half-loop Lie algebra is constructed by employing the anti-invariant Manin triple twist. Here we will consider the left Lie coideal structure. The right Lie coideal is obtained in a similar way.

\begin{Definition} [\cite{BelCra12}]
The anti-invariant Manin triple twist $\phi$ of $(\cL, \cL^+, \cL^-)$ is an automorphism of~$\cL$ satisfying:
\begin{itemize}\itemsep=0pt
\item $\phi$ is an involution;
\item $\phi(\cL^\pm)=\cL^\pm$;
\item $(\phi(x),y)_\cL=-(x,\phi(y))_\cL$ for all $x\in \cL^+$ and $y\in \cL^-$.
\end{itemize}
\end{Definition}

Given an involution $\theta$ of $\fg$ one may uniquely associate the anti-invariant Manin triple twist~$\phi$, that is the natural extension of~$\theta$ to the whole~$\cL$:
\begin{gather*}
\phi(x_a^{(k)})= (-1)^k\theta_a^{\; b} x_b^{(k)} \cong (-1)^k \theta(x_a) \ot u^k \qquad \text{for all}\quad k\in\ZZ,
\end{gather*}
satisfying the definition above and leading to the symmetric pair decomposition of the Manin triple $(\cL, \cL^+, \cL^-)$:
\begin{gather}
\cL^\pm=\cH^\pm\op \cM^\pm \quad\text{such that}\quad \phi(x)=x, \quad \forall \, x\in\cH^\pm \qquad \text{and}\nonumber\\ \phi(y)=-y, \quad \forall\, y\in \cM^\pm .\label{MTdec}
\end{gather}
From the anti-invariance of $\phi$ for $(\;,\;)_\cL$ it follows that
\begin{gather*}
 (\cH^-,\cH^+)_\cL=(\cM^-,\cM^+)_\cL=0.
\end{gather*}
Thus we must have $(\cH^\pm)^* \cong {\cM}^\mp$. This is easy to check:
\begin{gather*}
(\cH^+)^* \cong \big( \textstyle\bigoplus_{k\geq0} (\fh\op u\fm)\ot u^{2k}\big)^* \\
\hphantom{(\cH^+)^*}{}= \prod_{k\geq0} \big((\fh\op u\fm)\ot u^{2k}\big)^* = \prod_{k\geq1}\big( (u\fh\op \fm)\ot u^{-2k}\big) \cong \cM^- , \\
(\cH^-)^* \cong \big( \textstyle\bigoplus_{k\geq1} (\fh\op u\fm)\ot u^{-2k}\big)^* \\
\hphantom{(\cH^-)^*}{} = \prod_{k\geq1}\big( (\fh\op u\fm)\ot u^{-2k}\big)^* = \prod_{k\geq0}\big( (u\fh\op \fm)\ot u^{2k}\big) \cong \cM^+ .
\end{gather*}
This decomposition of the Manin triple allows us to construct a Lie coideal structure on $\cH^+$.

\begin{Definition} [\cite{BelCra12}]
Let $\phi$ be an anti-invariant Manin triple twist for $(\cL,\cL^+,\cL^-)$ which leads to the symmetric space decomposition~(\ref{MTdec}). Then the linear map $\tau\colon \cH^+\rightarrow \cM^+ \ot \cH^+$ is a left Lie coideal structure for $(\cH^+,\cM^+)$ if it is the dual of the following action of $\cH^-$~on~$\cM^-$,
\begin{gather} \label{t-dual}
\begin{aligned}
 \tau^* \colon \ \cH^-\ot \cM^- & \to \cM^- , \\
 x\ot y &\mapsto [x,y]_{\cL^-},
\end{aligned}
\end{gather}
for all $x\in\cH^-$ and $y\in\cM^-$.
\end{Definition}

Given as Manin triple $(\cL,\cL^+,\cL^-)$ and the twist $\phi$ the requirement \eqref{t-dual} uniquely fixes the Left coideal structure $\tau$ of $(\cL^+,\cH^+)$.

\begin{Proposition}The left Lie coideal structure of $(\cL^+,\cH^+)$, $\tau \colon \cH^+ \to \cM^+ \ot \cH^+$, is given by
\begin{gather}
\label{taunotid} \theta \neq \operatorname{id} \colon \ \tau(X_\alpha)=0, \qu \tau(B(Y_p))=[Y_p \ot 1, \Omega_X], \qquad \Omega_X=\mysum_{\al,\bet}(\kappa_\fh)^{\alpha \beta} X_\alpha \otimes X_\beta, \\
\label{tauid} \theta = \operatorname{id} \colon \ \tau(x_a)=0, \qquad \tau(G(x_a))=[J(x_a) \ot 1, \Omega_\fg] .
\end{gather}
\end{Proposition}

\begin{proof}The construction of the left Lie coideal structure $ \tau$ from the anti-invariant Manin triple twist is similar to the one of the Lie bi-algebra structure from the Manin triple. We have to consider the duality relation $(\tau(x),y\ot z)_\cL=(x,[y,z])_\cL$ with $x\in \cH^+$, $y\in \cH^-$ and $z\in \cM^-$.

Consider the case $\theta \neq id$. For the degree zero generators $X_\alpha^{(0)}=X_\alpha$, we have $(X_\alpha,[y,z])_\cL=0$ for all $y\in \cH^-$ and $z\in \cM^-$. This follows by similar arguments as for the degree zero generators of the half-loop Lie algebra. Hence we have
\begin{gather*}
\tau(X_\alpha)=0.
\end{gather*}

For elements $Y^{(1)}_p=B(Y_p)$ we have a non trivial paring $\big(B(Y_p), Y^{(-2)}_q\big)_\cL=-(\kappa_\fm)_{pq}$. Then, using relation $Y^{(-2)}_q= \sum_{\alpha} 2\fc_\fg^{-1} g_{q}^{\;\alpha p} \big[Y^{(-1)}_p,X^{(-1)}_\alpha\big]$ and the duality relation we obtain
\begin{gather*}
{\mysum_{\alpha}} \big(\tau(B(Y_p)),g_{q}^{\;\al r} Y_r^{(-1)}\ot X^{(-1)}_\alpha\big)_{\cL\ot\cL}=-\tfrac12 \fc_\fg(\kappa_\fm)_{pq}.
\end{gather*}
It is clear from properties of the pairing that $\tau(B(Y_p))= \sum_{\beta} v_{p}^{\;\beta s}Y_s\ot X_\beta$ for some $v_{p}^{\;\beta s}\in\CC$. We must have $\sum_{\alpha}v_{p}^{\;\alpha r}g_{q\al r}=-\frac{1}{2}\fc_\fg(\kappa_\fm)_{pq}$. Writing the second identity in~(\ref{ggid}) as $\mysum_\al g_p^{\;r\al} g_{q\al r} = \tfrac12 \fc_\fg (\ka_\fm)_{pq}$ we find $v_p^{\; \al r}=-g_p^{\;r\al}$. Hence
\begin{gather*}
\tau(B(Y_p))= -\mysum_{\alpha}g_{p}^{\;r\alpha} Y_r\ot X_\alpha=[Y_p\ot 1,\Omega_X].
\end{gather*}
The Lie coideal structure for the $\theta = \operatorname{id}$ case follows from the pairing $\big(G(x_a),x_b^{(-3)}\big)_\cL =- (\kappa_\fg)_{ab}$ by similar arguments.
\end{proof}

For completeness we recall the remark which was stated by one of the authors in~\cite{BelCra12}.

\begin{Remark}The notion of left (right) Lie coideal is related to the notion of co-isotropic subalgebra $\fh$ of a Lie bi-algebra $(\fg, \delta)$. It is a Lie subalgebra which is also a Lie coideal, meaning that $\delta (\fh) \subset \fh \wedge \fg$. We have $\delta(x)=\tau(x) + \tau'(x)$, for $x\in\fh$ with $\tau'=-\sigma\circ\tau$ the right ideal structure.
\end{Remark}

\section{Twisted Yangians as quantized Lie coideals} \label{Sec:5}

\subsection{Quantization of Lie bi-algebras and Lie coideals}

To obtain a quantization of the Lie coideal we need to introduce some additional notation. Recall the definition of a bi-algebra and of a Hopf algebra. A bi-algebra is a quintuple $(A,\mu,\imath,\Delta,\varepsilon)$ such that $(A,\mu,\imath)$ is an algebra and $(A,\Delta,\varepsilon)$ is a coalgebra; here $A$ is a $\CC$-module, $\mu \colon A \ot A \to A$ is the multiplication, $\Delta\colon A\to A\ot A$ is the comultiplication (coproduct), $\imath \colon \CC \to A$ is the unit and $\varepsilon \colon A \to \CC$ is the counit. A Hopf algebra is a bi-algebra with an antipode $S\colon A \to A$, an antiautomorphism of algebra. These maps are required to satisfy a set of compatibility relations (see, e.g., \cite[Section~4]{ChaPre94}). Since we are interested in left coideal subalgebras, we also need to recall the notion of coideal.

\begin{Definition} \label{D51} Let $\cA=(A,\Delta,\varepsilon)$ be a coalgebra. Then $\cB=(B,\LD,\epsilon)$ is a left coideal of $\cA$ if:
\begin{enumerate}\itemsep=0pt
\item[1)] $B$ is a submodule of $A$, i.e., there exists an injective homomorphism $\varphi\colon B \to A$;

\item[2)] co-action $\LD$ is a left coideal map $\LD \colon B \to A \ot B$ and is a homomorphism of modules;

\item[3)] coalgebra and coideal structures are compatible with each other, i.e., the following identities hold:
\begin{gather}
(\Delta \ot \operatorname{id})\circ \LD =(\operatorname{id} \ot \LD )\circ \LD , \label{co-co-ass}
\\
(\operatorname{id} \ot \varphi) \circ \LH = \Delta\circ\varphi; \label{co-inv}
\end{gather}
\item[4)] $\epsilon \colon B \to \CC$ is the counit.
\end{enumerate}
\end{Definition}

The relation \eqref{co-co-ass} is usually referred to as the coideal co-associativity of~$B$. We will refer to~\eqref{co-inv} as the coideal coinvariance. We call the map $\varphi$ the natural embedding $\varphi \colon B \hookrightarrow A$. The definition of a coideal can be naturally extended for a bi-algebra. Namely, let $\cA=(A,\mu,\imath,\Delta,\varepsilon)$ be a bi-algebra. Then $\cB=(B,m,i,\LD,\epsilon)$ is a left coideal of $\cA$ if: the triple $(B,m,i)$, where~$m$ is the multiplication and $i$ is the unit, is an algebra; $B$ is a subalgebra of $A$; the triple $(B,\LD,\epsilon)$ is a coideal of $(A,\Delta,\varepsilon)$. However note that $\epsilon$ is not a restriction of $\varepsilon$, that is $\varepsilon \circ \varphi \ne \epsilon$, since in general~$B$ may have a larger center than~$A$.

We are now ready to discuss quantization of Lie bi-algebras and their coideals. The next definition, a quantization of a Lie bi-algebra~${\tilde \fg}$, is due to V.G.~Drinfeld.

\begin{Definition}[\cite{Dri87}] \label{Qbial} Let $({\tilde \fg},\delta)$ be a Lie bi-algebra. We say that a quantized universal enveloping algebra $(\cU_{\hbar}({\tilde \fg}),\Delta_\hbar)$ is a quantization of $({\tilde \fg},\delta)$, or that $({\tilde \fg},\delta)$ is the quasi-classical limit of $(\cU_{\hbar}({\tilde \fg}),\Delta_\hbar)$, if it is a topologically free $\CC[[\hbar]]$ module and:
\begin{enumerate}\itemsep=0pt
\item[1)] $\cU_{\hbar}({\tilde \fg})/\hbar \cU_{\hbar}({\tilde \fg})$ is isomorphic to $\cU({\tilde \fg})$ as a Hopf algebra;
\item[2)] for any $x \in {\tilde \fg}$ and any $X \in \cU_{\hbar}({\tilde \fg})$ equal to $x$ $(\operatorname{mod} \, \hbar)$ one has
\begin{gather*}
\big(\Delta_\hbar(X)-\sigma\circ\Delta_\hbar(X)\big)/\hbar = \delta(x) \quad (\operatorname{mod}\, \hbar)
\end{gather*}
with $\sigma$ the permutation map $\sigma(a\ot b)=b\ot a$.
\end{enumerate}
\end{Definition}

Note that $(\cU_{\hbar}({\tilde \fg}),\Delta_\hbar)$ is a topological Hopf algebra over $\CC[[\hbar]]$ and is a topologically free $\CC[[\hbar]]$ module. Drinfeld also noted that for a given Lie-bialgebra $(\tilde\fg,\delta)$, there exists a unique extension of the map $\delta\colon \tilde\fg\to\tilde\fg\ot\tilde\fg$ to $\delta \colon \cU(\tilde\fg) \to \cU(\tilde\fg) \ot \cU(\tilde\fg)$ which turns $\cU(\tilde\fg)$ into a co-Poisson--Hopf algebra. The converse is also true. In such a way $(\cU_{\hbar}({\tilde \fg}),\Delta_\hbar)$ can be viewed as a quantization of $(\cU(\tilde\fg),\delta)$.

Let $\theta$ be an involution of $\tilde\fg$. We define a quantization of the Lie coideal $\big({\tilde \fg}^\theta,\tau\big)$ following~\cite{BelCra12}.

\begin{Definition} \label{D53} Let $({\tilde \fg},\delta)$ be a Lie bi-algebra and $\big({\tilde \fg}^\theta,\tau\big)$ be a left Lie coideal. We say that a left coideal subalgebra $\big(\cU_{\hbar}\big({\tilde \fg}^\theta\big),\LHb\big)$ is a quantization of $\big({\tilde \fg}^\theta,\tau\big)$, or that $\big({\tilde \fg}^\theta,\tau\big)$ is the quasi-classical limit of $\big(\cU_{\hbar}\big({\tilde \fg}^\theta\big),\LHb\big)$, if it is a topologically free $\CC[[\hbar]]$ module and:
\begin{enumerate}\itemsep=0pt
\item[1)] $(\cU_{\hbar}({\tilde \fg}),\Delta_\hbar)$ is a quantization of $(\tilde\fg,\delta)$;
\item[2)] $\cU_{\hbar}\big({\tilde \fg}^\theta\big)/ \hbar\cU_{\hbar}\big({\tilde \fg}^\theta\big)$ is isomorphic to $\cU\big({\tilde \fg}^\theta\big)$ as an algebra;
\item[3)] $\big(\cU_{\hbar}\big({\tilde \fg}^\theta\big),\LHb\big)$ is a left coideal of $(\cU_{\hbar}({\tilde \fg}),\Delta_\hbar)$;
\item[4)] for any $x \in {\tilde \fg}^\theta$ and any $X \in \cU_{\hbar}\big({\tilde \fg}^\theta\big)$ equal to $x$ $(\operatorname{mod} \, \hbar)$ one has
\begin{gather*}
\big(\LHb(X) -(\varphi(X) \ot 1 + 1 \ot X)\big)/\hbar = \tau(x) \quad (\operatorname{mod}\, \hbar)
\end{gather*}
modulo symmetric elements of $\cU\big(\tilde\fg^\theta\big)\ot\cU\big(\tilde\fg^\theta\big)$, with $\varphi$ the natural embedding $\cU_{\hbar}\big({\tilde \fg}^\theta\big) \hookrightarrow \cU_{\hbar}({\tilde \fg})$. The symmetric elements must be allowed: replacing the lift $X$ by $X+\hbar Y$ changes the left-hand side, modulo $\hbar$, by $\Delta(\bar y)-\bar y\ot1-1\ot\bar y$, where $\bar y\in\cU\big(\tilde\fg^\theta\big)$ is $Y$ modulo $\hbar$. This element is symmetric, and it is zero only if $\bar y\in\tilde\fg^\theta$.
\end{enumerate}
\end{Definition}

The Lie coideal structure $\tau$ can be extended to $\tau \colon \cU\big(\tilde\fg^\theta\big) \to \cU(\tilde\fg) \ot \cU\big(\tilde\fg^\theta\big)$ such that $\tau(a_1a_2)=\tau(a_1)\LH(a_2)+\LH(a_1)\tau(a_2)$ \cite{BelCra12}. However this does not turn $\big(\cU\big(\tilde\fg^\theta\big),\LH,\tau\big)$ into a co-Poisson--Hopf structure. Rather it would be a~``one sided coideal-Poisson'' extension of the one-sided Lie coideal $\big(\tilde\fg^\theta, \tau\big)$. In the case of the two-sided coideal structures, the associated Poisson structures are called ``Poisson homogeneous spaces''. We refer to~\cite{BFGP94} and~\cite{EtiKaz95} for details on quantization of such structures.

\subsection{Yangians and twisted Yangians in Drinfeld J presentation}

For any set of primitive elements $x_{i_1},x_{i_2},\ldots,x_{i_m}$ of any associative bi-algebra over $\CC$, set
\begin{gather*}
\{x_{i_1},x_{i_2},\ldots,x_{i_m}\} = \mfrac{1}{m!} \mysum_\pi x_{\pi(i_1)}x_{\pi(i_2)}\cdots x_{\pi(i_m)},
\end{gather*}
where the sum is over all permutations $\pi$ of the set of indices ${i_1},{i_2},\ldots,{i_m}$ and
\begin{gather} \label{xxx}
\lan x_{i_1},x_{i_2},\ldots,x_{i_m} \ran_n = \mfrac{1}{m!} \mysum_\pi x_{\pi(i_1)}\cdots x_{\pi(i_n)} \ot x_{\pi(i_{n+1})} \cdots x_{\pi(i_{m})}.
\end{gather}
such that
\begin{gather} \label{Dxxx}
\Delta(\{x_{i_1},x_{i_2},\ldots,x_{i_m}\}) = \mysum_{n=0}^m {\binom{m}{n}} \lan x_{i_1},x_{i_2},\ldots,x_{i_m} \ran_n,
\end{gather}
where $\binom{m}{n}$ denotes the binomial coefficient. Moreover, for any set of indices $i_1,\dots,i_m$, set
\begin{gather*}
a_{(i_1} a_{i_2} \cdots a_{i_m )} = \mysum_{\sigma\ \mathrm{cyclic}} a_{\sigma(i_1)}a_{\sigma(i_2)}\cdots a_{\sigma(i_m)},
\end{gather*}
where the sum is over the $m$ cyclic permutations of the indices $i_1,\dots,i_m$ $($and not over all $m!$ permutations$)$, e.g.,
$\al_{(ab}^{\qu d} \al_{c)d}^{\qu e} = \al_{ab}^{\sk d} \al_{cd}^{\sk e} + \al_{bc}^{\sk d} \al_{ad}^{\sk e} + \al_{ca}^{\sk d} \al_{bd}^{\sk e} $. The same convention applies when the round bracket encloses indices of a single tensor, e.g., $\mathcal{A}_{(abc)}^{ijk} = \mathcal{A}_{abc}^{ijk} + \mathcal{A}_{bca}^{ijk} + \mathcal{A}_{cab}^{ijk}$.

Next, let us recall the definition of the Yangian as a quantization of the Lie bi-algebra $(\fg[u],\delta)$.

\begin{Theorem}[\cite{Dri85}] \label{T:1} Let $\fg$ be a finite-dimensional complex simple Lie algebra. Fix a~$($non-zero$)$ invariant bilinear form on $\fg$ and a basis~$\{x_a\}$. There is, up to isomorphism, a unique homogeneous quantization $\cY(\fg):=\cU_\hbar(\fg[u])$ of $(\fg[u],\delta)$, with $\delta$ given by~\eqref{deltaJ}. It is topologically generated by elements~$x_a$, $\cJ(x_a)$ with the defining relations:
\begin{gather}
[x_a,x_b] = \alpha_{ab}^{\sk c} x_c, \qquad [x_a,\cJ(x_b)] = \alpha_{ab}^{\sk c} \cJ(x_c), \qquad
\cJ(\lambda x_a + \mu x_b) = \lambda \cJ(x_a) + \mu \cJ(x_b), \label{DT01} \\
[\cJ(x_a),\cJ([x_b,x_c])] + [\cJ(x_b),\cJ([x_c,x_a])] + [\cJ(x_c),\cJ([x_a,x_b])] = \tfrac{1}{4}\hbar^2 \mathcal{A}_{abc}^{ijk} \{x_i,x_j,x_k\}, \label{DT2}\\
[[\cJ(x_a),\cJ(x_b)],\cJ([x_c,x_d])] + [[\cJ(x_c),\cJ(x_d)],\cJ([x_a,x_b])] = \tfrac{1}{4}\hbar^2 \mathcal{B}_{abcd}^{ijk} \{x_i,x_j,\cJ(x_k)\} , \label{DT3}
\end{gather}
for all $x_a\in\fg$ and $\lambda,\mu\in\CC$. Here $\mathcal{A}_{abc}^{ijk} = \alpha_{a}^{\;il}\alpha_{b}^{\;jm}\alpha_{c}^{\;kn}\alpha_{lmn}$ and $\mathcal{B}_{abcd}^{ijk} = \alpha_{cd}^{\sk e} \mathcal{A}_{abe}^{ijk} + \alpha_{ab}^{\sk e} \mathcal{A}_{cde}^{ijk}$. The coproduct is
\begin{gather}
\Delta_\hbar(x_a) = x_a \ot 1 + 1 \ot x_a , \qquad
\Delta_\hbar(\cJ(x_a)) = \cJ(x_a) \ot 1 + 1 \ot \cJ(x_a) + \tfrac{1}{2} \hbar [x_a \ot 1, \Omega_\fg] . \label{dhYJx}
\end{gather}
The grading is $\deg(x_a)=0$, $\deg(\hbar)=1$, $\deg(\cJ(x_a))=1$. The counit is given by $\varepsilon_\hbar(x_a)=\varepsilon_\hbar(\cJ(x_a))=0$. The antipode is
\begin{gather}
S(x_a)=-x_a,\qquad S(\J(x_a))=-\J(x_a)+\tfrac{1}{4}\hbar \fc_\fg x_a. \label{SxJ}
\end{gather}
\end{Theorem}

An outline of the proof can be found in \cite[Section 12.1]{ChaPre94}. Let us make a remark on the Drinfeld terrific relations~\eqref{DT2} and~\eqref{DT3}, which are deformations of the relations~\eqref{DL2} and~\eqref{DL3}, respectively. The coproduct~\eqref{dhYJx} is a solution of the quantization condition satisfying the co-associativity property
\begin{gather*}
(1 \ot \Delta_\hbar) \circ \Delta_\hbar=(\Delta_\hbar \ot 1) \circ \Delta_\hbar .
\end{gather*}
The right-hand sides (rhs) of the terrific relations are tailored so that $\Delta_\hbar \colon \cY(\fg)\to \cY(\fg)\ot\cY(\fg)$ is a homomorphism of algebras. Indeed, consider the co-action on the left-hand side (lhs) of~\eqref{DT2}. The linear terms in $\hbar$ vanish due to the Jacobi identity. What remains are the quadratic terms in $\hbar$, that are cubic and totally symmetric in $x_a$. Hence the rhs of \eqref{DT2} must be of the form $\hbar^2 A_{abc}^{ijk} \{x_i,x_j,x_k\}$ for some set of coefficients $A_{abc}^{ijk}\in\CC$. By comparing the terms on the both sides and using the Jacobi identity one finds $A_{abc}^{ijk}= \tfrac14 \mathcal{A}_{abc}^{ijk}$. Relation \eqref{DT3} is obtained in a similar way.\footnote{Showing \eqref{DT3} is more complicated since it requires a heavy usage of~\eqref{DT2} and Jacobi identity. We were unable to locate the explicit proof of the Drinfeld terrific relations in the mathematical literature available to us.}

We now define twisted Yangians as quantizations of Lie coideals of twisted half-loop algebras for both $\theta\neq \operatorname{id}$ and $\theta=\operatorname{id}$ cases.

\begin{Theorem} \label{T:2}
Let $\big(\fg,\fg^\theta\big)$ be a symmetric pair decomposition of a finite-dimensional simple complex Lie algebra $\fg$ of $\operatorname{rank}(\fg)\geq2$ with respect to the involution $\theta$, such that $\fg^\theta$ is the positive eigenspace of $\theta$. Let $\{X_\al,Y_p\}$ be a symmetric space basis of $\fg$ with respect to $\theta$. There is, up to isomorphism, a unique homogeneous quantization $\cY\big(\fg,\fg^\theta\big)^{\rm tw}:=\cU_\hbar\big(\fg[u]^\theta\big)$ of $\big(\fg[u]^\theta,\tau\big)$, where $\tau$ is given by~\eqref{taunotid}. It is topologically generated by elements $X_\al$, $\B(Y_p)$ with the defining relations:
\begin{gather}
[X_\al,X_\bet]=f^{\sk \ga}_{\al\bet} X_\ga, \qquad [X_\al,\B(Y_p)] = g^{\sk q}_{\al p}\B(Y_q) , \nonumber\\
 \B(a Y_p + b Y_q) = a \B(Y_p) + b \B(Y_q), \label{[X,BY]} \\
[\B(Y_p),\B(Y_q)] + \mysum_\al (\bar\fc_{(\al)})^{-1} w_{pq}^{\sk\al}w_{\al}^{\;rs}[\B(Y_r),\B(Y_s)] = \hbar^2 \mysum_{\la,\mu,\nu} \Lambda_{pq}^{\la\mu\nu} \{X_\la,X_\mu,X_\nu\} , \label{H2}\\
[[\B(Y_p),\B(Y_q)],\B(Y_r)] + 2\fc_\fg^{-1} \mysum_{\al} (\ka_\fm)^{tu} w_{pq}^{\sk\al} g_{r\al}^{\sk s} [[\B(Y_s),\B(Y_t)],\B(Y_u)] \nonumber \\
\qquad{}= \hbar^2 \mysum_{\la,\mu,u} \Upsilon_{pqr}^{\la\mu u} \{X_\la,X_\mu,\B(Y_u)\} , \label{H3}
\end{gather}
for all $X_\al,Y_p\in\fg$ and $a,b\in\CC$. Here
\begin{gather}
\Lambda_{pq}^{\la\mu\nu}= \tfrac{1}{3} \Big( g^{\mu t}_{\sk p} g^{\la u}_{\sk q} + \mysum_{\al} (\bar\fc_{(\al)})^{-1} w_{pq}^{\sk \al} w_{\al}^{\; rs} g^{\mu t}_{\sk r} g^{\la u}_{\sk s} \Big) w_{tu}^{\sk \nu} ,
\label{H2L} \\
\Upsilon_{pqr}^{\la\mu u} = \tfrac{1}{2} \mysum_{\al} \Bigl( w_{st}^{\sk\al} g_{p}^{\;\la s} g_{q}^{\;\mu t} g_{\al r}^{\sk u} + \mysum_\bet w_{pq}^{\sk\al} f_\al^{\;\la\bet} g_{r}^{\;\mu s} g_{\bet s}^{\sk u} + (\la\leftrightarrow\mu)\Bigr) \nonumber\\
\hphantom{\Upsilon_{pqr}^{\la\mu u} =}{}  {} + \fc_\fg^{-1} \mysum_{\al,\ga} (\ka_\fm)^{vx} w_{pq}^{\sk\ga} g_{r\ga}^{\sk y} \Bigl( w_{st}^{\sk\al} g_{y}^{\;\la s} g_{v}^{\;\mu t} g_{\al x}^{\sk u} + \mysum_\bet w_{yv}^{\sk\al} f_\al^{\;\la\bet} g_{x}^{\;\mu s} g_{\bet s}^{\sk u} + (\la\leftrightarrow\mu) \Bigr)  . \label{H3U}
\end{gather}
The coideal structure is given by the co-action $\LHb \colon \cY\big(\fg,\fg^\theta\big)^{\rm tw} \to \cY(\fg) \ot \cY\big(\fg,\fg^\theta\big)^{\rm tw}$ such that
\begin{gather}
\LHb(X_\alpha) = X_\alpha \ot 1 + 1 \ot X_\alpha , \nonumber\\
\LHb(\B(Y_p)) =\varphi(\B(Y_p)) \ot 1 + 1 \ot \B(Y_p) + \hbar [Y_p \ot 1, \Omega_X], \label{dhBY}
\end{gather}
where the embedding $\varphi \colon \cY\big(\fg,\fg^\theta\big)^{\rm tw} \to \cY(\fg)$ is
\begin{gather}
\varphi(\B(Y_p)) =\cJ(Y_p) + \tfrac{1}{4}\hbar [Y_p,C_X] \label{fBY} .
\end{gather}
The grading is $\deg (X_\alpha)=0$, $\deg (\hbar)=1$, $\deg (\B(Y_p))=1$. The counit is $\epsilon_\hbar(X_\al) = \epsilon_\hbar(\B(Y_p)) = 0$ for all non-central $X_\al$. In the case when $\fg^\theta$ is reductive with its centre~$\fk$ generated by $X_z$, the counit $\epsilon_\hbar(X_z)=c$ with $c\in\CC$.
\end{Theorem}

In the case when $\fg^\theta$ is reductive, the one-dimensional representation of $\fg^\theta$ is parametrized by the free parameter $c\in\CC$. This parameter corresponds to the free boundary parameter of a~quantum integrable model with a twisted Yangian as the underlying symmetry algebra. For Lie algebras of type $A$, this parameter also appears in the solutions of the boundary intertwining equation leading to a one-parameter family of the boundary $S$-matrices satisfying the reflection equation~\cite{AACDFR04}. Similar results can also be deduced for types B, C and D.

\begin{Theorem} \label{T:3}
Let $\fg$ be a finite-dimensional simple complex Lie algebra of $\operatorname{rank}(\fg)\geq2$. Fix a~$($non-zero$)$ invariant bilinear form on~$\fg$ and a~basis $\{x_i\}$. There is, up to isomorphism, a~unique homogeneous quantization $\cY(\fg,\fg)^{\rm tw}:=\cU_\hbar(\fg[u^2])$ of $(\fg[u^2],\tau)$, where~$\tau$ is given by~\eqref{tauid}. It is topologically generated by elements~$x_i$, $\G(x_i)$ with the defining relations:
\begin{gather}
[x_a,x_b] = \alpha^{\sk c}_{ab} x_c, \qquad [x_a,\G(x_b)] = \alpha^{\sk c}_{ab} \G(x_c), \qquad \G(\la x_a+\mu x_b) = \la \G(x_a)+\mu \G(x_b),\!\!\! \label{xG1}\\
[\G(x_a),\G([x_b,x_c])] + [\G(x_b),\G([x_c,x_a])] + [\G(x_c),\G([x_a,x_b])] \nonumber\\
\qquad{} = \hbar^2 \Psi_{abc}^{ijk} \{x_i,x_j,\G(x_k)\} + \hbar^4 \big( \Phi_{abc}^{ijk}\{x_i,x_j,x_k\} + \overline{\Phi}_{abc}^{ijklm}\{x_i,x_j,x_k,x_l,x_m\} \big) , \label{H4}
\end{gather}
for all $x_a\in\fg$ and $\la,\mu\in\CC$. Here
\begin{gather}
\Psi_{abc}^{ijk} = \al_{(ab}^{\qu d}\al_{c)r}^{\qu k}\overline{h}_d^{rij} - \al_{dr}^{\sk k}\al_{(ab}^{\qu d}\overline{h}_{c)}^{rij} , \label{coefevenloop1}\\
\overline{\Phi}_{abc}^{\,ijklm} = \tfrac{1}{5}\big(\al_{rs}^{\sk i}\al_{(ab}^{\qu d}h_{c)}^{rjk}h_{d}^{slm}- \Psi_{abc}^{jkr} h_{r}^{ilm} \big) ,
\tag{5.20}\label{coefevenloop2}\\
\Phi_{abc}^{ijk}=\tfrac13\operatorname{Sym}_{ijk}\big(
 \al_{(ab}^{\qu d}\Pi_{c)d}^{ijk}-\Xi_{abc}^{ijk}\big) .
\tag{5.21}\label{coefevenloop3}
\end{gather}
Here
\begin{gather*}
\operatorname{Sym}_{ijk}F^{ijk}=\tfrac16\mysum_{\sigma\in S_3}
F^{\sigma(i)\sigma(j)\sigma(k)}
\end{gather*}
is the normalised symmetrisation of the upper indices, distinct from the
unnormalised cyclic round-bracket convention. The two tensors entering
\eqref{coefevenloop3} are
\begin{gather}
\Pi_{cd}^{ijk} = h_c^{buv}h_d^{frs}\al_{bf}^{i}k_{uv;rs}^{jk}
 +\big(\tfrac12\overline h_c^{juv}\overline h_d^{krs}\ell_{uv;rs}^{i}
 +h_c^{buk}\overline h_d^{frs}m_{b;rs}^{i}\al_{uf}^{j}\nonumber\\
\hphantom{\Pi_{cd}^{ijk} = h_c^{buv}h_d^{frs}\al_{bf}^{i}k_{uv;rs}^{jk}+\big(}
 {}-\overline h_c^{buv}h_d^{frk}m_{f;uv}^{i}\al_{rb}^{j}+(j\leftrightarrow k)\big),\nonumber\\
\Xi_{abc}^{ijk} = \Psi_{abc}^{pqr}h_r^{iuv}n_{pq;uv}^{jk}
 +\big({-}\tfrac16\Psi_{abc}^{pqr}h_r^{suk}\al_{ps}^{i}\al_{uq}^{j}
 +\Psi_{abc}^{pjr}\overline h_r^{kuv}m_{p;uv}^{i}+(j\leftrightarrow k)\big),
\tag{5.22}\label{PiXi}
\end{gather}
\addtocounter{equation}{3}
where
\begin{gather*}
h_a^{\;bcd} = \phi_a^{\;bcd} + 2\psi_a^{\;bcd}, \qquad \overline{h}_a^{\;bcd} = \phi_a^{\;bcd} - \psi_a^{\;bcd}, \qquad
\psi_a^{\;bcd} = \tfrac{1}{12}\big(\al_a^{\;jd} \al_j^{\;bc}+\al_a^{\;jc} \al_j^{\;bd}\big) ,\\
\phi_a^{\;bcd} = \tfrac{1}{24}\fc_\fg^{-1}\mysum_\pi\big(\al_a^{\;jk} \al_{j}^{\;\pi(d)r} \al_{k}^{\;\pi(b)s} \al_{sr}^{\sk \pi(c)}\big) ,
\end{gather*}
and
\begin{align*}
m_{a;bc}^{i}&=\tfrac1{12}\big(\al_{ca}^{r}\al_{br}^{i}+\al_{ba}^{r}\al_{cr}^{i}\big),\\
n_{ab;cd}^{ij}&=\tfrac1{24}\big(\al_{ac}^{i}\al_{bd}^{j}+\al_{ad}^{i}\al_{bc}^{j}\big)
 -\tfrac14\big(\delta_a^i m_{d;bc}^{j}+\delta_b^i m_{d;ac}^{j}\big)\\
&\hphantom{={}}
 +\tfrac1{16}\big(\delta_a^i\al_{bc}^{u}\al_{ud}^{j}+\delta_b^i\al_{ac}^{u}\al_{ud}^{j}\big)
 +(i\leftrightarrow j),\\
k_{ab;cd}^{ij}&=\tfrac18\big(\al_{ac}^{i}\al_{bd}^{j}+\al_{ad}^{i}\al_{bc}^{j}\big)
 -\tfrac14\big(\delta_a^i m_{d;bc}^{j}+\delta_b^i m_{d;ac}^{j}\big)\\
&\hphantom{={}}
 +\tfrac1{16}\big(\delta_a^i\al_{bc}^{u}\al_{ud}^{j}+\delta_b^i\al_{ac}^{u}\al_{ud}^{j}\big)
 +\tfrac12\big(\delta_c^i m_{d;ab}^{j}+\delta_d^i m_{c;ab}^{j}\big)
 +(i\leftrightarrow j),\\
\ell_{ab;cd}^{i}&=-\tfrac12\big(\al_{au}^{i}m_{d;bc}^{u}+\al_{bu}^{i}m_{d;ac}^{u}\big)
 +\tfrac18\big(\al_{bc}^{u}\al_{ud}^{v}\al_{av}^{i}+\al_{ac}^{u}\al_{ud}^{v}\al_{bv}^{i}\big).
\end{align*}
The coideal structure is given by the co-action $\LHb \colon \cY(\fg,\fg)^{\rm tw} \to \cY(\fg) \ot \cY(\fg,\fg)^{\rm tw}$ defined by
\begin{gather}
\LHb(x_a)= x_a \ot 1 + 1 \ot x_a , \label{dhx}\\
\LHb(\G(x_a)) = \varphi(\G(x_a))\ot1+1\ot\G(x_a) + \hbar [\cJ(x_a)\ot1,\Omega_{\fg}] \nonumber\\
\hphantom{\LHb(\G(x_a)) =}{} + \tfrac{1}{4}\hbar^2 \big( [[x_a\ot1,\Omega_\fg],\Omega_\fg] + \fc_\fg^{-1} \al_{a}^{\;bc}[[x_c\ot1,\Omega_\fg],[x_b\ot1,\Omega_\fg]] \big) , \label{dhGx}
\end{gather}
where the embedding $\varphi \colon \cY(\fg,\fg)^{\rm tw} \to \cY(\fg)$ is
\begin{gather}
\varphi(\G(x_a)) = \fc_\fg^{-1} \al_a^{\;bc} [\cJ(x_c),\cJ(x_b)] + \tfrac{1}{4}\hbar [\cJ(x_a),C_\fg]. \label{fGx}
\end{gather}
The grading on $\cY(\fg,\fg)^{\rm tw}$ is $\deg(x_a)=0$, $\deg (\hbar)=1$, $\deg(\G(x_a))=2$. The co-unit is $\epsilon_\hbar(x_i)=\epsilon_\hbar(\G(x_i))=0$.
\end{Theorem}

\begin{Remark} The co-action \eqref{dhGx} can be equivalently written as
\begin{gather}
\LHb(\G(x_a)) = \varphi(\G(x_a))\ot1+1\ot\G(x_a) \nonumber\\
\hphantom{\LHb(\G(x_a)) =}{} + \hbar [\cJ(x_a)\ot1,\Omega_{\fg}] + \hbar^2 \big( h_a^{\;bcd} x_b\ot \{x_c,x_d\} + \overline{h}_a^{\;bcd} \{x_c,x_d\} \ot x_b \big) . \label{dhGx2}
\end{gather}
\end{Remark}

\begin{Remark}The algebras $\cY(\fg,\fg^\theta)^{\rm tw}$ and $\cY(\fg,\fg)^{\rm tw}$ may be considered as flat deformations of the twisted current algebras $\cU(\fg[u]^\theta)$ and $\cU(\fg[u^2])$, respectively. It is clear that $\cY(\fg,\fg^\theta)^{\rm tw} / \hbar \cY(\fg,\fg^\theta)^{\rm tw}$ $ \cong \cU(\fg[u]^\theta)$ and $\cY(\fg,\fg)^{\rm tw} / \hbar \cY(\fg,\fg)^{\rm tw} \cong \cU(\fg[u^2])$. The flatness of the deformations then follows from the Poincar\'e--Birkhoff--Witt (PBW) theorem which is due the fact that both twisted Yangians can be embedded to the Yangian $\cY_\hbar(\fg)$; the PBW theorem for $\cY_\hbar(\fg)$ was demonstrated in~\cite{Lev93}.
\end{Remark}

The proof of Theorems~\ref{T:2} and~\ref{T:3} follows using the same arguments outlined in~\cite[Section~21.1]{ChaPre94}. The uniqueness of co-actions \eqref{dhBY} and (\ref{dhx}), (\ref{dhGx}) is demonstrated in Sections~\ref{sec:7.3} and~\ref{sec:7.4}; in particular, the co-action $\LHb$ is determined uniquely (up to an isomorphism discussed in the paragraph below) by the coideal compatibility identities (\ref{co-co-ass}), (\ref{co-inv}) and the property~(4) of Definition~\ref{D53}. The map~\eqref{fBY} is the MacKay twisted Yangian formula presented in~\cite{DMS01}. The challenging task is to obtain the horrific relations~\eqref{H2}, \eqref{H3} and~\eqref{H4}, which are quantizations of~\eqref{LH2}, \eqref{LH3} and \eqref{LH4}, respectively. A~proof of the first two relations is given in Section~\ref{sec:7.5}. Proving~\eqref{H4} is substantially more difficult. We have given an outline of a proof in Section~\ref{sec:7.6}.

Recall that the Yangian $\cY(\fg)$ has a one-parameter group of automorphisms $\varkappa_c$, with $c\in\CC$, given by $\varkappa_c(x_a)=x_a$ and $\varkappa_c(\J(x_a))=\J(x_a)+\hbar c x_a$, which is compatible with both algebra and Hopf algebra structure. An analogue of this automorphism for the twisted Yangians is a one-parameter family of co-actions $(\varkappa_c\ot \operatorname{id})\circ \LHb$, which for $\cY(\fg,\fg^\theta)^{tw}$ is equivalent to one-parameter family of embeddings $\varkappa_c\circ \varphi$; the later does not apply to $\cY(\fg,\fg)^{tw}$ due to the $\hbar$-order term in~\eqref{dhGx}.

\begin{Remark}The coefficient $\Psi_{abc}^{ijk}$ has the following compact
form for every classical simple $\fg$:
\begin{gather*}
\Psi_{abc}^{ijk} = \tfrac{1}{3}\mathcal{A}_{(abc)}^{ijk}
 + \al_{(ab}^{\qu d} \al_{c)l}^{\sk\; k} \phi_d^{\;lij}
 - \al_{dl}^{\sk k} \al_{(ab}^{\qu d} \phi_{c)}^{\;lij} .
\end{gather*}
Here $\mathcal{A}_{(abc)}^{ijk}$ is understood symmetrised in $i$ and $j$; this
is no restriction, since $\Psi_{abc}^{ijk}$ enters \eqref{H4} only through
$\Psi_{abc}^{ijk}\{x_i,x_j,\G(x_k)\}$, which is symmetric in $i$ and~$j$.
\end{Remark}

\section[Coideal subalgebras of the Yangian Y(sl\_3)]{Coideal subalgebras of the Yangian $\boldsymbol{\cY(\fsl_3)}$} \label{Sec:6}

In this section we present three examples of twisted Yangians $\cY\big(\fg,\fg^\theta\big)^{\rm tw}$ when $\fg=\fsl_3$ and $\fg^\theta=\fso_3$ or $\fg^\theta=\fgl_2$ (both $\theta\neq \operatorname{id}$ cases), or $\fg^\theta=\fsl_3$ ($\theta=\operatorname{id}$ case). For ease of notation we will denote generators of the former two algebras by symbols $\Ch$, $\Ce$, $\Cf$, $\Ck$ and $\CH$, $\CE$, $\CF$. We start by recalling the Chevalley--Serre presentation of the~$\fsl_3$ Lie algebra and its Yangian.

\subsection[The sl\_3 Lie algebra and the Yangian Y(sl\_3)]{The $\boldsymbol{\fsl_3}$ Lie algebra and the Yangian $\boldsymbol{\cY(\fsl_3)}$} Lie algebra $\fsl_3$ in the Chevalley--Serre presentation is generated by elements $e_i$, $f_i$, $h_i$ with $i=1,2$ subject to the defining relations
\begin{gather}
[e_i,f_j]=\delta_{ij}h_i,\qquad [h_i,e_j]=a_{ij} e_j, \qquad [h_i,f_j]=-a_{ij} f_j, \nonumber\\ [e_i,[e_i,e_{i\pm1}]]=0, \qquad [f_i,[f_i,f_{i\pm1}]]=0, \label{Lie:sl3-CS}
\end{gather}
where $a_{ii}=2$ and $a_{12}=-1=a_{21}$ are the matrix elements of the Cartan matrix of $\fsl_3$. The last two relations are called the Serre relations.

Vector space basis of $\fsl_3$ contains eight elements. In addition to the elements given above there are two more root vectors that correspond to the non-simple roots of $\fsl_3$, namely $e_3=[e_1,e_2]$ and $f_3=[f_2,f_1]$. The defining relations for the linear basis are obtained by dropping the Serre relations in \eqref{Lie:sl3-CS} and adding
\begin{gather}
[e_1,e_2]=e_3, \qquad [e_1,e_3]=[e_2,e_3]=0, \qquad [f_1,f_2]=-f_3, \qquad [f_1,f_3]=[f_2,f_3]=0,\nonumber \\
[e_1,f_3]=-f_2, \qquad [e_2,f_3]=f_1, \qquad [f_1,e_3]=e_2, \qu [f_2,e_3]=-e_1,\nonumber\\
[h_i,f_3]=-f_3, \qquad [h_i,e_3]=e_3,\qquad [e_3,f_3]=h_1+h_2.\label{Lie:sl3}
\end{gather}
The quadratic Casimir operator is $C_\fg=\sum_{1\le i\le 3} (e_i f_i+ f_i e_i) + \frac{2}{3} \sum_{1\le i\le j\le 2} h_i h_j$ and $\fc_{\fg}=6$.

\begin{Example}The Yangian $\cY(\fsl_3)$ is the unital associative algebra with sixteen generators $e_i$, $f_i$, $h_j$, $\cJ(e_i)$, $\cJ(f_i)$, $\cJ(h_j)$ with $i=1,2,3$ and $j=1,2$. The defining relations (\ref{DT01}), (\ref{DT2}) are given by (\ref{Lie:sl3-CS}), (\ref{Lie:sl3}) and
\begin{gather} \label{DT2:sl3}
[\cJ(h_1),\cJ(h_2)] = \tfrac{3}{4} \hbar^2 (\{e_3,f_1,f_2\}-\{e_1,e_2,f_3\}).
\end{gather}
The remaining relations defined by \eqref{DT2} are obtained by the adjoint action of degree-$0$ generators on~\eqref{DT2:sl3}. (The same applies to all horrific relations in the examples given below.) The Hopf algebra structure on $\cY(\fsl_3)$ is given by~(\ref{dhYJx}),~(\ref{SxJ}) with $\Omega_{\fg} = \sum_{1\le i\le 3} (e_i\ot f_i + f_i\ot e_i) + \frac{1}{3} \sum_{1\le i\le j\le 2} (h_i\ot h_j + h_j\ot h_i)$.
\end{Example}

\subsection[Twisted Yangian Y(sl\_3,so\_3)\^{}\{tw\}]{Twisted Yangian $\boldsymbol{\cY(\fsl_3,\fso_3)^{\rm tw}}$}

Let involution $\theta$ be defined by
\begin{gather*}
\theta \colon \ e_1 \mapsto -e_2, \qquad f_1 \mapsto - f_2, \qquad h_1 \mapsto h_2 .
\end{gather*}
The action of $\theta$ on the remaining elements of $\fg=\fsl_3$ is deduced from the constraint $\theta^2=\operatorname{id}$. The symmetric space basis is given by $\fg^\theta(=\fh)=\{\Ch=h_1+h_2,\Ce=e_1-e_2,\Cf=f_1-f_2\}$ and $\fm=\{h_1-h_2,e_1+e_2,f_1+f_2,e_3,f_3\}$. The positive eigenspace of $\theta$ forms the orthogonal subalgebra $\fso_3\subset\fsl_3$. We denote the generators of this subalgebra by $\Ch$, $\Ce$, $\Cf$. They generate the degree-0 subalgebra of the twisted Yangian $\cY(\fsl_3,\fso_3)^{\rm tw}$. We denote the degree-1 generators by $\CH$, $\CE$, $\CF$, $\CE_2$, $\CF_2$.

\begin{Example}\label{Y3o}The twisted Yangian $\cY(\fsl_3,\fso_3)^{\rm tw}$ is the unital associative algebra with eight generators $\Ch$,~$\Ce$,~$\Cf$, $\CH$, $\CE$, $\CF$, $\CE_2$, $\CF_2$. The defining relations are the degree-0 Lie relations (of the $\fso_3$ Lie algebra)
\begin{gather*}
[\Ce,\Cf]=\Ch,\qquad [\Ch,\Ce]=\Ce, \qquad [\Ch,\Cf]=-\Cf,
\end{gather*}
degree-1 Lie relations
\begin{alignat*}{6}
& [\Ch,\CE] =\CE, \qquad && [\Ch,\CE_2] =2\CE_2, \qquad && [\CH,\Ce] =3 \CE, \qquad && [\Ce,\CE] =2\CE_2,\qquad && [\Ce,\CE_2]= 0,& \\
& [\Ch,\CF] =-\CF, \qquad && [\Ch,\CF_2] =-2\CF_2, \qquad && [\CH,\Cf] =-3 \CF, \qquad && [\Cf,\CF] =-2\CF_2, \qquad && [\Cf,\CF_2]=0,&\\
& [\Cf,\CE_2] =\CE, \qquad && [\Ce,\CF_2] =-\CF, \qquad && [\Ce,\CF] =\CH, \qquad && [\CE,\Cf]=\CH, \qquad && [\CH,\Ch] =0,&
\end{alignat*}
degree-2 horrific relation
\begin{gather*}
[\CE,\CF]-[\CE_2,\CF_2] = \tfrac{1}{4} \hbar^2 \big( \{\Ch,\Ch,\Ch\} - 3 \{\Ce,\Cf,\Ch\} \big),
\end{gather*}
degree-3 horrific relation
\begin{gather*}
[[\CE,\CF],\CH] = \tfrac{3}{2} \hbar^2 \big(\{\CE_2,\Cf,\Cf\} - \{\CF_2,\Ce,\Ce\} \big) + \tfrac{15}{4} \hbar^2 \big(\{\CE,\Cf,\Ch\} - \{\CF,\Ce,\Ch\} \big).
\end{gather*}
The co-action is given by
\begin{gather*}
\LHb(\Cx) = \varphi(\Cx) \ot 1+1 \ot \Cx , \qquad \LHb(\CY) = \varphi(\CY)\ot 1 +1\ot \CY + \hbar [\varphi_0(\CY)\ot1,\Omega_\fh] ,
\end{gather*}
for all $\Cx\in\{\Ch,\Ce,\Cf\}$ and $\CY\in\{\CE,\CF,\CH,\CE_2,\CF_2\}$ such that
\begin{alignat*}{6}
& \varphi(\Ce) =e_1-e_2,\qquad && \varphi(\Cf) =f_1-f_2, \qquad && \varphi(\Ch)=h_1+h_2, &\\
& \varphi_0(\CE) =e_1+e_2, \qquad && \varphi_0(\CF)=f_1+f_2, \qquad && \varphi_0(\CH)=h_1-h_2, & \\
& \varphi_0(\CE_2)=e_3, \qquad && \varphi_0(\CF_2)=f_3,& &&
\end{alignat*}
and
\begin{alignat*}{3}
& \varphi(\CE) = \cJ(e_1)+\cJ(e_2) + \tfrac{1}{4} \hbar [e_1+e_2,C_X],\qquad && \varphi(\CE_2) = \cJ(e_3) + \tfrac{1}{4} \hbar [e_3,C_X],& \\
& \varphi(\CF) = \cJ(f_1)+\cJ(f_2) + \tfrac{1}{4} \hbar [f_1+f_2,C_X], \qquad && \varphi(\CF_2) = \cJ(f_3) + \tfrac{1}{4} \hbar [f_3,C_X] ,& \\
& \varphi(\CH)= \cJ(h_1)-\cJ(h_2) + \tfrac{1}{4} \hbar [h_1-h_2,C_X].\qquad &&&
\end{alignat*}
Here
\begin{gather*}
C_X = \tfrac{1}{2}\big(\Ce\Cf+ \Cf\Ce+ \Ch^2\big)\in \cU(\fh), \qquad \Omega_\fh = \tfrac12 (\varphi\ot \operatorname{id})\circ\big(\Ce\ot\Cf + \Cf\ot\Ce + \Ch\ot\Ch\big)\in\fg\ot\fh .
\end{gather*}
The co-unit is $\ep(\mathsf x)= \ep(\mathsf Y)=0$.
\end{Example}

\subsection[Twisted Yangian Y(sl\_3,gl\_2)\^{}\{tw\}]{Twisted Yangian $\boldsymbol{\cY(\fsl_3,\fgl_2)^{\rm tw}}$}

Let involution $\theta$ be defined by
\begin{gather*}
\theta \colon \ e_1\mapsto e_1, \qquad f_1\mapsto f_1,\qquad e_2\mapsto -e_2,\qquad f_2\mapsto -f_2, \qquad h_i\mapsto h_i,
\end{gather*}
In this case $\fg^\theta = \{\Ch=h_1,\Ck=2h_2+h_1,\Ce=e_1,\Cf=f_1\}\sim \fgl_2 $ and $\fm = \{e_2,f_2,e_3,f_3\}$. We denote the corresponding degree-1 generators by $\CE_2$, $\CF_2$, $\CE_3$, $\CF_3$.

\begin{Example} \label{Y3g}
The twisted Yangian $\cY(\fsl_3,\fgl_2)^{\rm tw}$ is the unital associative algebra with eight generators $\Ch$, $\Ce$, $\Cf$, $\Ck$, $\CE_2$, $\CF_2$, $\CE_3$, $\CF_3$. The defining relations are the degree-0 Lie relations (of the~$\fgl_2$ Lie algebra)\footnote{The standard $\fgl_2$ basis $\{e_{ij}\}_{i,j=1,2}$ with the defining relations $[e_{ij},e_{kl}]=\delta_{kj}e_{il}-\delta_{il}e_{kj}$ is obtained by setting $e_{11}=-(\Ch+\Ck)/2$, $e_{22}=(\Ch-\Ck)/2$, $e_{12}=\Cf$ and $e_{21}=\Ce$.}
\begin{gather*}
[\Ce,\Cf]=\Ch,\qquad [\Ch,\Ce]=2\Ce, \qquad [\Ch,\Cf]=-2\Cf, \qquad [\Ce,\Ck]= [\Cf,\Ck]= [\Ch,\Ck]=0,
\end{gather*}
degree-1 Lie relations
\begin{alignat*}{4}
& [\Ch,\CE_2] =-\CE_2, \qquad && [\Ch,\CE_3] =\CE_3, \qquad && [\Ck, \CE_i] = 3\CE_i,& \\
& [\Ch,\CF_2] =\CF_2, \qquad && [\Ch,\CF_3] =-\CF_3, \qquad && [\Ck, \CF_i] = -3\CF_i,& \\
& [\Ce,\CE_2] =\CE_3, \qquad && [\Ce,\CF_3] =-\CF_2, \qquad && [\Ce,\CF_2] = [\Ce,\CE_3] =0,& \\
& [\Cf,\CF_2] =-\CF_3, \qquad && [\Cf,\CE_3] =\CE_2, \qquad && [\Cf,\CE_2] =[\Cf,\CF_3]=0, &
\end{alignat*}
degree-2 horrific relations
\begin{gather*}
[\CE_2,\CE_3] = 0, \qquad [\CF_2,\CF_3] = 0,
\end{gather*}
degree-3 horrific relations
\begin{gather*}
[\CE_2,[\CE_2,\CF_3]] = 2 \hbar^2\{\CE_2,\Cf,\Ck\}, \qquad [\CF_2,[\CE_3,\CF_2]] = -2 \hbar^2\{\CF_2,\Ce,\Ck\}.
\end{gather*}
The co-action given by, for $\Cx\in\{\Ch,\Ce,\Cf,\Ck\}$,
\begin{gather*}
\LHb(\Cx) = \varphi(\Cx) \ot 1+1 \ot \Cx , \qquad \varphi(\Ch) = h_1,\\ \varphi(\Ck) = 2h_2+h_1, \qquad \varphi(\Ce) = e_1, \qquad \varphi(\Cf) = f_1
\end{gather*}
and, for $i=2,3$,
\begin{alignat*}{3}
& \LHb(\CE_i) = \varphi(\CE_i) \ot 1+1 \ot \CE_i + \hbar [e_i\ot1,\Omega_\fh] , \qquad && \varphi(\CE_i)=\cJ(e_i) + \tfrac{1}{4}\hbar [e_i,C_X], & \\
& \LHb(\CF_i) = \varphi(\CF_i) \ot 1+1 \ot \CF_i + \hbar [f_i\ot1,\Omega_\fh] , \qquad && \varphi(\CF_i)=\cJ(f_i) + \tfrac{1}{4}\hbar [f_i,C_X].&
\end{alignat*}
Here
\begin{gather*}
C_X = \Ce\Cf + \Cf\Ce + \tfrac12\Ch^2+\tfrac16\Ck^2\in \cU(\fh),\\ \Omega_\fh = (\varphi\ot \operatorname{id})\circ\big(\Ce\ot\Cf+\Cf\ot\Ce+\tfrac12 \Ch\ot\Ch+\tfrac16 \Ck\ot\Ck\big)\in\fg\ot\fh.
\end{gather*}
The co-unit is $\ep(\Ck)=c\in\CC$, $\ep(\Cx)=\ep(\CY)=0$ for all $\Cx\in\{\Ce,\Cf,\Ch\}$ and $\CY\in\{\CE_i,\CF_i\}$ with $i=2,3$.
\end{Example}

\subsection[Twisted Yangian Y(sl\_3,sl\_3)\^{}\{tw\}]{Twisted Yangian $\boldsymbol{\cY(\fsl_3,\fsl_3)^{\rm tw}}$}

In this case the involution is trivial, $\theta=\operatorname{id}$, hence $\fg^\theta = \fsl_3$ and $\fm=\varnothing$.

\begin{Example}\label{Y30}The twisted Yangian $\cY(\fsl_3,\fsl_3)^{\rm tw}$ is the unital associative algebra with sixteen generators $e_i$, $f_i$, $h_j$, $\G(e_i)$, $\G(f_i)$, $\G(h_j)$ with $i=1,2,3$ and $j=1,2$, obeying the standard $\fsl_3$ Lie algebra relations of the Cartan--Chevalley presentation and the standard degree-2 Lie relations~\eqref{xG1} and the following degree-4 horrific relation
\begin{gather*}
[\G(h_1),\G(h_2)] =\hbar^2 \big(\{f_1,f_2,\G(e_3)\}-\{e_1,e_2,\G(f_3)\}) - \tfrac{1}{2}\hbar^4 (\{f_1,f_2,e_3\}-\{e_1,e_2,f_3\}\big) \\
\hphantom{[\G(h_1),\G(h_2)] =}{} + \hbar^2 \big(\{e_3,f_2,\G(f_1)\}+\{e_3,f_1,\G(f_2)\}-\{f_3,e_2,\G(e_1)\}-\{f_3,e_1,\G(e_2)\}\big)\\
\hphantom{[\G(h_1),\G(h_2)] =}{} + \tfrac{1}{2}\hbar^2 \big( \{h_1,f_2,\G(e_2)\}-\{h_1,e_2,\G(f_2)\}- \{h_2,f_1,\G(e_1)\}+ \{h_2,e_1,\G(f_1)\}\big)
\\
\hphantom{[\G(h_1),\G(h_2)] =}{} + \tfrac{1}{4}\hbar^2 \big( \{h_1,e_1,\G(f_1)\}- \{h_1,f_1,\G(e_1)\}- \{h_2,e_2,\G(f_2)\}+ \{h_2,f_2,\G(e_2)\} \\
\hphantom{[\G(h_1),\G(h_2)] =}{} +\{h_1{-}h_2,f_3,\G(e_3)\}-\{h_1{-}h_2,e_3,\G(f_3)\}\big)
\\
\hphantom{[\G(h_1),\G(h_2)] =}{} + \tfrac{1}{4}\hbar^4 \big(\{e_1,e_1,e_2,f_1,f_3\}+\{e_1,e_2,e_2,f_2,f_3\}+\{e_1,e_2,e_3,f_3,f_3\} \\
\hphantom{[\G(h_1),\G(h_2)] =}{} -\{f_1,f_1,f_2,e_1,e_3\}-\{f_1,f_2,f_2,e_2,e_3\}-\{f_1,f_2,f_3,e_3,e_3\}\big) \\
\hphantom{[\G(h_1),\G(h_2)] =}{} +\tfrac{1}{12}\hbar^4 \big(\{e_1,e_2,h_1,h_1,f_3\}+\{e_1,e_2,h_1,h_2,f_3\}+\{e_1,e_2,h_2,h_2,f_3\} \\
\hphantom{[\G(h_1),\G(h_2)] =}{} -\{f_1,f_2,h_1,h_1,e_3\}-\{f_1,f_2,h_1,h_2,e_3\}-\{f_1,f_2,h_2,h_2,e_3\}\big) .
\end{gather*}

The co-action given by (\ref{dhx}), (\ref{dhGx}) with $\varphi(\G(x_i)) = \cK(x_i)+\tfrac{1}{4}\hbar [\cJ(x_i),C_\fg]$ and with $\cK(x_i)=\fc_\fg^{-1}\alpha_i^{\;kj}[\cJ(x_j),\cJ(x_k)]$, where\footnote{The non-zero structure constants $\alpha_i^{\;kj}$ for $\fsl_3$ in the Cartan--Chevalley basis can be read from here.}
\begin{gather*}
\cK(e_1) = \tfrac{1}{3}([\cJ(e_3),\cJ(f_2)]+[\cJ(h_1),\cJ(e_1)]) , \\ \cK(f_1) = \tfrac{1}{3}([\cJ(e_2),\cJ(f_3)] -[\cJ(h_1),\cJ(f_1)]) , \\
\cK(e_2) = \tfrac{1}{3}([\cJ(f_1),\cJ(e_3)]+[\cJ(h_2),\cJ(e_2)]) , \\ \cK(f_2)= \tfrac{1}{3}([\cJ(f_3),\cJ(e_1)]-[\cJ(h_2),\cJ(f_2)]) , \\
\cK(e_3) = \tfrac{1}{3}([\cJ(e_1),\cJ(e_2)]+[\cJ(h_1+h_2),\cJ(e_3)]) , \\ \cK(f_3)= \tfrac{1}{3}([\cJ(f_2),\cJ(f_1)]-[\cJ(h_1+h_2),\cJ(f_3)]) , \\
\cK(h_1) = \tfrac{1}{3}\big(2[\cJ(e_1),\cJ(f_1)]-[\cJ(e_2),\cJ(f_2)]+[\cJ(e_3),\cJ(f_3)]\big) , \\
\cK(h_2) = \tfrac{1}{3}\big({-}[\cJ(e_1),\cJ(f_1)]+2[\cJ(e_2),\cJ(f_2)]+[\cJ(e_3),\cJ(f_3)]\big) .
\end{gather*}
The co-unit is $\epsilon(x_i)= \epsilon(\G(x_i))=0$.
\end{Example}

\section{Proofs} \label{Sec:7}

In the remaining part of the paper we provide proofs omitted in the previous sections.

\subsection{A proof of Proposition \ref{P:31}} \label{sec:7.1}

We follow the arguments outlined in \cite[Section 12.1]{ChaPre94} and fill in the gaps. Recall that given a~simple Lie algebra $\fg$ the half-loop Lie algebra $\cL^+$ is generated by the elements~$x^{(k)}_a$ with $k\ge0$ and $a=1,\ldots,\dim(\fg)$, and satisfying~\eqref{Lp[x,x]}. Let $\wt{\cL}^+$ denote the algebra generated by the elements~$x_a$ and~$J(x_a)$ satisfying (\ref{DL01})--(\ref{DL3}). The map \mbox{$\varphi \colon \wt\cL^+ \to \cL^+$} given by $ x_a \mapsto x_a^{(0)}$, $J(x_a) \mapsto x_a^{(1)}$ is an algebra homomorphism. It is a direct computation to check that the image of (\ref{DL01})--(\ref{DL3}) holds in~$\cL^+$.

To prove that $\varphi$ is surjective we need to show that elements $x^{(0)}_a$ and $x^{(1)}_a$ generate the whole~$\cL^+$. Indeed, since $\fg$ is simple, we have $[\fg,\fg]=\fg$, or equivalently $\big[\cL^{(1)},\cL^{(1)}\big]=\cL^{(2)}$, where $\cL^{(k)} = \operatorname{span}\big\{x_a^{(k)}\big\}$. By the same arguments we must have, for all $k>0$, that $\big[\cL^{(1)},\cL^{(k)}\big]=\cL^{(k+1)}$. This shows that $\cL^+$ is generated by $\cL^{(0)}$ and $\cL^{(1)}$. Hence it only remains to show that~$\varphi$ is injective.

We know that $\cL^+ = \bigoplus_{k\ge 0} \cL^{(k)}$, with $\dim\big(\cL^{(k)}\big)=\dim(\fg)$, as a vector space. We need to show that an analogous statement holds for $\wt\cL^+$. Denote $J^{(0)}(x_a)=x_a$, $J^{(1)}(x_a)=J(x_a)$ and define recursively
\begin{gather}
J^{(k)}(x_a) = \fc_\fg^{-1} \al_a^{bc} \big[J^{(1)}(x_c), J^{(k-1)}(x_b)\big] \qquad \text{for all}\quad k\ge 1. \label{71:0}
\end{gather}
We need to show that
\begin{gather}
\big[J^{(\ell)}(x_a), J^{(k-\ell)}(x_b)\big] = \al_{ab}^{\sk c} J^{(k)}(x_c) \qquad \text{for all}\quad \ell \le k , \label{71:1}
\end{gather}
so that $\wt\cL^+ = \bigoplus_{k\ge 0} \wt\cL^{(k)}$ with $\wt\cL^{(k)} = \operatorname{span}\{J^{(k)}(x_a)\}$ and $\dim(\wt\cL^{(k)})=\dim(\fg)$. We first demonstrate that
\begin{gather}
\big[J^{(k-\ell)}(x_a), J^{(\ell)}(x_b)\big] = \big[J^{(k-\ell')}(x_a), J^{(\ell')}(x_b)\big] \qquad \text{for all}\quad \ell<\ell' \le k.\label{71:2}
\end{gather}
We proceed by induction on $k$, the result being clear when $k=0,1$, by \eqref{DL01}. For $k=2$ there are two inequivalent identities that need to be shown: with $\ell=0$, $\ell'=2$ and with $\ell=0$, $\ell'=1$. The first identity follows by~\eqref{71:0} and co-Jacobi identity:
\begin{gather*}
\big[J^{(2)}(x_a),J^{(0)}(x_b)\big] = \fc_\fg^{-1} \al_{a}^{\; dc} \big[\big[J^{(1)}(x_c),J^{(1)}(x_d)\big],J^{(0)}(x_b)\big]\\
\hphantom{\big[J^{(2)}(x_a),J^{(0)}(x_b)\big]}{} = \fc_\fg^{-1} \al_{a}^{\; dc} \big( \al_{cb}^{\sk e}\big[J^{(1)}(x_e),J^{(1)}(x_d)\big] + \al_{db}^{\sk e} \big[J^{(1)}(x_c),J^{(1)}(x_e)\big] \big)\\
\hphantom{\big[J^{(2)}(x_a),J^{(0)}(x_b)\big]}{} = \al_{ab}^{\sk c} J^{(2)}(x_c) .
\end{gather*}
In a similar way we obtain that $\big[J^{(0)}(x_a),J^{(2)}(x_b)\big] = \al_{ab}^{\sk c} J^{(2)}(x_c)$ thus yielding the required identity. To obtain the second identity we need to contract \eqref{DL2} with $\al_{d}^{\;cb}$ giving
\begin{gather*}
\fc_\fg \big[J^{(1)}(x_a),J^{(1)}(x_d)\big] + \al_{d}^{\;cb} \al_{ca}^{\sk e} \big[J^{(1)}(x_b),J^{(1)}(x_e)\big] + \al_{d}^{\;cb} \al_{ab}^{\sk e} \big[J^{(1)}(x_c),J^{(1)}(x_e)\big] = 0.
\end{gather*}
Now rename the indices $b\to c$, $c\to e$, $e\to b$ in the third term and use co-Jacobi identity together with~\eqref{71:0}. This gives $\big[J^{(1)}(x_a),J^{(1)}(x_d)\big] = \al_{ad}^{\sk c} J^{(2)}(x_c) = \big[J^{(0)}(x_a),J^{(2)}(x_d)\big]$. Next, assuming inductively~\eqref{71:2}, we have
\begin{gather*}
\big[J^{(s)}(x_c), \big[J^{(r-\ell)}(x_a), J^{(\ell)}(x_a)\big]\big]\\
 \qquad{}= \al_{ca}^{\sk b}\big(\big[J^{(r-\ell+s)}(x_b), J^{(\ell)}(x_a)\big] + \big[J^{(r-\ell)}(x_a), J^{(\ell+s)}(x_b)\big]\big) = 0
\end{gather*}
for all $\ell\le r\le k$, $r-\ell+s\le k$ and $\ell+s\le k$ giving
\begin{gather*}
\big[J^{(r-\ell+s)}(x_b), J^{(\ell)}(x_a)\big] = \big[J^{(\ell+s)}(x_b),J^{(r-\ell)}(x_a)\big].
\end{gather*}
Taking all the allowed $\ell$, $r$, $s$ such that $r+s\le k+1$ we obtain all the necessary identities thus completing the induction.

We are now ready to prove~\eqref{71:1}. We use induction on $k$. The cases with $k\le2$ are discussed above. Thus, by induction hypothesis, we have
\begin{gather*}
\big[J^{(1)}(x_d),\big[J^{(k-\ell)}(x_a), J^{(\ell)}(x_b)\big]\big] = \al_{ab}^{\sk c} \big[J^{(1)}(x_d),J^{(k)}(x_c)\big].
\end{gather*}
Now contract the lhs with $\al_e^{\;ad}$, use \eqref{71:0}, \eqref{71:2} and co-Jacobi identity:
\begin{gather*}
\al_e^{\;ad}\big[J^{(1)}(x_d),\big[J^{(k-\ell)}(x_a), J^{(\ell)}(x_b)\big]\big] \\
\qquad{} = \fc_\fg \big[J^{(k-\ell+1)}(x_e), J^{(\ell)}(x_b)\big] + \al_e^{\;ad} \al_{db}^{\sk f} \big[J^{(k-\ell)}(x_a), J^{(\ell+1)}(x_f)\big]\\
\qquad{} = \fc_\fg \big[J^{(k-\ell+1)}(x_e), J^{(\ell)}(x_b)\big] - \tfrac12 \fc_\fg \al_{eb}^{\sk d} J^{(k+1)}(x_d).
\end{gather*}
By doing the same for the rhs we get
\begin{gather*}
\al_e^{\;ad}\al_{ab}^{\sk c} \big[J^{(1)}(x_d),J^{(k)}(x_c)\big] = \tfrac12 \fc_\fg \al_{eb}^{\sk a} J^{(k+1)}(x_a)
\end{gather*}
yielding $\big[J^{(k-\ell+1)}(x_e), J^{(\ell)}(x_b)\big] = \al_{eb}^{\sk d} J^{(k+1)}(x_d)$, which completes the induction. Therefore $\wt\cL^+ = \bigoplus \wt\cL^{(k)}$. By setting $\varphi \colon J^{(k)}(x_a) \mapsto x^{(k)}_a$ we obtain a bijection of vector spaces, hence $\wt\cL^+ \cong \cL^+$.

\subsection{A proof of Proposition~\ref{P:32}} \label{sec:7.2}

The idea behind the proof is very similar to that of the prove above. Recall that $\cH^+$ is the subalgebra of $\cL^+$ generated by elements $X^{(2k)}_\al$ and $Y^{(2k+1)}_p$ with $k\ge0$ and satisfying~\eqref{LpXY}. Let~$\wt{\cH}^+$ denote the algebra generated by the elements $X_\al$ and $B(Y_p)$ satisfying (\ref{L[X,BY]})--(\ref{LH3}). The map $\psi \colon \wt{\cH}^+ \to \cH^+$ given by $X_\al \mapsto X_\al^{(0)}$, $B(Y_p) \mapsto Y_p^{(1)}$ is an algebra homomorphism. It is easy to see that the image of~\eqref{L[X,BY]} holds in~$\cH^+$. To see that the same is true for \eqref{LH2}~and~\eqref{LH3} we additionally need to use~\eqref{wwid} and~\eqref{ggid}.

To prove that $\psi$ is surjective we need to show that elements $X^{(0)}_\al$ and $Y^{(1)}_p$, where $X_\al$ runs over a basis for $\fh$ and $Y_p$ runs over a basis for $\fm$, generate the whole $\cH^+$. Indeed, since $\fg = \fm \op \fh$ is simple, we must have $[\fm,\fm]=\fh$ and $[\fm,\fh]=\fm$, or equivalently $\big[\cH^{(1)},\cH^{(1)}\big]=\cH^{(2)}$ and $\big[\cH^{(1)},\cH^{(2)}\big]=\cH^{(3)}$, where $\cH^{(2k)}=\operatorname{span}\big\{X_\al^{(2k)}\big\}$ and $\cH^{(2k+1)}=\operatorname{span}\big\{Y_p^{(2k+1)}\big\}$. By the same arguments we must have, for all $k>0$, that $\big[\cH^{(1)},\cH^{(k)}\big]= \cH^{(k+1)}$. This shows that $\cH^+$ is generated by $\cH^{(0)}$ and $\cH^{(1)}$. It remains to show that $\psi$ is injective.

We know that $\cH^+ = \bigoplus_{k\ge 0} \cH^{(k)}$, with $\dim\big(\cH^{(2k)}\big)=\dim(\fh)$ and $\dim\big(\cH^{(2k+1)}\big)=\dim(\fm)$, as a vector space. We need to show that an analogous statement holds for $\wt\cH^+$. Let $B^{(0)}(X_\al)=X_\al$, $B^{(1)}(Y_p)=B(Y_p)$ and define recursively, for $k\ge0$,
\begin{gather}
B^{(2k+2)}(X_\al) = (\bar{\fc}_{(\al)})^{-1} w_{\al}^{\;qp}\big[B^{(1)}(Y_p),B^{(2k+1)}(Y_q)\big] , \label{72:BX} \\
B^{(2k+1)}(Y_p) = 2\fc_\fg^{-1} \mysum_{\al} g_{p}^{\;\al q} \big[B^{(1)}(Y_q),B^{(2k)}(X_\al)\big]. \label{72:BY}	
\end{gather}
We need to show that, for all $\ell\le k$,
\begin{gather}
\big[B^{(2k-2\ell)}(X_\alpha),B^{(2\ell)}(X_\beta)\big]= {f}^{\sk \gamma}_{\alpha\beta} B^{(2k)}(X_\gamma), \label{B:XX} \\
\big[B^{(2k-2\ell)}(X_\al),B^{(2\ell+1)}(Y_p)\big]= g^{\sk q}_{\al p} B^{(2k+1)}(Y_q), \label{B:XY} \\
\big[B^{(2k-2\ell+1)}(Y_p),B^{(2\ell+1)}(Y_q)\big] = \mysum_{\alpha}w^{\sk \alpha}_{pq} B^{(2k+2)}(X_\alpha), \label{B:YY}
\end{gather}
so that $\wt\cH^+ = \bigoplus_{k\ge 0} \wt\cH^{(k)}$, with $\wt\cH^{(2k)}=\operatorname{span}\big\{B^{(2k)}(X_\al)\big\}$ and $\wt\cH^{(2k+1)}=\operatorname{span}\big\{B^{(2k+1)}(Y_p)\big\}$, as a vector space.

We start by showing (\ref{B:XX})--(\ref{B:YY}) for small $k$ that follow from (\ref{L[X,BY]})--(\ref{LH3}). First, notice that~\eqref{72:BX} implies that~\eqref{LH2} is equivalent to $\big[B^{(1)}(Y_p),B^{(1)}(Y_q)\big] = \sum_\al w_{pq}^{\sk \al} B^{(2)}(X_\al)$. Also, notice that, by~\eqref{inv:Y}, \eqref{wwid} and~\eqref{72:BX},
\begin{gather*}
\begin{split}
& \mysum_\al w_{pq}^{\sk\al} \big[B^{(2)}(X_\al),B^{(1)}(Y_r)\big]= \mysum_{\al} (\bar\fc_{(\al)})^{-1} w_{pq}^{\sk\al} w_{\al}^{\;st} \big[\big[B^{(1)}(Y_t),B^{(1)}(Y_s)],B^{(1)}(Y_r)\big] \\
& \hphantom{\mysum_\al w_{pq}^{\sk\al} \big[B^{(2)}(X_\al),B^{(1)}(Y_r)\big]}{} = \big[\big[B^{(1)}(Y_p),B^{(1)}(Y_q)\big],B^{(1)}(Y_r)\big] .
\end{split}
\end{gather*}
Then, by \eqref{inv:Y}, \eqref{ggid} and (\ref{72:BX}), (\ref{72:BY}),
\begin{gather*}
\mysum_\al w_{pq}^{\sk\al} g_{r\al}^{\sk s} B^{(3)}(Y_s) = \mysum_{\al,\bet} 2 (\fc_\fg)^{-1} w_{pq}^{\sk\al} g_{r\al}^{\sk s} g_s^{\;u\bet}\big[B^{(2)}(X_\bet),B^{(1)}(Y_u)\big] \\
\hphantom{\mysum_\al w_{pq}^{\sk\al} g_{r\al}^{\sk s} B^{(3)}(Y_s)}{}
= \mysum_{\al,\bet} 2 (\ka_\fm)^{tu}(\fc_{\fg} \bar\fc_{(\bet)})^{-1} w_{pq}^{\sk\al} g_{r\al}^{\sk\, s} w_{st}^{\sk\bet} w_{\bet}^{\;vx} \big[\big[B^{(1)}(Y_x),B^{(1)}(Y_v)\big],B^{(1)}(Y_u)\big] \\
\hphantom{\mysum_\al w_{pq}^{\sk\al} g_{r\al}^{\sk s} B^{(3)}(Y_s)}{} = 2\fc_\fg^{-1}(\ka_\fm)^{tu} w_{pq}^{\sk\al} g_{r\al}^{\sk\, s}\big[\big[B^{(1)}(Y_s),B^{(1)}(Y_t)\big],B^{(1)}(Y_u)\big] .
\end{gather*}
By combining the expressions above with \eqref{LH3} we obtain $\big[B^{(2)}(X_\al),B^{(1)}(Y_r)\big] = g_{\al r}^{\sk\, s} B^{(3)}(Y_s)$. Consequently, we have that
\begin{gather*}
w_\bet^{\;rp} \big[B^{(1)}(Y_p),\big[B^{(2)}(X_\al),B^{(1)}(Y_r)\big]\big]= w_\bet^{\;rp} g_{\al r}^{\sk\, s} \big[ B^{(1)}(Y_p), B^{(3)}(Y_s)\big]
\end{gather*}
giving
\begin{gather*}
w_\bet^{\;rp} g_{p\al}^{\sk q} \big[B^{(3)}(Y_q),B^{(1)}(Y_r)\big]+ \mysum_\ga w_\bet^{\;rp} w_{pr}^{\sk \ga} \big[B^{(2)}(X_\al),B^{(2)}(X_\ga)\big]\\
 \qquad{}= w_\bet^{\;pr} g_{\al p}^{\sk\, q} \big[ B^{(1)}(Y_r), B^{(3)}(Y_q)\big].
\end{gather*}
Using \eqref{MixJac}, \eqref{wwid}, (\ref{72:BX}), (\ref{72:BY}) and antisymmetry in $\al$ and $\beta$ we obtain
\begin{gather*}
\bar\fc_{(\bet)} \big[B^{(2)}(X_\al),B^{(2)}(X_\bet)\big]= \big( w_\bet^{\;pr} g_{\al p}^{\sk\, q} - w_\al^{\;pr} g_{\bet p}^{\sk q}\big) \big[ B^{(1)}(Y_r), B^{(3)}(Y_q)\big] = \bar\fc_{(\bet)} f_{\al\bet}^{\sk \ga} B^{(4)}(X_\ga) ,
\end{gather*}
which gives $\big[B^{(2)}(X_\al),B^{(2)}(X_\bet)\big]= f_{\al\bet}^{\sk \ga} B^{(4)}(X_\ga)$. Then, a direct computation using (\ref{72:BX}), (\ref{72:BY}) yields $\big[B^{(0)}(X_\al),B^{(2)}(X_{\bet})\big]= f_{\al\bet}^{\sk\ga} B^{(2)}(X_\ga)$ and $\big[B^{(0)}(X_\al),B^{(3)}(Y_r)\big] = g_{\al r}^{\sk\, s} B^{(3)}(Y_s)$ and finally $\big[B^{(0)}(X_\al),B^{(4)}(X_\bet)\big] = f_{\al\bet}^{\sk\ga} B^{(4)}(X_\ga)$.

Next we show the following identities that will be necessary in proving~\eqref{B:XY} and~\eqref{B:YY}:
\begin{gather}
\big[B^{(2k-2\ell)}(X_\al),B^{(2\ell+1)}(Y_p)\big] = \big[B^{(2k-2\ell')}(X_\al),B^{(2\ell'+1)}(Y_p)\big] , \label{BB:XY} \\
\big[B^{(2k-2\ell+1)}(Y_p),B^{(2\ell+1)}(Y_q)\big] = \big[B^{(2k-2\ell'+1)}(Y_p),B^{(2\ell'+1)}(Y_q)\big] . \label{BB:YY}
\end{gather}
We proceed by induction on $k$, the result being clear when $k=0$ and when $k=1$ for~\eqref{BB:XY}. Assuming induction we have
\begin{gather*}
\big[B^{(2s)}(X_\al),\big[B^{(2r-2\ell+1)}(Y_p),B^{(2\ell+1)}(Y_p)\big]\big] \\
\qquad{} = g_{\al p}^{\sk\; q} \big( \big[B^{(2r-2\ell+2s+1)}(Y_q),B^{(2\ell+1)}(Y_p)\big] + \big[B^{(2r-2\ell+1)}(Y_p),B^{(2\ell+2s+1)}(Y_q)\big] \big) = 0
\end{gather*}
for all $\ell\le r\le k$, $r-\ell+s\le k$ and $\ell+s\le k$, therefore we must have
\begin{gather*}
\big[B^{(2r+2s-2\ell+1)}(Y_q),B^{(2\ell+1)}(Y_p)\big]=\big[B^{(2s+2\ell+1)}(Y_q),B^{(2r-2\ell+1)}(Y_p)\big].
 \end{gather*}
Taking all the allowed $\ell$, $r$, $s$ such that $r+s\le k+1$ we obtain all the necessary identities thus completing the induction. Then, by~\eqref{72:BX}, it follows that, for all $\ell\le k$,
\begin{gather}
B^{(2k+2)}(X_\al) = (\bar{\fc}_{(\al)})^{-1} w_{\al}^{\;qp}\big[B^{(2\ell+1)}(Y_p),B^{(2k-2\ell+1)}(Y_q)\big]. \label{BX2}
\end{gather}

Relation \eqref{BB:XY} is proved similarly, using
\begin{gather*}
\big[B^{(2s+1)}(Y_q),\big[B^{(2r-2\ell+1)}(Y_p),B^{(2\ell+1)}(Y_p)\big]\big] \\
\qquad{} = \mysum_\al w_{qp}^{\sk \al}\big( \big[B^{(2r-2\ell+2s+2)}(X_\al),B^{(2\ell+1)}(Y_p)\big] + \big[B^{(2r-2\ell+1)}(Y_p),B^{(2\ell+2s+2)}(X_\al)\big] \big) = 0
\end{gather*}
and implies that, for all $\ell\le k$,
\begin{gather}
B^{(2k+1)}(Y_p) = 2 \fc_\fg^{-1} \mysum_{\al} g_{p}^{\;\al q} \big[B^{(2\ell+1)}(Y_q),B^{(2k-2\ell)}(X_\al)\big]. \label{BY2}
\end{gather}

We also need to show that, for all $0\le\ell\le k$ and all $X_\al$ in a simple ideal of $\fh$ (so that $\fc_{(\al)}\neq0$),
\begin{gather}
B^{(2k+2)}(X_\al) = (\fc_{(\al)} )^{-1} f_{\al}^{\;\ga\bet}\big[B^{(2k-2\ell)}(X_\bet),B^{(2\ell+2)}(X_\ga)\big] . \label{BX3}
\end{gather}
Again, we proceed by induction on $k$. The base of induction was already shown above. Hence, assuming inductively~\eqref{BX3}, we obtain
\begin{gather*}
 f_{\al}^{\;\ga\bet}\big[B^{(2k-2\ell)}(X_\bet),B^{(2\ell+2)}(X_\ga)\big]\\
 \qquad{} = (\bar{\fc}_{(\ga)} )^{-1} f_{\al}^{\;\ga\bet} w_{\ga}^{\;qp}\big[B^{(2k-2\ell)}(X_\bet),\big[B^{(2\ell'+1)}(Y_p),B^{(2\ell-2\ell'+1)}(Y_q)\big]\big]\\
\qquad{} =(\bar{\fc}_{(\ga)} )^{-1} f_{\al}^{\;\ga\bet} w_{\ga}^{\;qp} \big( g_{\bet p}^{\sk r}\big[B^{(2k-2\ell+2\ell'+1)}(Y_r),B^{(2\ell-2\ell'+1)}(Y_q)\big]\\
\qquad\quad{} +g_{\bet q}^{\sk r}\big[B^{(2\ell'+1)}(Y_p),B^{(2k-2\ell'+1)}(Y_r)\big] \big)\\
\qquad{} =(\bar{\fc}_{(\al)} )^{-1} \fc_{(\al)} w_{\al}^{\;rp}\big[B^{(2\ell+1)}(Y_p),B^{(2k-2\ell+1)}(Y_r)\big] .
\end{gather*}
Here $0\le\ell'\le\ell$, and we used (\ref{MixJac}), (\ref{wwid}) and (\ref{BY2}). By combining the expression above with \eqref{BX2} we complete the induction.

We are now in position to prove (\ref{B:XX})--(\ref{B:YY}) using induction on $k$. Notice that the base of induction for each of (\ref{B:XX})--(\ref{B:YY}) was already shown above. The induction step for~\eqref{B:XX} is analogous to the one in Section~\ref{sec:7.1}, thus we do not repeat the proof. That step uses~\eqref{BX3}, so it proves~\eqref{B:XX} when $X_\al$ lies in a simple ideal of $\fh$, and hence, by antisymmetry, when $X_\bet$ does. The remaining case arises when $\fh$ has a centre $\fk=\CC X_z$: we must show that ${\big[B^{(2a)}(X_z),B^{(2b)}(X_z)\big]=0}$, given~\eqref{B:XY} in degrees below $n=2a+2b$. Let $D_p^{\;q}=g_{zp}^{\sk q}$ and $M^{pq}=(\ka_\fm)^{pr}D_r^{\;s}D_s^{\;q}$, which is symmetric because $\ka_\fm$ is invariant under $\operatorname{ad} X_z$, and for odd $i$ set $\om_i=M^{pq}\big[B^{(i)}(Y_p),B^{(n-i)}(Y_q)\big]$, so that $\om_{n-i}=-\om_i$. For $b\ge1$, expanding $B^{(2b)}(X_z)$ by \eqref{72:BX}, where $\bar\fc_{(z)}=\fc_\fg$, and applying~\eqref{B:XY} gives
\begin{gather*}
\big[B^{(2a)}(X_z),B^{(2b)}(X_z)\big]=\fc_\fg^{-1}(\om_{2a+1}-\om_1).
\end{gather*}
The left-hand side is antisymmetric in $a$ and $b$, so $\om_{2a+1}=\om_{2a-1}+2\om_1$ for $1\le a<n/2$, whence $\om_{n-1}=(n-1)\om_1$. Since also $\om_{n-1}=-\om_1$, we get $\om_1=0$, so all $\om_i$ vanish, and so does the bracket. In the proof of \eqref{B:YY} below, the terms with $\ga=z$ vanish on both sides, since $f_z^{\;\al\bet}=0$ and $\fc_{(z)}=0$.

For \eqref{B:YY}, assuming induction, we have
\begin{gather*}
\big[B^{(2)}(X_\bet),\big[B^{(2k-2\ell+1)}(Y_p),B^{(2\ell+1)}(Y_q)\big]\big] = \mysum_{\al} w^{\sk \alpha}_{pq} \big[B^{(2)}(X_\bet),B^{(2k+2)}(X_\alpha)\big]. \end{gather*}
By \eqref{MixJac}, \eqref{ggid} and \eqref{BY2} we obtain
\begin{gather*}
\mysum_{\bet} g_{r}^{\;p\bet} \big[B^{(2)}(X_\bet),\big[B^{(2k-2\ell+1)}(Y_p),B^{(2\ell+1)}(Y_q)\big]\big]\\
\qquad{} = \mysum_{\bet} g_{r}^{\;p\bet} g_{\bet p}^{\sk s} \big[B^{(2k-2\ell+3)}(Y_s),B^{(2\ell+1)}(Y_q)\big] + \mysum_{\bet} g_{r}^{\;p\bet} g_{\bet q}^{\sk\, s} \big[B^{(2k-2\ell+1)}(Y_p),B^{(2\ell+3)}(Y_s)\big]\\
\qquad{} = \tfrac12 \fc_\fg \big[B^{(2k-2\ell+3)}(Y_r),B^{(2\ell+1)}(Y_q)\big] - \tfrac12 \mysum_{\bet} \bar\fc_{(\bet)} w_{rq}^{\sk\bet} B^{(2k+4)}(X_\bet) .
\end{gather*}
By \eqref{MixJac} and \eqref{BX2} we get
\begin{gather*}
\mysum_{\al,\bet} g_{r}^{\;p\bet} w^{\sk \alpha}_{pq} \big[B^{(2)}(X_\bet),B^{(2k+2)}(X_\alpha)\big] = \tfrac12 \mysum_{\al,\bet,\ga} w_{rq}^{\sk \ga} f_{\ga}^{\; \al\bet} \big[B^{(2)}(X_\bet),B^{(2k+2)}(X_\alpha)\big] \\
\hphantom{\mysum_{\al,\bet} g_{r}^{\;p\bet} w^{\sk \alpha}_{pq} \big[B^{(2)}(X_\bet),B^{(2k+2)}(X_\alpha)\big]}{} = \tfrac12 \mysum_\ga \fc_{(\ga)} w_{rq}^{\sk\ga} B^{(2k+4)}(X_\ga)
\end{gather*}
giving $\big[B^{(2k-2\ell+3)}(Y_r),B^{(2\ell+1)}(Y_q)\big] = \sum_{\bet} w_{rq}^{\sk\bet} B^{(2k+4)}(X_\bet)$, which combined with~\eqref{BB:YY} completes the induction.

Finally, for \eqref{B:XY}, assuming induction, we have
\begin{gather}
\mysum_\bet g_q^{\;\bet p} \big[B^{(1)}(Y_p),\big[B^{(2k-2\ell)}(X_\al),B^{(2\ell)}(X_\bet)\big]\big] = \mysum_\bet g_q^{\;\bet p} f_{\al\bet}^{\sk \ga} \big[B^{(1)}(Y_p),B^{(2k)}(X_\ga)\big]. \label{72:YXX}
\end{gather}
Using \eqref{ggid} and \eqref{BY2} we compute
\begin{gather*}
\mysum_\bet g_q^{\;\bet p}\big[B^{(1)}(Y_p),\big[B^{(2k-2\ell)}(X_\al),B^{(2\ell)}(X_\bet)\big]\big]\\
\qquad{} = \mysum_\bet g_q^{\;\bet p} g_{p\al}^{\sk r}\big[B^{(2k-2\ell+1)}(Y_r),B^{(2\ell)}(X_\bet)\big] + \mysum_\bet g_q^{\;\bet p} g_{p\bet}^{\sk r}\big[B^{(2k-2\ell)}(X_\al),B^{(2\ell+1)}(Y_r)\big]\\
\qquad{} = \mysum_\ga g_q^{\;\ga r} g_{r\al}^{\sk p}\big[B^{(2k-2\ell+1)}(Y_p),B^{(2\ell)}(X_\ga)\big] + \tfrac12 \fc_\fg \big[B^{(2k-2\ell)}(X_\al),B^{(2\ell+1)}(Y_q)\big].
\end{gather*}
Subtracting the first term in the last equality above from the rhs of \eqref{72:YXX} and using \eqref{MixJac}, \eqref{BY2} we get
\begin{gather*}
\mysum_\ga \Big( \mysum_\bet g_q^{\;\bet p} f_{\al\bet}^{\sk \ga} - g_q^{\;\ga r} g_{r\al}^{\sk p} \Big) \big[B^{(1)}(Y_p),B^{(2k)}(X_\ga)\big]\\
 \qquad {}= \mysum_\ga g_{\al q}^{\sk\, r} g_r^{\; \ga p} \big[B^{(1)}(Y_p),B^{(2k)}(X_\ga)\big] = \tfrac12 \fc_\fg g_{\al q}^{\sk\, r} B^{(2k+1)}(Y_r)
\end{gather*}
yielding $\big[B^{(2k-2\ell)}(X_\al),B^{(2\ell+1)}(Y_q)\big] = g_{\al q}^{\sk\, r} B^{(2k+1)}(Y_r)$, which completes the induction.

Thus we have proved that $\wt\cH^+ = \bigoplus \wt\cH^{(k)}$. By setting $\varphi \colon B^{(2k)}(X_\al) \mapsto X^{(2k)}_\al$, $B^{(2k+1)}(Y_p) \mapsto Y^{(2k+1)}_p$ we obtain a bijection of vector spaces, hence $\wt\cH^+ \cong \cH^+$.

\subsection[A proof of the uniqueness of the co-action of Y(g,g\^{}theta)\^{}\{tw\}]{A proof of the uniqueness of the co-action of $\boldsymbol{\cY(\fg,\fg^\theta)^{\rm tw}}$} \label{sec:7.3}

We need to show that the map $\LHb$ given by \eqref{dhBY} together with \eqref{fBY} is the unique map satisfying properties given by Definition \ref{D51} and property~(4) of Definition~\ref{D53}. The latter property defines the co-action up to the first order in $\hbar$,
\begin{gather*}
 \LHb(x) =\varphi(x) \ot 1 + 1 \ot x + \hbar \tau(x)+\cO\big(\hbar^2\big) ,
\end{gather*}
with $x \in \cY(\fg,\fh)^{\rm tw}$ and $\varphi$ the natural embedding $\cY(\fg,\fh)^{\rm tw}\hookrightarrow\cY(\fg)$.

For the degree-$0$ generators of $\cH^+$ the Lie coideal structure is trivial and the minimal form of the co-action is given by
\begin{gather}
 \LHb(X_\alpha) =\varphi(X_\alpha) \ot 1 + 1 \ot X_\alpha . \label{co1}
\end{gather}
The coideal compatibility identities \eqref{co-co-ass} and \eqref{co-inv} and the requirement $\varphi(X_\alpha)|_{\hbar\to 0}=X_\alpha$ implies that the natural inclusion for $X_\al$ is $\varphi \colon X_\alpha \in \cY(\fg,\fh)^{\rm tw} \mapsto X_\alpha \in \cY(\fg)$ and that~\eqref{co1} is indeed the unique co-action satisfying the required properties.

For the degree-$1$ generators of $\cH^+$ the Lie coideal structure is non-trivial. The minimal co-action is
\begin{gather}
 \LHb(\B(Y_p)) =\varphi(\B(Y_p)) \ot 1 + 1 \ot \B(Y_p)+\hbar [Y_p \ot 1, \Omega_\fh] . \label{co2}
\end{gather}
As previously, the co-action must satisfy \eqref{co-co-ass} and \eqref{co-inv} and the requirement $\varphi(\B(Y_p))|_{\hbar\to 0}=B(Y_p)$. By~\eqref{co-inv} we have
\begin{gather*}
\Delta_{\hbar}(\varphi(\B(Y_p)))=\varphi(\B(Y_p)) \ot 1+1\ot \varphi(\B(Y_p))+\hbar [Y_{p}\ot1,\Omega_\fh] .
\end{gather*}
Consider the ansatz $\varphi(\B(Y_p)) =\cJ(Y_p)+\hbar F^{(0)}_p$ with some degree-$0$ element $F^{(0)}_p \in \cY(\fg)$. Recall that $\Dh(\cJ(Y_p)) = \cJ(Y_p)\ot 1+1\ot \cJ(Y_p)+\frac{1}{2}\hbar [Y_p\ot1, \Omega_\fg]$. We rewrite the equality above as
\begin{gather*}
\Delta_{\hbar}\big(\cJ(Y_p)+\hbar F^{(0)}_p\big)-\big(\big(\cJ(Y_p)+\hbar F^{(0)}_p\big)\ot 1+1 \ot \big(\cJ(Y_p)+\hbar F^{(0)}_p\big)\big)=\hbar [Y_p \ot 1, \Omega_\fh] ,
\end{gather*}
which equates to
\begin{gather*}
\Delta_{\hbar}\big(F^{(0)}_p\big)-\big(F^{(0)}_p\ot 1+1 \ot F^{(0)}_p\big) = \tfrac{1}{2}[Y_p \ot 1, 2\Omega_\fh- \Omega_\fg]\\
\qquad{} = \tfrac{1}{2}\mysum_{\alpha}\big(g_{p}^{\;\alpha q} Y_q\ot X_\alpha-g_{p}^{\;q\alpha} X_\alpha \ot Y_q\big)
=\tfrac{1}{2}\mysum_{\alpha}g_{p}^{\;\alpha q} (Y_q\ot X_\alpha+ X_\alpha \ot Y_q) .
\end{gather*}
Thus we find that $F^{(0)}_p=\frac{1}{4}\sum_{\alpha}g_{p}^{\;\alpha q}(Y_q X_\alpha+ X_\alpha Y_q) + c Y_p= \tfrac{1}{4}[Y_{p},C_{X}] + c Y_p$ for any $c\in\CC$, which we can set to $c=0$ without loss of generality (this is equivalent to the composition $\varkappa_{-c} \circ \varphi$ with the shift-automorphism $\varkappa_c \colon \cJ(x) \mapsto \cJ(x) + \hbar c x$ of $\cY_\hbar(\fg)$ for any $x\in\fg$) giving
\begin{gather*}
 \varphi(\B(Y_p)) = \cJ(Y_p) + \tfrac{1}{4}\hbar [Y_p,C_X] ,
\end{gather*}
or in other words is there is a one-parameter family of embeddings that are equivalent to each other via the shift-automorphism of~$\cY(\fg)$.

It remains to verify that the coideal compatibility identities \eqref{co-co-ass} and \eqref{co-inv} hold, which imply that \eqref{co2} is the unique co-action, and to show that $\varphi$ is an algebra homomorphism, namely $\varphi([b_i,b_j])=[\varphi(b_i),\varphi(b_j)]$ and $\LHb([b_i,b_j])=[\LHb(b_i),\LHb(b_j)]$ for all $b_i,b_j\in\cY(\fg,\fh)^{\rm tw}$, which follow by a direct computation.

\subsection[A proof of the uniqueness of the co-action of Y(g,g)\^{}\{tw\}]{A proof of the uniqueness of the co-action of $\boldsymbol{\cY(\fg,\fg)^{\rm tw}}$} \label{sec:7.4}

Similarly as before, we need to show that the map $\LHb$ given by \eqref{dhx} and \eqref{dhGx} together with \eqref{fGx} is the unique map satisfying properties given by Definition~\ref{D51} and property~(4) of Definition~\ref{D53}. We first demonstrate an additional lemma that provides us with necessary identities.

\begin{Lemma} \label{L1}The following identities hold in $U(\fg)$ and $U(\fg)\ot U(\fg)$:
\begin{gather}
[[x_i\ot1,\Omega_\fg],\Omega_\fg] = \tfrac{1}{2}(\al_i^{\;jc} \al_j^{\;ab}+\al_i^{\;jb} \al_j^{\;ac})(x_a\ot \{x_b,x_c\} - \{x_b,x_c\} \ot x_a),\label{ID1} \\
\al_i^{\;jk} [[x_k\ot1,\Omega_\fg],[x_j\ot1,\Omega_\fg]] = \al_i^{\;jk} \al_{j}^{\;cr} \al_{k}^{\;bs} \al_{sr}^{\sk a} (x_a\ot \{x_b,x_c\} + \{x_b,x_c\} \ot x_a), \label{ID2} \\
6\al_i^{\;jk} \al_{j}^{\;cr} \al_{k}^{\;bs} \al_{sr}^{\sk a} = \al_i^{\;jk} \mysum_\pi \al_{j}^{\;\pi(c)r} \al_{k}^{\;\pi(b)s} \al_{sr}^{\sk \pi(a)} + \fc_\fg\big(\al_i^{\;jc} \al_j^{\;ab}+\al_i^{\;jb} \al_j^{\;ac}\big). \label{ID3}
\end{gather}
\end{Lemma}

\begin{proof}
Recall that
\begin{gather*}
[a\otimes b,c\otimes d]=[a,c]\otimes \{b,d\}+\{a,c\}\otimes [b,d].
\end{gather*}
For the first identity we have
\begin{gather*}
[[x_i\ot1,\Omega_\fg],\Omega_\fg] = \al_i^{\;bj}\eta^{kc}[x_j\ot x_b,x_k\ot x_c]\\
\hphantom{[[x_i\ot1,\Omega_\fg],\Omega_\fg]}{} = \al_i^{\;bj}\eta^{kc} (\al_{jk}^{\sk a} x_a \ot \{x_b,x_c\} + \al_{bc}^{\sk a} \{x_k,x_j\} \ot x_a)\\
\hphantom{[[x_i\ot1,\Omega_\fg],\Omega_\fg]}{} = \al_i^{\;bj}\eta^{kc} (\al_{jk}^{\sk a} x_a \ot \{x_b,x_c\} - \al_{jk}^{\sk a} \{x_c,x_b\} \ot x_a) \qquad \text{by ren. }k,j\leftrightarrow c,b\\
\hphantom{[[x_i\ot1,\Omega_\fg],\Omega_\fg]}{}= \al_i^{\;bj} \al_{j}^{\; ca} (x_a \ot \{x_b,x_c\} - \{x_c,x_b\} \ot x_a) .
\end{gather*}
For the second identity we have
\begin{gather*}
\al_i^{\;jk} [[x_k\ot1,\Omega_\fg],[x_j\ot1,\Omega_\fg]] =
\al_i^{\;jk} \al_k^{\;br}\al_j^{\;cs} [x_r \ot x_b, x_s \ot x_c]\\
\hphantom{\al_i^{\;jk} [[x_k\ot1,\Omega_\fg],[x_j\ot1,\Omega_\fg]]}{} = \al_i^{\;jk} \al_j^{\;cs} \al_k^{\;br}(\al_{rs}^{\sk a} x_a \ot \{x_b,x_c\} + \al_{bc}^{\sk a} \{x_s ,x_r\} \ot x_a )\\
\hphantom{\al_i^{\;jk} [[x_k\ot1,\Omega_\fg],[x_j\ot1,\Omega_\fg]]}{} = \al_i^{\;jk} \al_j^{\;cs} \al_k^{\;br} \al_{rs}^{\sk a} ( x_a \ot \{x_b,x_c\} + \{x_c,x_b\} \ot x_a )
\end{gather*}
by renaming $b,c\leftrightarrow r,s$. The third identity is obtained using the following auxiliary identities:
\begin{gather*}
\al_i^{\;jk} \al_{j}^{\;cr} \al_{k}^{\;bs} \al_{sr}^{\sk a} = \al_i^{\;jk} \al_{j}^{\;cr} \al_{k}^{\;as} \al_{sr}^{\sk b} + \tfrac{1}{2}\fc_\fg \al_i^{\;jc} \al_j^{\; ab} ,\\
\al_i^{\;jk} \al_{j}^{\;cr} \al_{k}^{\;bs} \al_{sr}^{\sk a} = \al_i^{\;jk} \al_{k}^{\;bs} \al_{j}^{\;ar} \al_{sr}^{\sk c} + \tfrac{1}{2}\fc_\fg \al_i^{\;jb} \al_j^{\; ac} .
\end{gather*}
The first auxiliary identity follows by multiple application of the Jacobi identity. In particular,
\begin{gather*}
\begin{split}
& \al_i^{\;jk} \al_{j}^{\;cr} \big(\al_{k}^{\;bs} \al_{sr}^{\sk a}-\al_{k}^{\;as} \al_{sr}^{\sk b}\big) = \al_i^{\;jk} \al_{j}^{\;cr} \al_{rk}^{\sk s} \al_s^{\; ab}= \tfrac12 \big( \al_i^{\;jk} \al_{j}^{\;cr} + \al_i^{\;jr} \al_{j}^{\;kc} \big) \al_{rk}^{\sk s} \al_s^{\; ab}\\
& \hphantom{\al_i^{\;jk} \al_{j}^{\;cr} (\al_{k}^{\;bs} \al_{sr}^{\sk a}-\al_{k}^{\;as} \al_{sr}^{\sk b}) }{}
 = \tfrac12 \al_i^{\;jc} \al_{j}^{\;kr} \al_{rk}^{\sk s} \al_s^{\; ab} = \tfrac{1}{2} \fc_\fg \al_i^{\;jc} \al_j^{\; ab} .
 \end{split}
\end{gather*}
The second auxiliary identity follows by the $b\leftrightarrow c$ symmetry and renaming $j,s\leftrightarrow k,r$ in the first term on the right hand side. We also have that
\begin{gather*}
\al_i^{\;jk} \mysum_\pi \al_{j}^{\;\pi(c)r} \al_{k}^{\;\pi(b)s} \al_{sr}^{\sk \pi(a)} = 2 \al_i^{\;jk} \big(
\al_{j}^{\;cr} \al_{k}^{\;bs} \al_{sr}^{\sk a} +
\al_{j}^{\;cr} \al_{k}^{\;as} \al_{sr}^{\sk b} +
\al_{j}^{\;ar} \al_{k}^{\;bs} \al_{sr}^{\sk c}\big) .
\end{gather*}
Hence taking a double sum of the auxiliary identities and adding $2\al_{j}^{\;cr} \al_{k}^{\;bs} \al_{sr}^{\sk a}$ to both sides of the resulting equality reproduces~\eqref{ID3} as required.
\end{proof}

The degree-0 generators have a trivial Lie coideal structure, hence, by the same arguments as before,
\begin{gather*}
\LHb(x_i) = \varphi(x_i) \ot 1 + 1 \ot x_i ,
\end{gather*}
and $\varphi \colon x_i\in\cY(\fg,\fg)^{\rm tw} \mapsto x_i\in\cY(\fg)$ is the natural inclusion. For the degree-2 generators we have
\begin{gather*}
\LHb(\G(x_i)) = \varphi(\G(x_i))\ot 1+1 \ot \G(x_i)+ \hbar [\cJ(x_i)\ot1,\Omega_{\fg}]+ \hbar^2 W^{(0)}_i ,
\end{gather*}
with $W^{(0)}_i$ being a degree-0 element in $\cY(\fg) \ot \cY(\fg,\fg)^{\rm tw}$, which can also be viewed as an element in $\cY(\fg) \ot \cY(\fg)$ via the natural inclusion. Coinvariance~\eqref{co-inv} implies
\begin{gather} \label{DfG}
\Delta_\hbar\big( \varphi(\G(x_i))\big)= \varphi(\G(x_i))\ot 1+1 \ot \varphi(\G(x_i))+ \hbar [\cJ(x_i)\ot1,\Omega_{\fg}]+ \hbar^2 W^{(0)}_i .
\end{gather}
We also must have $\varphi(\G(x_i))|_{\hbar\to 0}= \fc_\fg^{-1} \al_i^{\;jk} [J(x_k),J(x_j)]$. Let $K^{(2)}_i = \fc_\fg^{-1} \al_i^{\;jk} [\cJ(x_k),\cJ(x_j)]$ and choose the ansatz $\varphi(\G(x_i)) = K^{(2)}_i + \hbar F^{(1)}_i$ with some degree-1 element $F^{(1)}_i\in \cY(\fg)$. (We do not need to consider degree-0 elements in the ansatz since all $x_i\in \cY(\fg,\fg)^{\rm tw}$.) We then have
\begin{gather}
\Dh\big(K^{(2)}_i\big)= K^{(2)}_i \ot 1 + 1 \ot K^{(2)}_i \nonumber\\
\hphantom{\Dh\big(K^{(2)}_i\big)=}{} + \tfrac{1}{2}\hbar [\cJ(x_i) \ot 1 - 1 \ot \cJ(x_i), \Omega_\fg] + \tfrac{1}{4}\fc_\fg^{-1}\hbar^2 \al_i^{\;jk}[[x_k\ot1,\Omega_\fg],[x_j\ot1,\Omega_\fg]] , \label{DK2}
\end{gather}
which is obtained by making use of the following computation
\begin{gather*}
\al_i^{\;jk} [\cJ(x_k)\ot1+1\ot\cJ(x_k),[x_j\ot1,\Omega_\fg]] = \al_i^{\;jk}\al_j^{\;ab} \al_{kb}^{\sk c} ( \cJ(x_c)\ot x_a - x_a \ot\cJ(x_c))\\
\qquad{} = \tfrac12 \al_i^{\;jk} (\al_j^{\;ab} \al_{kb}^{\sk c} - \al_k^{\;ab} \al_{jb}^{\sk c} ) ( \cJ(x_c)\ot x_a - x_a \ot\cJ(x_c))\\
\qquad{} = \tfrac12 \al_i^{\;jk} \al_{kj}^{\sk b} \al_{b}^{\; ac} ( \cJ(x_c)\ot x_a - x_a \ot\cJ(x_c))\\
\qquad{} = \tfrac12 \fc_\fg \al_{i}^{\; ac} ( \cJ(x_c)\ot x_a - x_a \ot\cJ(x_c)) = \tfrac12 \fc_\fg [\cJ(x_i)\ot 1 - 1 \ot\cJ(x_i), \Omega_\fg ] .
\end{gather*}
By comparing the terms linear in $\hbar$ in \eqref{DfG} and~\eqref{DK2} we find $F^{(1)}_i = \tfrac{1}{4}[\cJ(x_i),C_\fg]$. Note that
\begin{gather*}
\Dh \big(F^{(1)}_i\big) = F^{(1)}_i\ot1 + 1\ot F^{(1)}_i + \tfrac{1}{2} [\cJ(x_i)\ot1+1\ot \cJ(x_i),\Omega_\fg] + \tfrac{1}{4}\hbar[[x_i\ot1,\Omega_\fg],\Omega_\fg] .
\end{gather*}
Lastly, by comparing the terms quadratic in $\hbar$ in \eqref{DfG} and using the expressions for $K^{(2)}_i$ and~$F^{(1)}_i$ we find
\begin{gather*} %\label{W0}
4 W^{(0)}_i = [[x_i\ot1,\Omega_\fg],\Omega_\fg] + \fc_\fg^{-1} \al_i^{\;jk} [[x_k\ot1,\Omega_\fg],[x_j\ot1,\Omega_\fg]] .
\end{gather*}
as required.

The coideal co-associativity identity \eqref{co-co-ass} in this case is not that straightforward, thus we demonstrate it explicitly. We need to show that the expression below equates to zero:
\begin{gather*}
\big((\Dh\ot \operatorname{id})\circ \LHb-(\operatorname{id}\ot \LHb)\circ \LHb\big) (\G(x_i)) \\
\qquad{} = \big(\Delta_\hbar (\varphi(\G(x_i)))-\big(\varphi(\G(x_i))\ot 1+1 \ot \varphi(\G(x_i))+\hbar [\cJ(x_i)\ot1,\Omega_{\fg}]\big)\big)\ot1 \\
\qquad\quad{} + \hbar^2\big((\Dh\ot \operatorname{id})\big(W^{(0)}_i\big)-(\operatorname{id}\ot\LHb)\big(W^{(0)}_i\big) - 1\ot W^{(0)}_i + \tfrac{1}{2} \al_i^{\;jk} ([x_k\ot 1, \Omega_\fg]) \ot x_j\big) .
\end{gather*}
Using the fact $\varphi(x_i)=x_i$ the relation above gives a constraint
\begin{gather}
(\Dh\ot \operatorname{id})\big(W^{(0)}_i\big) -(\operatorname{id} \ot \LHb)\big(W^{(0)}_i\big) + W^{(0)}_i\ot1 - 1\ot W^{(0)}_i\nonumber\\
 \qquad{}+ \tfrac{1}{2} \al_{i}^{\;jk} ([x_k\ot 1, \Omega_\fg]) \ot x_j = 0 . \label{W-constr}
\end{gather}
By Lemma~\ref{L1}, the explicit form of $W^{(0)}_i$ is equal to
\begin{gather}
W^{(0)}_i = \tfrac{1}{8}\big(\al_i^{\;jc} \al_j^{\;ab}+\al_i^{\;jb} \al_j^{\;ac}\big)(x_a\ot \{x_b,x_c\} - \{x_b,x_c\} \ot x_a) \nonumber\\
\hphantom{W^{(0)}_i =}{} + \tfrac{1}{24} \Big( \fc_\fg^{-1} \al_i^{\;jk} \mysum_\pi \al_{j}^{\;\pi(c)r} \al_{k}^{\;\pi(b)s} \al_{sr}^{\sk \pi(a)} + \al_i^{\;jc} \al_j^{\;ab}+\al_i^{\;jb} \al_j^{\;ac}\Big)\nonumber\\
\hphantom{W^{(0)}_i =}{}\times (x_a\ot \{x_b,x_c\} + \{x_b,x_c\} \ot x_a) \nonumber\\
\hphantom{W^{(0)}_i}{}
= \tfrac{1}{12}\big(\al_i^{\;jc} \al_j^{\;ab}+\al_i^{\;jb} \al_j^{\;ac}\big)( 2 x_a\ot \{x_b,x_c\} - \{x_b,x_c\} \ot x_a) \nonumber\\
\hphantom{W^{(0)}_i=}{} + \tfrac{1}{24}\fc_\fg^{-1} \al_i^{\;jk} \mysum_\pi \al_{j}^{\;\pi(c)r} \al_{k}^{\;\pi(b)s} \al_{sr}^{\sk \pi(a)} (x_a\ot \{x_b,x_c\} + \{x_b,x_c\} \ot x_a) .\label{W0E}
\end{gather}
Denote $W^{(0)}_i = h_i^{\;abc} x_a\ot \{x_b,x_c\} + \overline{h}_i^{\;abc} \{x_b,x_c\} \ot x_a$. Then \eqref{W-constr} is equivalent to
\begin{gather*}
4\overline{h}_i^{\;abc} x_c\ot x_b \ot x_a - 4 h_i^{\;abc} x_a\ot x_b\ot x_c + \al_{i}^{\;jc}\al_j^{\;ab} x_a \ot x_b \ot x_c = 0 ,
\end{gather*}
giving $4\overline{h}_i^{\;cba} - 4 h_i^{\;abc} + \al_{i}^{\;jc}\al_j^{\;ab}=0$. Using \eqref{W0E} we find
\begin{gather*}
4\overline{h}_i^{\;cba} - 4 h_i^{\;abc} + \al_{i}^{\;jc}\al_j^{\;ab}= -\tfrac{1}{3}\big(\al_i^{\;ja}\al_j^{\;cb}+\al_i^{\;jb}\al_j^{\;ca}\big)-\tfrac{2}{3}\big(\al_i^{\;jc}\al_j^{\;ab}+\al_i^{\;jb}\al_j^{\;ac}\big) + \al_i^{\;jc}\al_j^{\;ab} \\
\hphantom{4\overline{h}_i^{\;cba} - 4 h_i^{\;abc} + \al_{i}^{\;jc}\al_j^{\;ab}}{}
= -\tfrac{1}{3} \big(\al_i^{\;ja}\al_j^{\;cb} + \al_i^{\;jb}\al_j^{\;ac}\big) +\tfrac{1}{3} \al_i^{\;jc}\al_j^{\;ab} = 0
\end{gather*}
by the co-Jacobi identity. By reversing the computations above and using \eqref{co-inv} one can deduce that co-action $\LHb$ is the unique map satisfying the required properties up to the shift $(\varkappa_c \ot \operatorname{id})\circ\LHb$, i.e., there exists a one-parameter family of co-actions that are equivalent to each other via the shift-automorphism $\varkappa_c$ of $\cY(\fg)$.

It remains to verify that the map $\varphi$ is an algebra homomorphism, $\varphi([b_i,b_j])=[\varphi(b_i),\varphi(b_j)]$ and $\LHb([b_i,b_j])=[\LHb(b_i),\LHb(b_j)]$ for all $b_i,b_j\in\cY(\fg,\fg)^{\rm tw}$, which follow by a lengthy but direct computation.

\subsection{A proof of relations (\ref{H2}) and (\ref{H3})} \label{sec:7.5}

Recall the horrific degree-2 and degree-3 relations of the twisted Yangian $\cY\big(\fg,\fg^\theta\big)^{\rm tw}$ from the Theorem \ref{T:2}:
\begin{gather*}
[\B(Y_p),\B(Y_q)] + \mysum_\al (\bar\fc_{(\al)})^{-1} w_{pq}^{\sk\al}w_{\al}^{\;rs} [\B(Y_r),\B(Y_s)] = \hbar^2 \mysum_{\la,\mu,\nu}\Lambda_{pq}^{\la\mu\nu} \{X_\la,X_\mu,X_\nu\} , \\
[[\B(Y_p),\B(Y_q)],\B(Y_r)] + 2 \fc_\fg^{-1} \mysum_{\al} (\ka_\fm)^{tu} w_{pq}^{\sk\al} g_{r\al}^{\sk s} [[\B(Y_s),\B(Y_t)],\B(Y_u)] \\
\qquad{} = \hbar^2 \mysum_{\la,\mu,u} \Upsilon_{pqr}^{\la\mu u} \{X_\la,X_\mu,\B(Y_u)\} ,
\end{gather*}
where, by \eqref{H2L} and \eqref{H3U},
\begin{gather*}
\Lambda_{pq}^{\la\mu\nu} = \tfrac{1}{3} \Big( g^{\mu t}_{\sk p} g^{\la u}_{\sk q} + \mysum_{\al} (\bar\fc_{(\al)})^{-1} w_{pq}^{\sk \al} w_{\al}^{\; rs} g^{\mu t}_{\sk r} g^{\la u}_{\sk s} \Big) w_{tu}^{\sk \nu} , \\
\Upsilon_{pqr}^{\la\mu u} = \tfrac{1}{2} \mysum_{\al} \Big( w_{st}^{\sk\al} g_{p}^{\;\la s} g_{q}^{\;\mu t} g_{\al r}^{\sk u} + \mysum_\bet w_{pq}^{\sk\al} f_\al^{\;\la\bet} g_{r}^{\;\mu s} g_{\bet s}^{\sk u} + (\la\leftrightarrow\mu)\Big) \\
\hphantom{\Upsilon_{pqr}^{\la\mu u} =}{} {} + \fc_\fg^{-1} \mysum_{\al,\ga} (\ka_\fm)^{vx} w_{pq}^{\sk\ga} g_{r\ga}^{\sk y} \Big( w_{st}^{\sk\al} g_{y}^{\;\la s} g_{v}^{\;\mu t} g_{\al x}^{\sk u} + \mysum_\bet w_{yv}^{\sk\al} f_\al^{\;\la\bet} g_{x}^{\;\mu s} g_{\bet s}^{\sk u} + (\la\leftrightarrow\mu) \Big) .
\end{gather*}
The coideal structure is given by the co-action $\LHb \colon \cY\big(\fg,\fg^\theta\big)^{\rm tw} \to \cY(\fg) \ot \cY\big(\fg,\fg^\theta\big)^{\rm tw}$ defined in~\eqref{dhBY}:
\begin{gather*}
\LHb(X_\alpha) = X_\alpha \ot 1 + 1 \ot X_\alpha ,\qquad \LHb(\B(Y_p)) =\varphi(\B(Y_p)) \ot 1 + 1 \ot \B(Y_p) + \hbar [Y_p \ot 1, \Omega_X],
\end{gather*}
where, by \eqref{fBY}, $\varphi(\B(Y_p)) =\cJ(Y_p) + \tfrac{1}{4}\hbar [Y_p,C_X]$.

Denote the horrific degree-2 and degree-3 relations as
\begin{gather*}
H^{(2)}_{pq} = \hbar^2 \Lambda_{pq}^{\la\mu\nu}\{X_\la,X_\mu,X_\nu\} \qquad \text{and}\qquad
H^{(3)}_{pqr} = \hbar^2 \Upsilon_{pqr}^{\la\mu u}\{X_\la,X_\mu,\B(Y_u)\},
\end{gather*}
respectively. We need to show that co-action $\LHb$ extends to a homomorphism of algebras $\cY\big(\fg,\fg^\theta\big)^{\rm tw} \to \cY(\fg) \ot \cY\big(\fg,\fg^\theta\big)^{\rm tw}$. This can be checked by a direct computation. We will demonstrate this in several steps.

{\it Degree-2}. Set $B_{pq}=[\cB(Y_p),\cB(Y_q)]$. We have that
\begin{gather*}
\LHb(B_{pq}) = \varphi(B_{pq})\ot 1 + 1 \ot B_{pq} + \hbar [ \varphi(\B(Y_p)) \ot 1 + 1 \ot \B(Y_p) , [Y_q \ot 1, \Omega_X]] \\
\hphantom{\LHb(B_{pq}) =}{} + \hbar [ [Y_p \ot 1, \Omega_X] , \varphi(\B(Y_q)) \ot 1 + 1 \ot \B(Y_q) ] + \hbar^2 [ [Y_p \ot 1, \Omega_X] , [Y_q \ot 1, \Omega_X]] .
\end{gather*}
The $\hbar$-order terms in the right hand side of the equality above equate to
\begin{gather*}
 \hbar\Big[\cJ(Y_p) + \tfrac{1}{4}\hbar \mysum_\al g_p^{\;\al r} (Y_r X_\al + X_\al Y_r) \ot 1 + 1 \ot \B(Y_p) , \mysum_\al g_q^{\;\al r} Y_r \ot X_\al\Big] - (p\leftrightarrow q)\\
\qquad{} = \hbar \mysum_{\al,\bet} g_q^{\;\al r} w_{pr}^{\sk \bet}\cJ(X_\bet) \ot X_\al + \hbar \mysum_{\mu} g_q^{\;\mu t} g_{p\mu}^{\sk s} Y_t \ot \B(Y_s) \\
\qquad\quad{} + \tfrac12\hbar^2 \mysum_{\al,\bet} g_p^{\;\al r} g_q^{\;\bet s} \Big( \mysum_\ga w_{rs}^{\sk \ga} \{X_\ga, X_\al\} + g_{\al s}^{\sk\, t} \{Y_r, Y_t\}\Big) \ot X_\bet - (p\leftrightarrow q) ,
\end{gather*}
where $(p\leftrightarrow q)$ denote all the terms with indices $p$ and $q$ interchanged. In a similar way the $\hbar^2$-order term equates to
\begin{gather*}
\hbar^2 \mysum_{\al,\bet} g_p^{\;\al r}g_q^{\;\bet s} [Y_r \ot X_\al , Y_s \ot X_\bet ] \\
\qquad{} = \tfrac12\hbar^2 \mysum_{\al,\bet,\ga} \big(g_p^{\;\al s}g_q^{\;\bet t} - g_q^{\;\al s}g_p^{\;\bet t}\big) \big( w_{st}^{\sk \ga} X_\ga \ot \{ X_\bet , X_\al \} + f_{\al\bet}^{\sk\ga} \{Y_s , Y_t\} \ot X_\ga\big) .
\end{gather*}
Using mixed Jacobi identities we write
\begin{gather*}
 \mysum_{\al,\bet} \big( g_q^{\;\al r} w_{pr}^{\sk \bet} - g_p^{\;\al r} w_{qr}^{\sk \bet}\big) \cJ(X_\bet) \ot X_\al = \mysum_{\al,\bet,\mu} w_{pq}^{\sk \mu} f_\mu^{\;\al\bet} \cJ(X_\bet) \ot X_\al, \\
\mysum_{\mu} \big(g_q^{\;\mu t} g_{p\mu}^{\sk s} - g_p^{\;\mu t} g_{q\mu}^{\sk s}\big) Y_t \ot \B(Y_s) = \mysum_{\mu} w_{pq}^{\sk \mu} w_{\mu}^{\; st} Y_t \ot \B(Y_s)
\end{gather*}
and
\begin{gather*}
\mysum_{\al,\bet} \big( g_{\al s}^{\sk\, t} \big(g_p^{\;\al r} g_q^{\;\bet s} - g_q^{\;\al r} g_p^{\;\bet s}\big) \{Y_r, Y_t\} \ot X_\bet + \mysum_\ga f_{\al\bet}^{\sk\ga} \big(g_p^{\;\al r}g_q^{\;\bet t} - g_q^{\;\al r}g_p^{\;\bet t} \big) \{Y_r , Y_t\} \ot X_\ga \big) \\
\qquad {} = \mysum_{\al,\bet} \big( g_p^{\;\al r} \big( g_{\al s}^{\sk\, t} g_q^{\;\bet s} + \mysum_\ga f_{\al\ga}^{\sk\bet} g_q^{\;\ga t} \big) - g_q^{\;\al r} \big( g_{\al s}^{\sk\, t} g_p^{\;\bet s} + \mysum_\ga f_{\al\ga}^{\sk\bet} g_p^{\;\ga t} \big) \big) \{Y_r , Y_t\} \ot X_\bet \\
\qquad{} = \mysum_{\al,\bet} \big( g_p^{\;\al r} g_{s}^{\;t\bet} g_{q\al}^{\sk s} - g_q^{\;\al r} g_{s}^{\; t\bet} g_{p\al}^{\sk s} \big) \{Y_r , Y_t\} \ot X_\bet = \mysum_{\al,\mu} w_{pq}^{\sk \mu} w_\mu^{\;rs } g_{s}^{\; t\al} \{Y_r , Y_t\} \ot X_\al .
\end{gather*}
Using the identity
\begin{gather*}
\mysum_{\al,\bet,\ga} \big( w_{rs}^{\sk \ga} \big(g_p^{\;\al r} g_q^{\;\bet s} - g_q^{\;\al r} g_p^{\;\bet s}\big) \{X_\ga, X_\al\} \ot X_\bet + w_{st}^{\sk \ga} \big(g_p^{\;\al s}g_q^{\;\bet t} - g_q^{\;\al s}g_p^{\;\bet t}\big) X_\ga \ot \{ X_\bet , X_\al \} \big) \\
\qquad{} = \mysum_{\al,\bet,\ga} w_{st}^{\sk \ga} \big(g_p^{\;\al s} g_q^{\;\bet t} + g_p^{\;\bet s} g_q^{\;\al t}\big) (\{X_\ga, X_\al\} \ot X_\bet + X_\ga \ot \{ X_\bet , X_\al \} )
\end{gather*}
and combing all the expressions above we arrive to the following result,
\begin{gather}
\LHb(B_{pq}) = \varphi(B_{pq})\ot 1 + 1 \ot B_{pq} \nonumber\\
\qquad{}
 + \hbar \mysum_\mu w_{pq}^{\sk \mu} \Big( w_{\mu}^{\; st} Y_t \ot \cB(Y_s) + \mysum_{\al,\bet} f_{\mu}^{\;\al\bet} \cJ(X_\bet) \ot X_\al + \tfrac{1}{2}\hbar \mysum_{\al} w_{\mu}^{\;\,su}g_{u}^{\;t\al} \{Y_s,Y_t\} \ot X_\al \Big) \nonumber\\
\qquad{} + \tfrac{1}{2}\hbar^2 \mysum_{\al,\bet,\ga} w_{st}^{\sk\ga} \big(g_{p}^{\;\al s}g_{q}^{\;\bet t}+g_{p}^{\;\bet s}g_{q}^{\;\al t}\big) (\{X_\ga,X_\al\}\ot X_\bet + X_\ga \ot \{ X_\bet, X_\al \} ) \label{91:1}
\end{gather}
and, by similar computations,
\begin{gather*}
\varphi(B_{pq}) = [\J(Y_p),\J(Y_q)] + \tfrac{1}{2}\hbar \mysum_{\mu} w_{pq}^{\sk\mu} \Big( w_{\mu}^{\;st} \{\J(Y_s),Y_t\} + \mysum_{\al,\bet} f_{\mu}^{\;\al\bet} \{\J(X_\bet),X_\al\} \Big) \\
\hphantom{\varphi(B_{pq}) =}{} + \tfrac{1}{4}\hbar^2 g_{p}^{\;\al s}g_{q}^{\;\bet t} [\{Y_s,X_\al\},\{Y_t,X_\bet\}] .
\end{gather*}
In a similar way as it was done in the proof of Lemma \ref{L1}, the last line of \eqref{91:1} can be split as
\begin{gather}
 \tfrac{1}{6}\hbar^2 \mysum_{\al,\bet,\ga,\la,\mu} w_{pq}^{\sk \mu} f_\mu^{\;\la\ga}f_\la^{\;\bet\al}(\{X_\ga,X_\bet\}\ot X_\al - 2 X_\al \ot \{ X_\bet, X_\ga \} ) \nonumber\\
\qquad{} + \hbar^2 \mysum_{\al,\bet,\ga} g_{p}^{\;\al s}g_{q}^{\;\bet t}w_{st}^{\sk \ga} (\lan X_\al,X_\bet,X_\ga\ran_1 + \lan X_\al,X_\bet,X_\ga\ran_2 ) . \label{91:2}
\end{gather}
Here we have used the same notation as in \eqref{xxx}. The terms in the second line of \eqref{91:1} and in the first line of~\eqref{91:2} do not contribute to $\LHb\big(H^{(2)}_{pq}\big)$, since $w_{pq}^{\sk \mu} + \sum_\al (\bar\fc_{(\al)})^{-1} w_{pq}^{\sk \al}w^{\;rs}_{\al} w_{rs}^{\sk \mu} = 0$. What remains are the terms in the second line of~\eqref{91:2} giving
\begin{gather}
\LHb\big(H^{(2)}_{pq}\big) = \varphi\big(H^{(2)}_{pq}\big)\ot 1 + 1 \ot H^{(2)}_{pq} +\hbar^2 \mysum_{\al,\bet,\ga} \Big(g_p^{\; \al s} g_q^{\; \bet t} + \mysum_\la (\bar\fc_{(\la)})^{-1} w_{pq}^{\sk\; \la} w_\la^{\; uv} g_u^{\; \al s} g_{v}^{\;\bet t}\Big) \nonumber\\
\hphantom{\LHb\big(H^{(2)}_{pq}\big) =}{} \times w_{st}^{\sk \ga} (\lan X_\al,X_\bet,X_\ga\ran_1 + \lan X_\al,X_\bet,X_\ga\ran_2 ) . \label{91:3}
\end{gather}
By \eqref{Dxxx} we have
\begin{gather}
\LHb(\{X_\al,X_\bet,X_\ga\}) - \big(\{X_\al,X_\bet,X_\ga\} \ot 1 + 1 \ot \{X_\al,X_\bet,X_\ga\} \big) \nonumber\\
 \qquad{} = 3\big(\lan X_\al,X_\bet,X_\ga\ran_1 + \lan X_\al,X_\bet,X_\ga\ran_2 \big) . \label{91:4}
\end{gather}
Comparing \eqref{91:3} with \eqref{91:4} gives \eqref{H2L}, as required.

{\it Degree-3}. Set $B_{pqr}=[B_{pq},\cB(Y_r)]$. We write
\begin{gather}
\LHb(B_{pqr})=\varphi(B_{pqr})\ot1+1\ot B_{pqr} + \hbar B_{pqr}^{(2)} + \hbar^2 B_{pqr}^{(1)} + \hbar^3 B_{pqr}^{(0)}, \label{91:5}
\end{gather}
where $B_{pqr}^{(i)}$ denote elements of degree-$i$. By consistency, the elements of degree-2 in $\LHb\big(H^{(3)}_{pqr}\big)$ must vanish altogether, since they are of a linear order in $\hbar$. This can be shown explicitly. We have
\begin{gather*}
B_{pqr}^{(2)} = \mysum_\al [[\J(Y_p),\J(Y_q)]\ot 1 + 1 \ot [\B(Y_p),\B(Y_q)] , g_r^{\;\al t} Y_t\ot X_\al ] \\
\hphantom{B_{pqr}^{(2)} =}{} + \mysum_\mu w_{pq}^{\sk \mu} \Big[ w_{\mu}^{\; st} Y_t \ot \cB(Y_s) + \mysum_{\al,\bet} f_{\mu}^{\;\al\bet} \cJ(X_\bet) \ot X_\al , \J(Y_r)\ot1+1\ot\B(Y_r)\Big ] .
\end{gather*}
Denote this expression as $B_{pqr}^{(2)}=B_{pqr}^{(2,0)}+B_{pqr}^{(0,2)}+B_{pqr}^{(1,1)}$, where $B_{pqr}^{(i,j)}$ represents elements of degree-$i$ in the left tensor factor and of degree-$j$ in the right tensor factor. Then
\begin{gather}
B_{pqr}^{(2,0)} = \mysum_{\al,\bet} g_r^{\;\al t} \big(w_{pt}^{\sk \bet}[\J(X_\bet),\J(Y_q)]+w_{qt}^{\sk\bet}[\J(Y_p),\J(X_\bet)]\big)\ot X_\al \nonumber\\
\hphantom{B_{pqr}^{(2,0)} =}{} + \mysum_{\al,\bet,\mu} w_{pq}^{\sk \mu} f_{\mu}^{\;\al\bet} [\cJ(X_\bet),\cJ(Y_r)] \ot X_\al \nonumber\\
\hphantom{B_{pqr}^{(2,0)} }{}= 2\mysum_{\al,\bet,\ga} \Big(g_r^{\;\al t} (w_{pt}^{\sk \bet}g_{\bet q}^{\sk s}-w_{qt}^{\sk\bet}g_{\bet p}^{\sk s}) + \mysum_{\mu} w_{pq}^{\sk \mu} f_{\mu}^{\;\al\bet} g_{\bet r}^{\sk s} \Big) \fc_\fg^{-1} g_s^{\;\ga u} [\cJ(Y_u),\cJ(X_\ga)] \ot X_\al \nonumber\\
\hphantom{B_{pqr}^{(2,0)} =}{} + \tfrac{1}{2}\hbar^2\fc_\fg^{-1} \mysum_{\al,\bet,\ga} \Big( g_r^{\;\al t} (w_{pt}^{\sk \bet} g_{q}^{\ga u} - w_{qt}^{\sk \bet} g_{p}^{\ga u}) + \mysum_\mu w_{pq}^{\sk \mu} f_{\mu}^{\;\al\bet} g_{r}^{\ga u} \Big)\nonumber\\
\hphantom{B_{pqr}^{(2,0)} =}{} \times\mysum_{ijk} \mathcal{A}_{\bet u \ga}^{ijk} \{x_i,x_j,x_k\}\ot X_\al \nonumber\\
\hphantom{B_{pqr}^{(2,0)} }{}
= 2\mysum_{\al,\bet,\ga} w_{pq}^{\sk \bet} g_{\bet r}^{\sk\, t} g_{t}^{\;\al s} g_s^{\;\ga u} \fc_\fg^{-1} [\cJ(Y_u),\cJ(X_\ga)] \ot X_\al + \mathcal{O}_1\big(\hbar^2\big) , \label{91:6}
\end{gather}
where $\mathcal{O}_1(\hbar^2)$ is a short--hand notation for the terms in the third line, and we have used the mixed Jacobi identities and the Drinfeld terrific relation \eqref{DT2} in the form
\begin{gather*}
[\J(X_\bet),\J(Y_q)] + {2}{\fc_\fg^{-1}} g_{\bet q}^{\sk\, s} g_{s}^{\;\ga u} [\J(X_\ga),\J(Y_u)] = \tfrac12 \hbar^2\fc_\fg^{-1} \mysum_{\ga} \mysum_{ijk} g_{q}^{\;\ga u} \mathcal{A}_{\bet u \ga}^{ijk} \{x_i,x_j,x_k\} .
\end{gather*}
Here the sum $\sum_i x_i$ spans all of the symmetric space basis, $x_i=\{Y_p,X_\al\}$; the same applies to~$x_j$ and~$x_k$. In a similar way we find
\begin{gather}
B_{pqr}^{(0,2)} = \mysum_\al g_r^{\;\al t} Y_t \ot \big(g_{p\al}^{\sk u} [\B(Y_u),\B(Y_q)] + g_{q\al}^{\sk u} [\B(Y_p),\B(Y_u)]\big) \nonumber\\
\hphantom{B_{pqr}^{(0,2)} =}{} + \mysum_\al w_{pq}^{\sk \al} w_{\al}^{\; st} Y_t \ot [\cB(Y_s),\B(Y_r)] \nonumber\\
\hphantom{B_{pqr}^{(0,2)}}{}= \mysum_{\al,\bet} \Big(g_r^{\;\al t} (g_{p\al}^{\sk u} w_{uq}^{\sk\bet} - g_{q\al}^{\sk u} w_{up}^{\sk\bet}) + w_{pq}^{\sk \al} w_{\al}^{\; ut} w_{ur}^{\sk\bet}\Big) (\bar\fc_{(\bet)})^{-1} w_\bet^{\;vs} Y_t \ot [\B(Y_s),\B(Y_v)] \nonumber\\
\hphantom{B_{pqr}^{(0,2)} =}{} + \hbar^2 \mysum_{\al,\la,\mu,\nu} \Big(g_r^{\;\al t} (g_{p\al}^{\sk u} \La_{uq}^{\la\mu\nu} - g_{q\al}^{\sk u} \La_{up}^{\la\mu\nu}) + w_{pq}^{\sk \al} w_{\al}^{\; st} \La_{sr}^{\la\mu\nu} \Big) Y_t \ot \{X_\la,X_\mu,X_\nu\} \nonumber\\
\hphantom{B_{pqr}^{(0,2)}}{}= \mysum_{\al,\bet} w_{pq}^{\sk\al} g_{\al r}^{\sk\; u} g_u^{\;\bet t} w_\bet^{\;vs} (\bar\fc_{(\bet)})^{-1} Y_t \ot [\B(Y_s),\B(Y_v)] + \mathcal{O}_2\big(\hbar^2\big) , \label{91:7}
\end{gather}
where we have used the mixed Jacobi identities and \eqref{H2}. The notation for $\mathcal{O}_2(\hbar^2)$ is the same as before. The remaining elements give
\begin{gather}
B_{pqr}^{(1,1)} = \mysum_{\mu,\beta} w_{pq}^{\sk \mu} \big( w_{\mu}^{\; st} w_{tr}^{\sk\beta} + \mysum_{\al} f_{\mu}^{\;\al\bet} g_{\al r}^{\sk s} \big) \cJ(X_\bet) \ot \B(Y_s) \nonumber\\
\hphantom{B_{pqr}^{(1,1)}}{} = \mysum_{\al,\beta} w_{pq}^{\sk \al} g_{\al r}^{\sk\, t} g_t^{\;s\bet} \cJ(X_\bet) \ot \B(Y_s) . \label{91:8}
\end{gather}
Now it is easy to check that leading terms in (\ref{91:6})--(\ref{91:8}) do not contribute to $\LHb\big(H_{pqr}^{(3)}\big)$, since (by~\eqref{ggid})
\begin{gather}
\mysum_\bet w_{pq}^{\sk \bet} g_{\bet r}^{\sk\, t} + 2\fc_{\fg}^{-1} \mysum_{\al,\bet} (\ka_\fm)^{uv} w_{pq}^{\sk\al} g_{r\al}^{\sk\, s} w_{su}^{\sk \bet} g_{\bet v}^{\sk\, t} = 0 . \label{91:9}
\end{gather}
Next step is to consider the degree-1 term in \eqref{91:5}. We write $B_{pqr}^{(1)}=B_{pqr}^{(1,0)}+B_{pqr}^{(0,1)}$. Using \eqref{91:1}, \eqref{91:2} and \eqref{fBY} we find
\begin{gather}
B_{pqr}^{(1,0)} = \tfrac{1}{2} \mysum_{\al,\mu} w_{pq}^{\sk \mu} \Big(\! \mysum_{\bet,\ga} f_{\mu}^{\;\al\bet} g_{r}^{\;\ga s} [\cJ(X_\bet) \ot X_\al,2Y_s\ot X_\ga + \{Y_s,X_\ga\}\ot1]\nonumber\\
\hphantom{B_{pqr}^{(1,0)} =}{} + w_{\mu}^{\;su}g_{u}^{\;t\al} [\{Y_s,Y_t\},\J(Y_r)] \ot X_\al \Big) \nonumber\\
\hphantom{B_{pqr}^{(1,0)} =}{} + \tfrac{1}{2}\mysum_{\ga,\mu} w_{pq}^{\sk\mu} g_{r}^{\;\ga u}\Big( w_{\mu}^{\;st} [\{\J(Y_s),Y_t\}, Y_u ] + \mysum_{\al,\bet} f_{\mu}^{\;\al\bet} [\{\J(X_\bet),X_\al\}, Y_u] \Big) \ot X_\ga \nonumber\\
\hphantom{B_{pqr}^{(1,0)} =}{} + \tfrac{1}{2} \mysum_{\al,\bet,\ga} w_{st}^{\sk\ga} \big(g_{p}^{\;\al s}g_{q}^{\;\bet t}+g_{p}^{\;\bet s}g_{q}^{\;\al t}\big) [\{X_\ga,X_\al\}\ot X_\bet + X_\ga \ot \{ X_\bet, X_\al \},\J(Y_r)\ot1] \nonumber\\
\hphantom{B_{pqr}^{(1,0)}}{} = \tfrac{1}{2} \mysum_{\ga,\la,\mu} w_{pq}^{\sk \mu} \Big( \mysum_{\al,\bet} g_r^{\;\bet s} \big( 2 f_\mu^{\;\al\la} f_{\al\bet}^{\sk\ga} + f_\mu^{\;\ga\al} f_{\al\bet}^{\sk \la}\big) + w_{tr}^{\sk\la}\big(w_\mu^{\;tu}g_u^{\;s\ga}+ w_\mu^{\;su}g_u^{\;t\ga}\big) \nonumber\\
\hphantom{B_{pqr}^{(1,0)} =}{} + g_r^{\;\ga u}\big(w_\mu^{\;ts}w_{tu}^{\sk\la}+\mysum_{\al}f_\mu^{\;\al\la}g_{\al u}^{\sk s}\big)\Big) \{Y_s,\J(X_\la)\} \ot X_\ga \nonumber\\
\hphantom{B_{pqr}^{(1,0)} =}{} + \tfrac{1}{2}\mysum_{\al,\bet,\ga} \big( \mysum_{\mu} w_{pq}^{\sk\mu} f_\mu^{\;\al\bet} g_{r}^{\;\ga s} g_{\bet s}^{\sk u} + w_{st}^{\sk\bet} g_{p}^{\;\al s} g_{q}^{\;\ga t} g_{\bet r}^{\sk u} \big) \nonumber\\
\hphantom{B_{pqr}^{(1,0)} =}{} \times ( 2 \cJ(Y_u) \ot \{X_\al,X_\ga\} + \{X_\al,\J(Y_u)\}\ot X_\ga + \{X_\ga,\J(Y_u)\}\ot X_\al )\nonumber\\
\hphantom{B_{pqr}^{(1,0)} =}{} + \tfrac{1}{2}\mysum_{\al,\bet,\ga} w_{ts}^{\sk\ga}\big( w_{pq}^{\sk\bet} w_{\bet}^{\;ut} g_{r}^{\;\al s} + g_{r\bet}^{\sk u} \big(g_{p}^{\;\bet s}g_{q}^{\;\al t}+g_{p}^{\;\al s} g_{q}^{\;\bet t}\big) \big) \{X_\ga,\J(Y_u)\}\ot X_\al . \label{91:10}
\end{gather}
The first sum in the last equality can be simplified using the mixed Jacobi identities, giving
\begin{gather*}
\tfrac{1}{2} \mysum_{\ga,\la,\mu} w_{pq}^{\sk \mu} g_{\mu r}^{\sk t} \big( g_t^{\;u\la} g_u^{\ga s} + \mysum_\al g_t^{\;\al s} f_\al^{\ga\la} \big) \{Y_s,\J(X_\la)\} \ot X_\ga .
\end{gather*}
This component does not contribute to $\LHb\big(H_{pqr}^{(3)}\big)$ due to \eqref{91:9}. The second degree-1 element, $B_{pqr}^{(0,1)}$, is
\begin{gather}
B_{pqr}^{(0,1)} = \tfrac{1}{2} \mysum_{\al,\bet} w_{pq}^{\sk \bet} \Big( w_{\bet}^{\; st} g_{r}^{\;\al u} [Y_t \ot \cB(Y_s),2 Y_u\ot X_\al+\{Y_u,X_\al\}\ot1] \nonumber\\
\hphantom{B_{pqr}^{(0,1)} =}{} + w_{\bet}^{\;su}g_{u}^{\;t\al} \{Y_s,Y_t\} \ot [X_\al,\B(Y_r)] \Big) \nonumber\\
\hphantom{B_{pqr}^{(0,1)} =}{} + \tfrac{1}{2} \mysum_{\al,\bet,\ga} w_{st}^{\sk\ga} \big(g_{p}^{\;\al s}g_{q}^{\;\bet t}+g_{p}^{\;\bet s}g_{q}^{\;\al t}\big) [\{X_\ga,X_\al\}\ot X_\bet + X_\ga \ot \{ X_\bet, X_\al \},1\ot\B(Y_r)] \nonumber\\
\hphantom{B_{pqr}^{(0,1)} }{} = \tfrac{1}{2} \mysum_{\al,\bet} w_{pq}^{\sk \bet} \big( 2 w_{\bet}^{\; st} g_{r}^{\;\al u} g_{s\al}^{\sk v} + w_{\bet}^{\; vs} g_{r}^{\;\al u} g_{s\al}^{\sk t} + w_{\bet}^{\;ts}g_{s}^{\;u\al} g_{\al r}^{\sk v} \big) \{Y_t,Y_u\} \ot \B(Y_v) \nonumber\\
\hphantom{B_{pqr}^{(0,1)} =}{} + \mysum_{\al,\bet,\ga} w_{ts}^{\sk\ga}\big( w_{pq}^{\sk \bet} w_{\bet}^{\; ut} g_{r}^{\;\al s} + g_{r\bet}^{\sk u} \big(g_{p}^{\;\al s}g_{q}^{\;\bet t}+g_{p}^{\;\bet s} g_{q}^{\;\al t}\big) \big)\nonumber\\
\hphantom{B_{pqr}^{(0,1)} =}{} \times \Big( X_\ga\ot\{X_\al,\B(Y_u)\} + \tfrac{1}{2} \{X_\al,X_\ga\}\ot\B(Y_u) \Big) . \label{91:11}
\end{gather}
The first sum in the last equality above can be simplified to $\frac{1}{2} \sum_{\al,\bet} w_{pq}^{\sk \bet} g_{r\bet}^{\sk\, s} g_{s}^{\;\al t} w_{\al}^{\;uv} \{Y_t,Y_u\} \ot \B(Y_v)$, and does not contribute to $\LHb\big(H_{pqr}^{(3)}\big)$ due to \eqref{91:9}. Set
\begin{gather*}
A_{pqr}^{\al\ga u} = \tfrac{1}{2}\mysum_{\bet} \big(\mysum_{\mu} w_{pq}^{\sk\mu} f_\mu^{\;\al\bet} g_{r}^{\;\ga s} g_{\bet s}^{\sk u} + w_{st}^{\sk\bet} g_{p}^{\;\al s} g_{q}^{\;\ga t} g_{\bet r}^{\sk u} + (\al\leftrightarrow\ga)\big) , \\
B_{pqr}^{\al\ga u} = \tfrac{1}{2}\mysum_{\bet} w_{ts}^{\sk\ga}\big( w_{pq}^{\sk \bet} w_{\bet}^{\; ut} g_{r}^{\;\al s} + g_{r\bet}^{\sk u} \big(g_{p}^{\;\al s}g_{q}^{\;\bet t}+g_{p}^{\;\bet s} g_{q}^{\;\al t}\big) \big) .
\end{gather*}
Using the mixed Jacobi identities \eqref{MixJac} we find that
\begin{gather*}
\Upsilon_{pqr}^{\al\ga u} = A_{pqr}^{\al\ga u} + 2\fc_\fg^{-1} \mysum_\mu (\ka_\fm)^{tv}
 w_{pq}^{\sk\mu} g_{r\mu}^{\sk s} A_{stv}^{\al\ga u} = B_{pqr}^{\al\ga u} + 2{\fc_\fg^{-1}} \mysum_\mu (\ka_\fm)^{tv} w_{pq}^{\sk\mu} g_{r\mu}^{\sk s} B_{stv}^{\al\ga u} %\label{G}
\end{gather*}
satisfying $\Upsilon_{pqr}^{\al\ga u}=\Upsilon_{pqr}^{\ga\al u}$.
The contribution of the remaining sums in~\eqref{91:10} and~\eqref{91:11} to $\LHb\big(H_{pqr}^{(3)}\big)$~is
\begin{gather}
\begin{split}
& \mysum_{\al,\ga} \Upsilon_{pqr}^{\al\ga u} \big( \cJ(Y_u) \ot \{X_\al,X_\ga\} + 2\{X_\al,\J(Y_u)\}\ot X_\ga \\
& \qquad {} + 2X_\ga\ot\{X_\al,\B(Y_u)\} + \{X_\al,X_\ga\}\ot\B(Y_u) \big) .\end{split}
 \label{91:12}
\end{gather}
Using \eqref{Dxxx} and \eqref{xxx} we find that
\begin{gather*}
\LHb(\{X_\al,X_\ga,\B(Y_u)\}) - (\varphi(\{X_\al,X_\ga,\B(Y_u)\}) \ot 1 + 1 \ot \{X_\al,X_\ga,\B(Y_u)\} + \mathcal{O}_3(\hbar)) \\
\qquad{} = \J(Y_u)\ot\{X_\al,X_\ga\} + \{X_\al,X_\ga\}\ot\B(Y_u) \\
\qquad\quad{} + (\{\J(Y_u),X_\al\}\ot X_\ga + X_\ga\ot \{X_\al,\B(Y_u)\} + (\al \leftrightarrow \ga)) ,
\end{gather*}
where $\mathcal{O}_3(\hbar)$ denotes linear in $\hbar$ terms that appear due to \eqref{fBY}. Comparing the relation above with~\eqref{91:12} and using symmetry $\Upsilon_{pqr}^{\al\ga u}=\Upsilon_{pqr}^{\ga\al u}$ gives~\eqref{H3U}, as required. The remaining elements in $\LHb\big(H_{pqr}^{(3)}\big)$, that we have denoted by $\hbar^3 B_{pqr}^{(0)}$, are cubic in $\hbar$ and, by summing with $\hbar \mathcal{O}_1(\hbar^2)$ of~\eqref{91:6} and $\hbar \mathcal{O}_2(\hbar^2)$ of~\eqref{91:7}, and repeating similar steps as we did above, give precisely the element $\hbar^2 \mathcal{O}_3(\hbar)$, thus the relations~\eqref{H2} and~\eqref{H3} are compatible with the co-action~\eqref{dhBY}.

\subsection{An outline of the proof of relation~(\ref{H4})} \label{sec:7.6}

Applying $\LHb$ to both sides of \eqref{H4}, we determine the coefficients
$\Psi_{abc}^{ijk}$, $\overline{\Phi}_{abc}^{\,ijklm}$ and $\Phi_{abc}^{ijk}$ by
comparing three components of the two images. A complete proof would also check
that all the remaining components agree, as was done for~\eqref{H2} and
\eqref{H3} in Section~\ref{sec:7.5}; we do not do this here.
We shall use the following elementary identities.

\begin{Lemma} \label{L2} Let $A$ be an associative algebra over
$\CC$.  For arbitrary elements $x_a,x_b,\ldots\in A$, the following
identities hold:
\begin{gather*}
\{x_i,\{x_j,x_k\}\} = \{x_i,x_j,x_k\}+\tfrac{1}{12} [x_{(j},[x_{k)},x_i]] ,\\ % \label{L2:1} \\
\{x_{j},x_{k},x_l,x_{m}\} = \tfrac{1}{3}\{x_{j},x_{(k},\{x_l,x_{m)}\}\} - \tfrac{1}{36} [[x_{j},x_{(k}],\{x_l,x_{m)}\}] ,\\ % \label{L2:2} \\
\{1\ot x_i,1\ot x_j,x_a\ot \{x_b,x_c\}\} = x_a \ot \{x_i,x_j,\{x_b,x_c\}\} , \\
\{x_i \ot 1,1\ot x_j,x_a\ot \{x_b,x_c\}\} = \{x_i, x_a\}\ot \{x_j, \{x_b,x_c\}\}-\tfrac{1}{12}[x_a,x_i]\otimes[x_j, \{x_b,x_c\}], \\
\{x_i \ot 1,1\ot x_j,\{x_b,x_c\} \ot x_a\} = \{x_i, \{x_b,x_c\}\}\ot \{x_j,x_a\}-\tfrac{1}{12} [\{x_b,x_c\},x_i] \ot[x_j, x_a]. %\label{L2:5}
\end{gather*}
\end{Lemma}

We apply the lemma after decomposing each tensor factor according to
PBW degree. Set $q(a,b)=\{a,b\}$, and let
$\operatorname{pr}_r\colon U(\fg)\to S^r\fg$ denote the projection induced
by the PBW symmetrisation $S(\fg)\xrightarrow{\sim}U(\fg)$.  The projections
needed below are the following universal maps, for $a,b,c,d\in\fg$:
\begin{gather*}
M(a;b,c) :=\operatorname{pr}_1\{a,q(b,c)\}
=\tfrac1{12}\big([b,[c,a]]+[c,[b,a]]\big),\\
N(a,b;c,d) :=\operatorname{pr}_2\{a,b,q(c,d)\}\\
\hphantom{N(a,b;c,d)\,}{} =\tfrac1{12}\{[a,c],[b,d]\}+\tfrac1{12}\{[a,d],[b,c]\}
 +\tfrac1{12}\{a,[b,[c,d]]\}-\tfrac16\{a,[c,[b,d]]\}\\
 \hphantom{N(a,b;c,d) :=}{}+\tfrac1{12}\{b,[a,[c,d]]\}-\tfrac16\{b,[c,[a,d]]\},\\
K(a,b;c,d) :=\operatorname{pr}_2\{q(a,b),q(c,d)\}\\
\hphantom{K(a,b;c,d)\,}{}
=\tfrac14\{[a,c],[b,d]\}+\tfrac14\{[a,d],[b,c]\}
 +\tfrac1{12}\{a,[b,[c,d]]\}-\tfrac16\{a,[c,[b,d]]\}\\
 \hphantom{K(a,b;c,d) :=}{}+\tfrac1{12}\{b,[a,[c,d]]\}-\tfrac16\{b,[c,[a,d]]\}
 +\tfrac1{12}\{c,[a,[b,d]]\}+\tfrac1{12}\{c,[b,[a,d]]\}\\
\hphantom{K(a,b;c,d) :=}{}+\tfrac1{12}\{d,[a,[b,c]]\}+\tfrac1{12}\{d,[b,[a,c]]\},\\
L(a,b;c,d) :=\operatorname{pr}_1[q(a,b),q(c,d)]\\
 \hphantom{L(a,b;c,d)\,}{}
 =\tfrac1{12}[a,[b,[c,d]]]-\tfrac16[a,[c,[b,d]]]
 +\tfrac1{12}[b,[a,[c,d]]]-\tfrac16[b,[c,[a,d]]].
\end{gather*}
The underlying decompositions already hold in the free associative algebra
generated by $a$, $b$, $c$, $d$.  For example,
the PBW decomposition of a nested Jordan product is
\begin{gather*}
\{a,b,q(c,d)\}=\{a,b,c,d\}+N(a,b;c,d).
\end{gather*}
Thus its degree-two component must be retained when the co-action is
projected to PBW bidegree $(1,2)$.

Let
\begin{gather*}
Q_c=h_c^{buv}x_b\ot q(x_u,x_v)
 +\overline h_c^{buv}q(x_u,x_v)\ot x_b
\end{gather*}
be the coefficient of $\hbar^2$ in \eqref{dhGx2}, so that
$\LHb(\G(x_c))=\varphi(\G(x_c))\ot1+1\ot\G(x_c)+\hbar[\cJ(x_c)\ot1,\Omega_\fg]+\hbar^2Q_c$.
The left-hand side of \eqref{H4} is the sum of
$\al_{ab}^{\sk d}[\G(x_c),\G(x_d)]$ over the cyclic permutations of $a$, $b$, $c$.

The coefficient $\Psi$ comes from the terms of the form $\{x_i,x_j\}\ot\G(x_k)$ at
order $\hbar^2$. On the left-hand side they arise only from
$[Q_c,1\ot\G(x_d)]+[1\ot\G(x_c),Q_d]$, which contributes
$\smash{\big(\al_{cr}^{\sk k}\overline{h}_d^{rij}}-\smash{\al_{dr}^{\sk k}\overline{h}_c^{rij}\big)
\{x_i,x_j\}\ot\G(x_k)}$ to $\LHb([\G(x_c),\G(x_d)])$. On the right-hand side they arise only from
$\Psi_{abc}^{ijk}\{x_i\ot1,x_j\ot1,1\ot\G(x_k)\}=\Psi_{abc}^{ijk}\{x_i,x_j\}\ot\G(x_k)$.
Comparing the two gives \eqref{coefevenloop1}.

The coefficient $\overline\Phi$ comes from the terms
$x_i\ot\{x_j,x_k,x_l,x_m\}$, of PBW bidegree~$(1,4)$, at order~$\hbar^4$.
On the left-hand side they arise only from the term
\begin{gather*}
h_c^{rjk}h_d^{slm}[x_r,x_s]\ot\{q(x_j,x_k),q(x_l,x_m)\}
\end{gather*}
of $[Q_c,Q_d]$, which contributes
$\al_{rs}^{\sk i}\al_{(ab}^{\qu d}h_{c)}^{rjk}h_{d}^{slm}\,x_i\ot\{x_j,x_k,x_l,x_m\}$.
On the right-hand side, the term $\Psi_{abc}^{jkr}\{x_j,x_k,\G(x_r)\}$ contributes
$\Psi_{abc}^{jkr}h_r^{ilm}\,x_i\ot\{x_j,x_k,x_l,x_m\}$, by the third identity of
Lemma~\ref{L2}, and by \eqref{Dxxx} the quintic term contributes
$5\overline{\Phi}_{abc}^{\,ijklm}\,x_i\ot\{x_j,x_k,x_l,x_m\}$. Comparing them up to
symmetrisation of the upper indices gives \eqref{coefevenloop2}.

The coefficient $\Phi$ comes from the terms of PBW bidegree $(1,2)$ at order
$\hbar^4$, and this is where the projections above are needed. Write
$\operatorname{pr}_{1,2}=\operatorname{pr}_1\ot\operatorname{pr}_2$.
The other terms of $\LHb(\G(x_c))$ do not contribute here, on either side of
\eqref{H4}: every product
involving them has a tensor factor containing $\G(x)$, $\varphi(\G(x))$ or a single
$\cJ(x)$, or a right tensor factor of PBW degree at most one, except for the terms
$[\cJ(x_e),\cJ(x_f)]\ot\{x_u,x_v\}$ of
$\hbar^2\big[[\cJ(x_c)\ot1,\Omega_\fg],[\cJ(x_d)\ot1,\Omega_\fg]\big]$. Since
$H_2(\fg)=0$, the kernel of the bracket $\Lambda^2\fg\to\fg$ is spanned by the
elements $a\wedge[b,c]+b\wedge[c,a]+c\wedge[a,b]$, so
$[\cJ(x_e),\cJ(x_f)]-\al_{ef}^{\sk g}\cK(x_g)$, with $\cK$ as in
Example~\ref{Y30}, is a combination of left-hand sides of~\eqref{DT2}. It is
therefore $\hbar^2$ times a symmetric cubic in the~$x_a$, which has no component
of PBW degree one. Hence, at order~$\hbar^4$, the contribution of the left-hand
side of~\eqref{H4} to this bidegree is obtained from
\begin{align*}
\operatorname{pr}_{1,2}[Q_c,Q_d]
={}&h_c^{buv}h_d^{frs}[x_b,x_f]\ot K(x_u,x_v;x_r,x_s)\\
&+\overline h_c^{buv}\overline h_d^{frs}L(x_u,x_v;x_r,x_s)\ot q(x_b,x_f)\\
&+h_c^{buv}\overline h_d^{frs}M(x_b;x_r,x_s)\ot[q(x_u,x_v),x_f]\\
&+\overline h_c^{buv}h_d^{frs}M(x_f;x_u,x_v)\ot[x_b,q(x_r,x_s)].
\end{align*}

In the notation of Theorem~\ref{T:3}, the four maps have components
\begin{alignat*}{3}
& M(x_a;x_b,x_c)=m_{a;bc}^{i}x_i,\qquad&&
N(x_a,x_b;x_c,x_d)=n_{ab;cd}^{ij}q(x_i,x_j),& \\
& K(x_a,x_b;x_c,x_d)=k_{ab;cd}^{ij}q(x_i,x_j),\qquad&&
L(x_a,x_b;x_c,x_d)=\ell_{ab;cd}^{i}x_i.&
\end{alignat*}
The component expressions displayed there follow from the Jacobi consequence
\begin{gather*}
\tfrac1{12}\big([b,[c,d]]-2[c,[b,d]]\big)
=-\tfrac12M(d;b,c)+\tfrac18[[b,c],d],
\end{gather*}
and from the normalisation of the symmetric product recorded in the second
identity below.  We shall also use
\begin{gather*}
[q(x_u,x_v),x_f]
=\al_{uf}^{j}q(x_j,x_v)+\al_{vf}^{j}q(x_u,x_j),\qquad
q(x_b,x_f)=\tfrac12\big(\delta_b^j\delta_f^k+\delta_b^k\delta_f^j\big)q(x_j,x_k).
\end{gather*}
Substituting these component expressions into the preceding projection,
using the symmetry of the last two indices of $h$ and $\overline h$, and
contracting the Kronecker deltas gives
\begin{align*}
\operatorname{pr}_{1,2}[Q_c,Q_d]
={}&\Big[h_c^{buv}h_d^{frs}\al_{bf}^{i}k_{uv;rs}^{jk}
 +\big(
 \tfrac12\overline h_c^{juv}\overline h_d^{krs}\ell_{uv;rs}^{i}\\
&\hphantom{\Big[{}+\big({}}
 +h_c^{buk}\overline h_d^{frs}m_{b;rs}^{i}\al_{uf}^{j}
 -\overline h_c^{buv}h_d^{frk}m_{f;uv}^{i}\al_{rb}^{j}
 +(j\leftrightarrow k)\big)\Big]
 x_i\ot\{x_j,x_k\}\\
={}&\Pi_{cd}^{ijk}x_i\ot\{x_j,x_k\}.
\end{align*}

We next compute the same PBW component on the right-hand side of
\eqref{H4}.  The contribution of the $\hbar^2Q_k$ part of
$\LHb(\{x_i,x_j,\G(x_k)\})$ is
\begin{align*}
T_{ijk}={}&h_k^{abc}\Big(
x_a\ot N(x_i,x_j;x_b,x_c)
+\tfrac1{12}[x_i,x_a]\ot[x_j,q(x_b,x_c)]
+\tfrac1{12}[x_j,x_a]\ot[x_i,q(x_b,x_c)]\Big)\\
&\hphantom{h_k^{abc}\Big(}{} +\overline h_k^{abc}\Big(
M(x_i;x_b,x_c)\ot q(x_j,x_a)
+M(x_j;x_b,x_c)\ot q(x_i,x_a)\Big).
\end{align*}
Since $\Psi_{abc}^{pqr}=\Psi_{abc}^{qpr}$, the terms interchanged by
$p\leftrightarrow q$ combine.  Applying the same Jacobi and delta
contractions as above therefore yields
\begin{align*}
\begin{aligned}
\Psi_{abc}^{pqr}T_{pqr}
={}&\Big[\Psi_{abc}^{pqr}h_r^{iuv}n_{pq;uv}^{jk}\!
 +\big(\Psi_{abc}^{pjr}\overline h_r^{kuv}m_{p;uv}^{i}\!
-\tfrac16\Psi_{abc}^{pqr}h_r^{suk}\al_{ps}^{i}\al_{uq}^{j}
 +(j\leftrightarrow k)\big)\Big]
 x_i\ot\{x_j,x_k\}\\ 
={}&\Xi_{abc}^{ijk}x_i\ot\{x_j,x_k\}.
\end{aligned}
\end{align*}

The quintic term in \eqref{H4} makes no contribution to PBW
bidegree $(1,2)$.  Indeed, PBW symmetrisation is a coalgebra isomorphism, so
\begin{gather*}
\operatorname{pr}_{1,2}\LHb(\{x_i,x_j,x_k,x_l,x_m\})=0.
\end{gather*}
Equivalently, the expansion obtained from Lemma~\ref{L2} contains the
component
\begin{gather*}
\tfrac1{36}x_{(i}\ot[[x_j,x_{(k}],\{x_l,x_{m))}\}],
\end{gather*}
which cancels its companion term with coefficient $-\tfrac1{36}$.

It remains to compare the cubic terms.  The PBW bidegree-$(1,2)$
component of the coproduct of
$\Phi_{abc}^{ijk}\{x_i,x_j,x_k\}$ is
$3\Phi_{abc}^{ijk}x_i\ot\{x_j,x_k\}$.  Hence the equality of the two sides
of \eqref{H4} at order~$\hbar^4$ gives
\begin{align*}
3\Phi_{abc}^{ijk}x_i\ot\{x_j,x_k\}
=\al_{(ab}^{\qu d}\operatorname{pr}_{1,2}[Q_{c)},Q_d]
-\Psi_{abc}^{pqr}T_{pqr}=\big(\al_{(ab}^{\qu d}\Pi_{c)d}^{ijk}-\Xi_{abc}^{ijk}\big)
x_i\ot\{x_j,x_k\}.
\end{align*}
Only the completely symmetric part of the upper indices contributes to the
symmetric cubic product.  Applying $\operatorname{Sym}_{ijk}$ and dividing
by three gives exactly \eqref{coefevenloop3}, completing the determination
of the cubic coefficient $\Phi$, and with it of all three coefficients of~\eqref{H4}.

\subsection*{Acknowledgements}

The authors would like to thank P.~Baseilhac, N.~Cramp\'e, N.~Guay, N.~MacKay and J.~Ohayon for discussions and their interest in this work. The authors also thank the anonymous referees for their valuable remarks and suggestions. V.R.~acknowledges the UK EPSRC for the Postdoctoral Fellowship under grant EP/K031805/1. S.B.~is supported by a public grant as part of the ``Investissement d'avenir'' project, reference ANR-11-LABX-0056-LMH, LabEx LMH.

\pdfbookmark[1]{References}{ref}
\LastPageEnding


\begin{thebibliography}{99}
\footnotesize\itemsep=0pt

\bibitem{Ara62}
Araki S., On root systems and an infinitesimal classification of irreducible
 symmetric spaces, \textit{J.~Math. Osaka City Univ.} \textbf{13} (1962),
 1--34.

\bibitem{AACDFR04}
Arnaudon D., Avan J., Cramp\'e N., Doikou A., Frappat L., Ragoucy E., General
 boundary conditions for the $sl(N)$ and $sl(M|N)$ open spin chains,
 \href{https://doi.org/10.1088/1742-5468/2004/08/P08005}{\textit{J.~Stat. Mech. Theory Exp.}} \textbf{2004} (2004), P08005, 35~pages,
 \mbox{\href{http://arxiv.org/abs/math-ph/0406021}{math-ph/0406021}}.

\bibitem{Bas05a}
Baseilhac P., Deformed {D}olan--{G}rady relations in quantum integrable models,
 \href{https://doi.org/10.1016/j.nuclphysb.2004.12.016}{\textit{Nuclear Phys.~B}} \textbf{709} (2005), 491--521,
 \href{http://arxiv.org/abs/hep-th/0404149}{hep-th/0404149}.

\bibitem{Bas05b}
Baseilhac P., An integrable structure related with tridiagonal algebras,
 \href{https://doi.org/10.1016/j.nuclphysb.2004.11.014}{\textit{Nuclear Phys.~B}} \textbf{705} (2005), 605--619,
 \href{http://arxiv.org/abs/math-ph/0408025}{math-ph/0408025}.

\bibitem{BasBel10}
Baseilhac P., Belliard S., Generalized {$q$}-{O}nsager algebras and boundary
 affine {T}oda field theories, \href{https://doi.org/10.1007/s11005-010-0412-6}{\textit{Lett. Math. Phys.}} \textbf{93} (2010),
 213--228, \href{http://arxiv.org/abs/0906.1215}{arXiv:0906.1215}.

\bibitem{BasBel11}
Baseilhac P., Belliard S., The central extension of the reflection equations
 and an analog of {M}iki's formula, \href{https://doi.org/10.1088/1751-8113/44/41/415205}{\textit{J.~Phys.~A: Math. Theor.}}
 \textbf{44} (2011), 415205, 13~pages, \href{http://arxiv.org/abs/1104.1591}{arXiv:1104.1591}.

\bibitem{BasBel13}
Baseilhac P., Belliard S., The half-infinite {XXZ} chain in {O}nsager's
 approach, \href{https://doi.org/10.1016/j.nuclphysb.2013.05.003}{\textit{Nuclear Phys.~B}} \textbf{873} (2013), 550--584,
 \href{http://arxiv.org/abs/1211.6304}{arXiv:1211.6304}.

\bibitem{Bax72}
Baxter R.J., Partition function of the eight-vertex lattice model, \href{https://doi.org/10.1016/0003-4916(72)90335-1}{\textit{Ann.
 Physics}} \textbf{70} (1972), 193--228.

\bibitem{Bax82}
Baxter R.J., Exactly solved models in statistical mechanics, Academic Press,
 Inc., London, 1982.

\bibitem{BelCra12}
Belliard S., Cramp\'e N., Coideal algebras from twisted {M}anin triples,
 \href{https://doi.org/10.1016/j.geomphys.2012.05.008}{\textit{J.~Geom. Phys.}} \textbf{62} (2012), 2009--2023, \href{http://arxiv.org/abs/1202.2312}{arXiv:1202.2312}.

\bibitem{BelFom12}
Belliard S., Fomin V., Generalized {$q$}-{O}nsager algebras and dynamical
 {$K$}-matrices, \href{https://doi.org/10.1088/1751-8113/45/2/025201}{\textit{J.~Phys.~A: Math. Theor.}} \textbf{45} (2012), 025201,
 17~pages, \href{http://arxiv.org/abs/1106.1317}{arXiv:1106.1317}.

\bibitem{BFGP94}
Bonneau P., Flato M., Gerstenhaber M., Pinczon G., The hidden group structure
 of quantum groups: strong duality, rigidity and preferred deformations,
 \href{https://doi.org/10.1007/BF02099415}{\textit{Comm. Math. Phys.}} \textbf{161} (1994), 125--156.

\bibitem{ChaPre94}
Chari V., Pressley A., A~guide to quantum groups, Cambridge University Press,
 Cambridge, 1994.

\bibitem{GCM13}
Chen H., Guay N., Ma X., Twisted {Y}angians, twisted quantum loop algebras and
 affine {H}ecke algebras of type {$BC$}, \href{https://doi.org/10.1090/S0002-9947-2014-05994-1}{\textit{Trans. Amer. Math. Soc.}}
 \textbf{366} (2014), 2517--2574.

\bibitem{Che84}
Cherednik I.V., Factorizing particles on a half line, and root systems,
 \href{https://doi.org/10.1007/BF01038545}{\textit{Theoret. and Math. Phys.}} \textbf{61} (1984), 977--983.

\bibitem{DMS01}
Delius G.W., MacKay N.J., Short B.J., Boundary remnant of {Y}angian symmetry
 and the structure of rational reflection matrices, \href{https://doi.org/10.1016/S0370-2693(01)01275-8}{\textit{Phys. Lett.~B}}
 \textbf{522} (2001), 335--344, \href{http://arxiv.org/abs/hep-th/0109115}{hep-th/0109115}.

\bibitem{Dri85}
Drinfeld V.G., Hopf algebras and the quantum {Y}ang--{B}axter equation,
 \textit{Dokl. Akad. Nauk SSSR} \textbf{32} (1985), 254--258.

\bibitem{Dri87}
Drinfeld V.G., Quantum groups, in Proceedings of the {I}nternational {C}ongress
 of {M}athematicians, {V}ols.~1,~2 ({B}erkeley, {C}alif., 1986), Amer. Math.
 Soc., Providence, RI, 1987, 798--820.

\bibitem{Dri92}
Drinfeld V.G., On some unsolved problems in quantum group theory, in
 Quantum Groups ({L}eningrad, 1990), \href{https://doi.org/10.1007/BFb0101175}{\textit{Lecture Notes in Math.}}, Vol.~1510, Springer, Berlin, 1992, 1--8.

\bibitem{EtiKaz96}
Etingof P., Kazhdan D., Quantization of {L}ie bialgebras.~{I}, \href{https://doi.org/10.1007/BF01587938}{\textit{Selecta
 Math.~(N.S.)}} \textbf{2} (1996), 1--41, \href{http://arxiv.org/abs/q-alg/9506005}{q-alg/9506005}.

\bibitem{EtiKaz95}
Etingof P., Kazhdan D., Quantization of {P}oisson algebraic groups and
 {P}oisson homogeneous spaces, in Sym\'etries Quantiques ({L}es {H}ouches,
 1995), North-Holland, Amsterdam, 1998, 935--946, \href{http://arxiv.org/abs/q-alg/9510020}{q-alg/9510020}.

\bibitem{FadRes86}
Faddeev L.D., Reshetikhin N.Yu., Integrability of the principal chiral field
 model in {$1+1$} dimension, \href{https://doi.org/10.1016/0003-4916(86)90201-0}{\textit{Ann. Physics}} \textbf{167} (1986),
 227--256.

\bibitem{FadTak81}
Faddeev L.D., Takhtadzhyan L.A., Spectrum and scattering of excitations in the
 one-dimensional isotropic Heisenberg model, \href{https://doi.org/10.1007/BF01087245}{\textit{J.~Sov. Math.}}
 \textbf{24} (1984), 241--267.

\bibitem{GR}
Guay N., Regelskis V., Twisted {Y}angians for symmetric pairs of types {B},
 {C}, {D}, \href{https://doi.org/10.1007/s00209-016-1649-2}{\textit{Math.~Z.}} \textbf{284} (2016), 131--166,
 \href{http://arxiv.org/abs/1407.5247}{arXiv:1407.5247}.

\bibitem{GRW1}
Guay N., Regelskis V., Wendlandt C., Twisted {Y}angians of small rank,
 \href{https://doi.org/10.1063/1.4947112}{\textit{J.~Math. Phys.}} \textbf{57} (2016), 041703, 28~pages,
 \href{http://arxiv.org/abs/1602.01418}{arXiv:1602.01418}.

\bibitem{GRW2}
Guay N., Regelskis V., Wendlandt C., Representations of twisted Yangians of
 types~B,~C,~D:~I, \href{https://doi.org/10.1007/s00029-017-0306-x}{\textit{Selecta Math.~(N.S.)}}, {t}o appear,
 \href{http://arxiv.org/abs/1605.06733}{arXiv:1605.06733}.

\bibitem{Hel78}
Helgason S., Differential geometry, {L}ie groups, and symmetric spaces,
 \textit{Pure and Applied Mathematics}, Vol.~80, Academic Press, Inc., New
 York~-- London, 1978.

\bibitem{Kol12}
Kolb S., Quantum symmetric {K}ac--{M}oody pairs, \href{https://doi.org/10.1016/j.aim.2014.08.010}{\textit{Adv. Math.}}
 \textbf{267} (2014), 395--469, \href{http://arxiv.org/abs/1207.6036}{arXiv:1207.6036}.

\bibitem{Koo93}
Koornwinder T.H., Askey--{W}ilson polynomials as zonal spherical functions on
 the {${\rm SU}(2)$} quantum group, \href{https://doi.org/10.1137/0524049}{\textit{SIAM~J. Math. Anal.}} \textbf{24}
 (1993), 795--813.

\bibitem{Let02}
Letzter G., Coideal subalgebras and quantum symmetric pairs, in New Directions
 in {H}opf Algebras, \textit{Math. Sci. Res. Inst. Publ.}, Vol.~43, Cambridge
 University Press, Cambridge, 2002, 117--165, \href{http://arxiv.org/abs/math.QA/0103228}{math.QA/0103228}.

\bibitem{Lev93}
Levendorski\u{\i} S.Z., On {PBW} bases for {Y}angians, \href{https://doi.org/10.1007/BF00739587}{\textit{Lett. Math.
 Phys.}} \textbf{27} (1993), 37--42.

\bibitem{Mac02}
MacKay N.J., Rational {$K$}-matrices and representations of twisted {Y}angians,
 \href{https://doi.org/10.1088/0305-4470/35/37/302}{\textit{J.~Phys.~A: Math. Gen.}} \textbf{35} (2002), 7865--7876,
 \href{http://arxiv.org/abs/math.QA/0205155}{math.QA/0205155}.

\bibitem{Mac04}
MacKay N.J., Introduction to {Y}angian symmetry in integrable field theory,
 \href{https://doi.org/10.1142/S0217751X05022317}{\textit{Internat.~J. Modern Phys.~A}} \textbf{20} (2005), 7189--7217,
 \href{http://arxiv.org/abs/hep-th/0409183}{hep-th/0409183}.

\bibitem{MacReg10}
MacKay N., Regelskis V., Yangian symmetry of the {$Y=0$} maximal giant
 graviton, \href{https://doi.org/10.1007/JHEP12(2010)076}{\textit{J.~High Energy Phys.}} \textbf{2010} (2010), no.~12, 076,
 17~pages, \href{http://arxiv.org/abs/1010.3761}{arXiv:1010.3761}.

\bibitem{MacReg11}
MacKay N., Regelskis V., Achiral boundaries and the twisted {Y}angian of the
 {D}5-brane, \href{https://doi.org/10.1007/JHEP08(2011)019}{\textit{J.~High Energy Phys.}} \textbf{2011} (2011), no.~8, 019,
 22~pages, \href{http://arxiv.org/abs/1105.4128}{arXiv:1105.4128}.


\bibitem{Mo92}
Molev A.I., Representations of twisted {Y}angians, \href{https://doi.org/10.1007/BF00420754}{\textit{Lett. Math. Phys.}}
 \textbf{26} (1992), 211--218.

\bibitem{Mo97}
Molev A.I., Finite-dimensional irreducible representations of twisted
 {Y}angians, \href{https://doi.org/10.1063/1.532551}{\textit{J.~Math. Phys.}} \textbf{39} (1998), 5559--5600,
 \href{http://arxiv.org/abs/q-alg/9711022}{q-alg/9711022}.

\bibitem{MNO96}
Molev A.I., Nazarov M.L., Olshanskii G.I., Yangians and classical {L}ie
 algebras, \href{https://doi.org/10.1070/RM1996v051n02ABEH002772}{\textit{Russ. Math. Surv.}} \textbf{51} (1996), 205--282,
 \href{http://arxiv.org/abs/hep-th/9409025}{hep-th/9409025}.

\bibitem{MolRag02}
Molev A.I., Ragoucy E., Representations of reflection algebras, \href{https://doi.org/10.1142/S0129055X02001156}{\textit{Rev.
 Math. Phys.}} \textbf{14} (2002), 317--342, \href{http://arxiv.org/abs/math.QA/0107213}{math.QA/0107213}.

\bibitem{MRS03}
Molev A.I., Ragoucy E., Sorba P., Coideal subalgebras in quantum affine
 algebras, \href{https://doi.org/10.1142/S0129055X03001813}{\textit{Rev. Math. Phys.}} \textbf{15} (2003), 789--822,
 \href{http://arxiv.org/abs/math.QA/0208140}{math.QA/0208140}.

\bibitem{Oha10}
Ohayon J., Quantization of coisotropic subalgebras in complex semisimple {L}ie
 algebras, \href{http://arxiv.org/abs/1005.1371}{arXiv:1005.1371}.

\bibitem{Ols90}
Olshanskii G.I., Twisted {Y}angians and infinite-dimensional classical {L}ie
algebras, in Quantum Groups ({L}eningrad, 1990), \href{https://doi.org/10.1007/BFb0101183}{\textit{Lecture Notes in Math.}}, Vol.~1510, Springer, Berlin, 1992, 104--119.

\bibitem{Reg12}
Regelskis V., Reflection algebras for ${\mathfrak{sl}}(2)$ and
 ${\mathfrak{gl}}(1|1)$, \href{http://arxiv.org/abs/1206.6498}{arXiv:1206.6498}.

\bibitem{FRT90}
Reshetikhin N.Yu., Takhtadzhyan L.A., Faddeev L.D., Quantization of {L}ie groups
 and {L}ie algebras, \textit{Leningrad Math.~J.} \textbf{1} (1990), 193--225.

\bibitem{Skl88}
Sklyanin E.K., Boundary conditions for integrable quantum systems,
 \href{https://doi.org/10.1088/0305-4470/21/10/015}{\textit{J.~Phys.~A: Math. Gen.}} \textbf{21} (1988), 2375--2389.

\bibitem{FST79}
Sklyanin E.K., Takhtadzhyan L.A., Faddeev L.D., Quantum inverse problem
 method.~{I}, \href{https://doi.org/10.1007/BF01018718}{\textit{Theoret. and Math. Phys.}} \textbf{40} (1979), 688--706.

\bibitem{Yan67}
Yang C.N., Some exact results for the many-body problem in one dimension with
 repulsive delta-function interaction, \href{https://doi.org/10.1103/PhysRevLett.19.1312}{\textit{Phys. Rev. Lett.}} \textbf{19}
 (1967), 1312--1315.

\bibitem{Zam11}
Zambon M., A construction for coisotropic subalgebras of {L}ie bialgebras,
 \href{https://doi.org/10.1016/j.jpaa.2010.04.026}{\textit{J.~Pure Appl. Algebra}} \textbf{215} (2011), 411--419,
 \href{http://arxiv.org/abs/0810.5160}{arXiv:0810.5160}.

\end{thebibliography}
\end{document}